\documentclass[10pt]{article}

\usepackage{graphicx}
\usepackage{dcolumn}
\usepackage{amsmath}
\usepackage{xcolor}
\usepackage{authblk}
\usepackage[letterpaper, margin=1in]{geometry}

\begin{document}


\title{Steiner triangular drop dynamics} 



\author[1]{Elizabeth Wesson}
\author[2]{Paul Steen}
\affil[1]{Center for Applied Mathematics, Cornell University}
\affil[2]{School of Chemical and Biomolecular Engineering, Cornell University}

\date{\today}

\maketitle 

\begin{abstract}
Steiner's circumellipse is the unique geometric regularization of any triangle to a circumscribed ellipse with the same centroid, a regularization that motivates our introduction of the Steiner triangle as a minimal model for liquid droplet dynamics.  The Steiner drop is a deforming triangle with one side making sliding contact against a  planar basal support.  The center of mass of the triangle is governed by Newton's law.  The resulting dynamical system  lives in a four dimensional phase space and exhibits a rich one-parameter family of dynamics.  Two invariant manifolds are identified with ``bouncing" and ``rocking" periodic motions; these intersect at the stable equilibrium and are surrounded by nested quasiperiodic motions. We study the inherently interesting dynamics and also find that this model, however minimal, can capture space-time symmetries of more realistic continuum drop models.
\end{abstract}

\textbf{Polygons are the simplest of planar shapes and triangles the simplest of all polygons.  Here, we introduce a triangle with variable base length and height as a minimal model of a sliding and deforming liquid droplet, inspired by MJ Steiner's (French geometer) classic theorem relating triangles to ellipses.  The Steiner drop's shape is related to Newtonian dynamics through a constitutive relationship that penalizes against sharp angles. The resulting dynamical system with a four dimensional phase space exhibits fixed points, periodic and quasiperiodic orbits, all related to two competing invariant manifolds.  The symmetries of the dynamical system allow us to characterize the Steiner drop motions near equilibrium as ``bouncing,'' ``rocking,'' or quasiperiodic combinations of the two. By associating each triangle with its unique smallest circumscribed ellipse, the motions of the Steiner drops can be considered equivalently as motions of elliptical drops. To relate to more realistic drops, we also compute the center-of-mass motions of normal modes of perturbed spherical cap liquid drops and find that they also can be classified as bouncing or rocking. Thus the minimal model predicts the qualitative behavior of spherical-cap drop oscillations.}

\section{Introduction}  We introduce the Steiner triangular model of a sliding and deforming liquid drop shaped by surface tension. In nature, of course, surface tension does not allow for corners.  This deficiency is overlooked in favor of the simplicity of the triangle, whose deformations have just two degrees of freedom.  The base can lengthen and shorten, consistent with a moving contact line, and the apex can move about in the plane. 
We seek a reduced order (ODE) model that captures  the space-time dynamics of higher order (PDE) continuum descriptions.  Our model is termed \emph{minimal} since, on removing one degree of freedom, the model becomes trivial dynamically. 

The triangular geometry can be regularized to an ellipse.  Every triangle is uniquely associated with a circumscribed ellipse whose center coincides with the centroid of the triangle and whose area is minimal, the so-called `Steiner ellipse' \cite{steiner1829,weisstein}.  Conversely, every ellipse has infinitely many maximal inscribed triangles, but only two of these have a side parallel to a given line, and these are reflections of each other. Thus each ellipse can be uniquely associated with a maximal inscribed triangle whose base is horizontal, which we refer to as its `Steiner triangle' (not to be confused with the usage of the term `Steiner triangle' as the Civian triangle of a Steiner point). This one-to-one mapping between Steiner triangles and ellipses yields an equivalence between triangular and elliptical drops.  Motions can be equivalently viewed as those of triangles or of ellipses.  

To complete the definition of the Steiner triangular drop, one needs to introduce a constitutive relationship between force and deformation that mimics surface tension.  While we accept corners as natural to triangles, we do assess a penalty on the sharpness of the corners. Sharply acute angles are penalized by a pressure that blows up in the zero angle limit. This can be thought of as the analog of the Young-Laplace relationship, familiar from the continuum model for normal stress balance across an interface. For a minimal model, the exact nature of the force-deformation behavior is not important as long as there is a robust penalty for flattening the triangle.

Our interest is in the dynamical regime where liquid inertia is shaped by surface tension and dissipation can be ignored, the so-called capillary ballistic regime where the dynamics are conservative. Sources of dissipation that are neglected are bulk dissipation due to the inherent viscosity of the liquid and possibly contact line dissipation due to drag of the moving contact line \cite{davis,deGennes,bostwick2015stability}.  Experiments with water drops on substrates fall into this regime \cite{sharpEtAl,sharp,chang2013substrate}, and PDE models have been successful at predicting observed frequencies of vibration \cite{basaran,lyubimovNonAxi,lyubimovBehavior,fayzra}.

For conservative dynamics, the motion of the drop's center of mass decouples from its deformation relative to the inertial frame.  In this paper, we restrict our study to motions in the inertial frame, consistent with a focus on a drop's natural vibrations.  This eliminates the translational degree-of-freedom introduced by the sliding motion. Motions relative to the center of mass are solved for and  motion of the contact points, $A$ and $B$ in Figure~\ref{fig:defn-sketch}(a), relative to the substrate are preserved up to a (Galilean) translational velocity.

The relevant configurational variables are the coordinates of the triangle's center of mass relative to its base. Newton's law is applied to each of these, which leads to a first-order dynamical system of nonlinear ODEs of dimension four.  Having a four dimensional phase-space means that, {\it a priori}, our minimal model  can possibly exhibit long-term dynamics that includes fixed points, periodic, and quasi-periodic orbits.  Surprisingly, this richness is actually realized.

Liquid drops that deform and translate along a solid support are found widely in nature and in application \cite{kumar2015liquid,josserand2016drop}.  The contact line represents the locus where liquid, solid and gas phases meet and is characteristic of partially wetting liquids such as water on a silicon wafer.  Drop translation requires contact lines to move, and rapidly moving contact lines remain a modeling challenge in the continuum description \cite{snoeijer2013moving,sui2014numerical}.  In our $2D$ model, contact lines become contact points, the points $A$ and $B$, Figure~\ref{fig:defn-sketch}(a). These points move independently according to the capillary forces acting on them and the Newtonian inertial dynamics of the center of mass. 

From an application viewpoint, there is a growing need to simulate large populations of droplets that are translating, deforming and colliding, splitting up (atomizing) or perhaps joining up (coalescing) \cite{tanaka,macner2014condensation,xu2017collective}. Having a minimal model that captures contact line motions with some fidelity provides a computationally economical path to studying the dynamics of large populations.

Vibrated drops of water exhibit capillary ballistic motions associated with a series of mode shapes and resonant frequencies \cite{chang2015dynamics}. The shapes and frequencies, or spectrum, of the free liquid sphere modeled as a continuum, were first reported by Rayleigh \cite{rayleigh18}. By symmetry, vibrations the hemispherical drop maintaining a $\pi/2$ contact angle against a planar substrate are  part of Rayleigh's family of solutions. Studies of the dynamical response of the hemispherical sessile drop subject to other contacting conditions are more recent \cite{noblin,sharp,xia}. Analyses have been extended to spherical cap drops both with pinned and moving contact lines\cite{bostwick2015stability}.  The space-time symmetries represented are non-trivial.  In this paper using the center-of-mass metric, we show  that all these responses map into one of three classes of motions: steady, bouncing or rocking. Our Steiner drops faithfully predict these classes.

This paper is organized as follows. In Sections  \ref{sec:model}-\ref{sec:CM} we introduce the model, derive the equations of motion, and convert to center-of-mass coordinates. The resulting dynamical system and its symmetries are described in Section \ref{sec:dynsys}. In Sections \ref{sec:eqpts}-\ref{sec:stability} we identify the equilibria of the system and analyze their stability and bifurcation. In Section \ref{sec:sym} we use the symmetries of the system to identify two-dimensional invariant manifolds filled with ``bouncing'' and ``rocking'' periodic orbits. Section \ref{sec:torus} describes the torus trajectories in the region surrounding the stable equilibrium. Finally, in Section \ref{sec:fluid}, we show that for small perturbations, the center-of-mass motion of sessile spherical-cap drops can be characterized as ``bouncing,'' ``rocking,'' or a combination of those. Thus the triangular drop model qualitatively predicts motions of fluid drops.

\section{The model}\label{sec:model}
\begin{figure} 
\centering
 \includegraphics[width=0.45\textwidth]{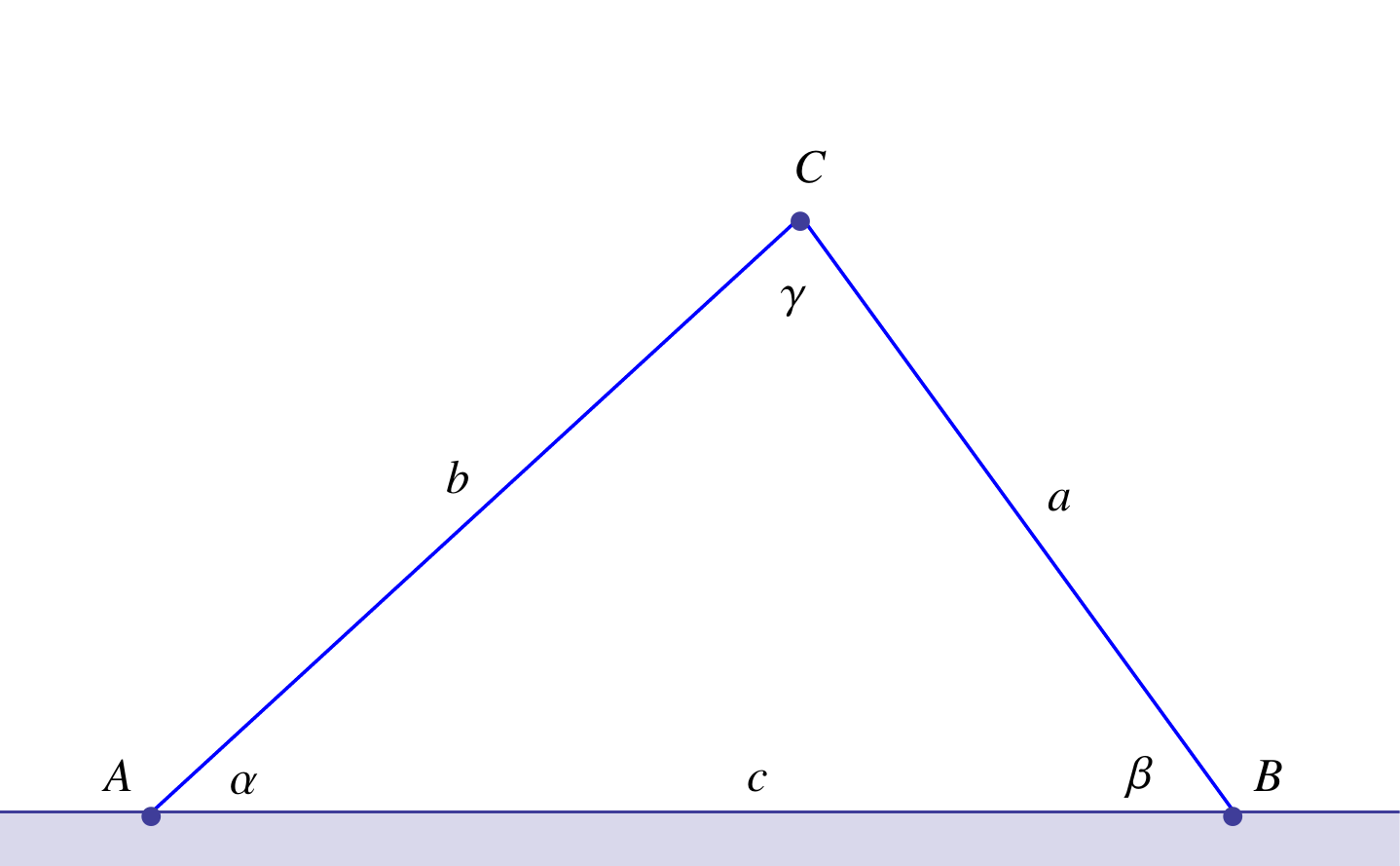}
 \includegraphics[width=0.45\textwidth]{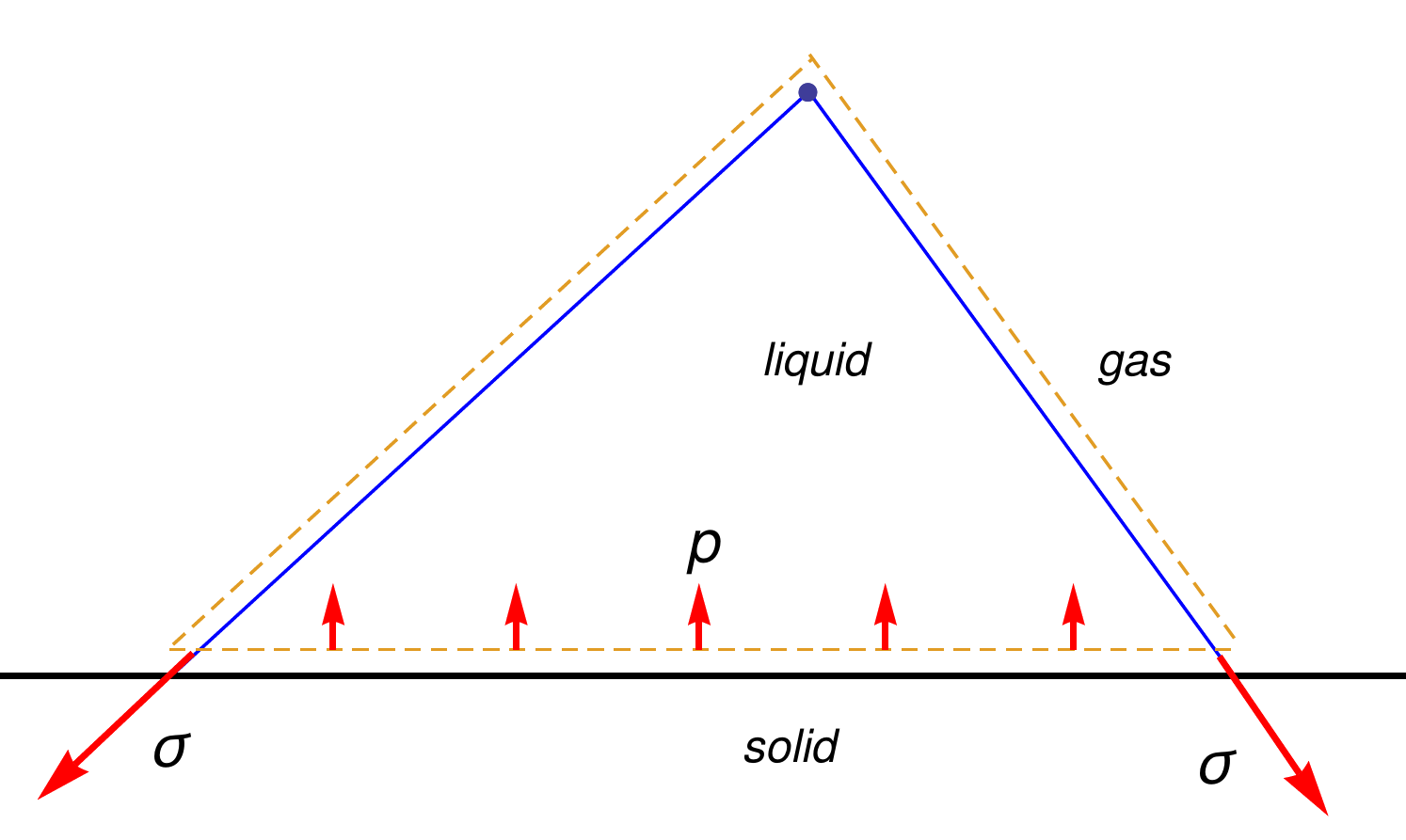}
 \caption{\textit{(a) Schematic of the triangular drop on the substrate. (b) Free-body diagram.}}
 \label{fig:defn-sketch}
\end{figure}

The triangle $ABC$  sits with its base on the $x$ axis which models the planar support. The triangle has vertices $A,B,C$, with corresponding interior angles $\alpha,\beta,\gamma$ and opposite edges $a,b,c$, respectively.  Side lengths are denoted by $a,b,c$: Figure~\ref{fig:defn-sketch} (a). We shall speak of a constant volume $V$ and a length scale $\ell \equiv V^{1/3}$ even though, for our Steiner triangles, it is an area $\ell^2$ which is preserved to model liquid incompressibility.  Alternatively, one can think of the Steiner drops as prismatic with a uniform depth of dimension $\ell$. In either case, one ends up with a description of a deforming $2D$ object with physical units that correspond to a $3D$ drop, for convenience in formulation, and without loss of generality.

The instantaneous configuration of the triangle is determined by two angles and the length scale $\ell$.  In order to apply Newton's law, the instantaneous forces acting on the object must be identified.  

This is done using the free-body diagram method, Figure \ref{fig:defn-sketch} (b), where we imagine a planar section just above the support which cuts the interface and is `shrink-wrapped' around the rest of the triangle.  The net force on the triangle thus has two contributions, one from the surface tension $\sigma$, acting tangent to the liquid/gas interface at the contact line, $\mathbf{F_{surf}}$, and one from the pressure $p$ acting plane-normal to the support, $\mathbf{F_{pres}}$.  The surface tension contribution is of constant magnitude, with varying components determined by the contact angle, $\alpha$ or $\beta$ via the Young-Dupre equation at the points $A$ and $B$.  The pressure contribution is also proportional to surface tension but depends on the global configuration in a way that needs to be modeled, 
\begin{equation}
    p=\sigma \kappa(\alpha,\beta;\ell),
    \label{eq:press}
\end{equation}
as described below.  These contributions are analogous to those exerted on a small supported drop in a full three dimensional continuum description \cite{bostwick2015stability}. Thus we have

\begin{align}
\mathbf{F_{surf}} & = \sigma \ell ([\cos\beta,-\sin\beta] + [-\cos\alpha,-\sin\alpha]) \label{eq:fsurf}\\
\mathbf{F_{pres}} & = \sigma \ell \kappa \ [0,c].
\end{align}

The net force gives the acceleration relative to the inertial frame of the center of mass of the triangle, $\mathbf{x}$, according to Newton's law: 
\begin{equation}
\rho V \frac{d^2\mathbf{x}}{dt^2} = \mathbf{F} = \mathbf{F_{surf}} + \mathbf{F_{pres}}.
\label{eq:newton1}
\end{equation}
Inertia, proportional to $\rho V$, and capillarity, proportional to $\sigma$, compete to determine the dynamics.  Note that the support is fixed in the lab frame and that, relative to the lab frame, the net force is non-zero, in both directions.  In contrast, the inertial frame is attached to the droplet center of mass, a frame in which the net force is zero and with respect to which there is no relative motion.

Our guide to pressure closure is the Laplace pressure across a curved surface, which penalizes against small radii of curvature.  Analogously, we choose a functional form for the curvature $\kappa$, and hence pressure \eqref{eq:press}, that penalizes against small internal angles:

\begin{equation}
 \kappa= k\left( \frac{1}{\sin\alpha} +\frac{1}{\sin\beta} + \frac{1}{\sin\gamma} \right).
\label{eq:pres}
\end{equation}

Here $k$ is constant, determined as now described.  For the fully $3D$ droplet, equilibrium requires that the integral of the horizontal force component around the closed curve representing the contact line vanishes.  The corresponding condition here is that the sum of horizontal force components at vertices $A$ and $B$ be zero. This implies that, at equilibrium, $\alpha = \beta \equiv \alpha_0$ and $c \equiv c_0$, the length of side $c$. The equilibrium contact angle $\alpha_0$ is readily measured in real drops, and makes a good point of comparison for the model.

\begin{figure}
\centering
 \includegraphics[width=0.45\textwidth]{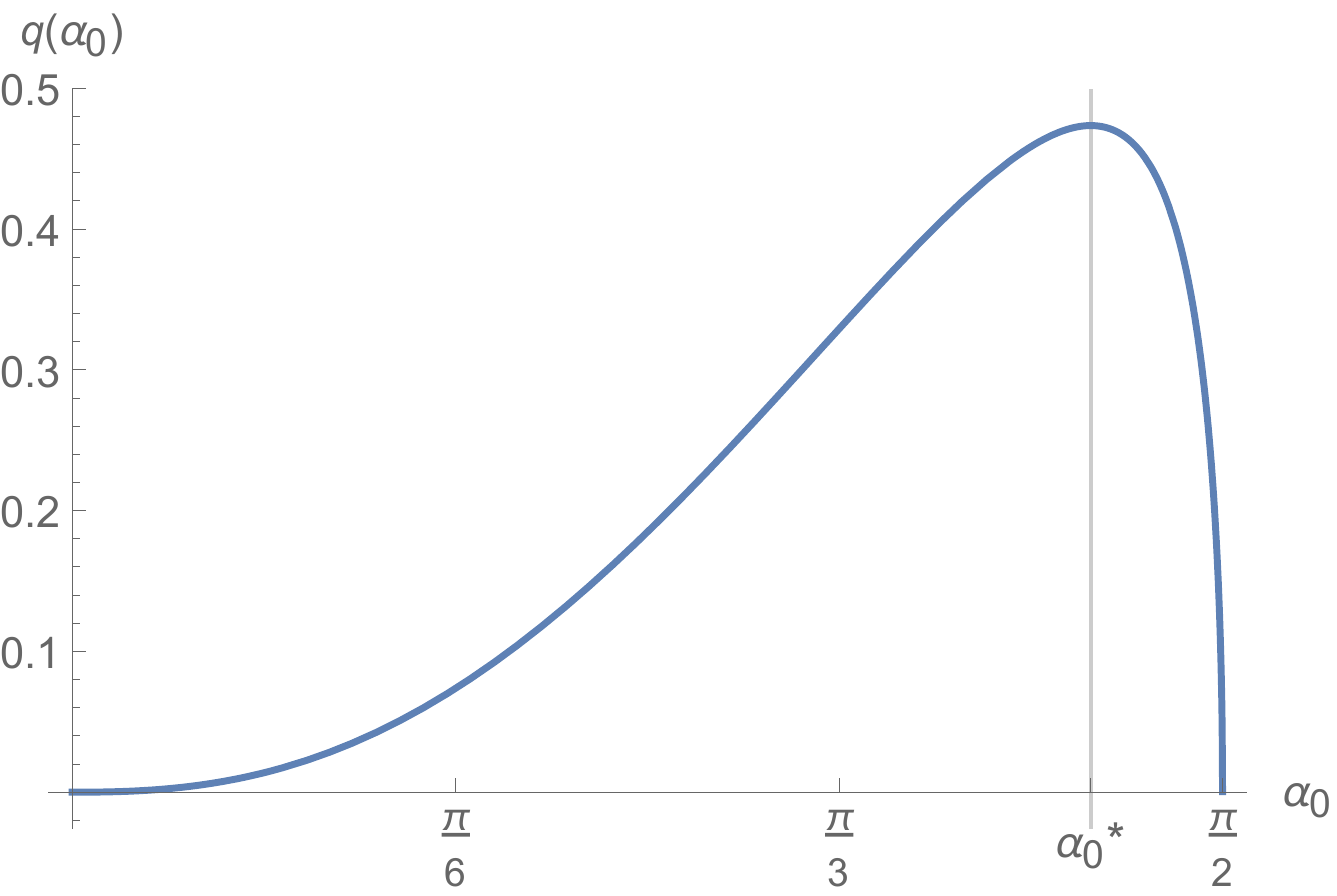}
 \caption{\textit{The dimensionless pressure coefficient $q(\alpha_0)$ as a function of $\alpha_0$.}}
 \label{fig:k(alpha0)}
\end{figure}

All equilibria, $(\alpha, \beta, c) = (\alpha_0, \alpha_0, c_0)$ are evidently isosceles triangles.   Using the identity $\alpha + \beta + \gamma = \pi$ to eliminate $\gamma$ in \eqref{eq:pres} and equating the vertical force component to zero in \eqref{eq:newton1} as required by equilibrium, one finds,

\begin{equation}
 k  =\frac{4 \sin ^2 \alpha _0}{c_0 (4 + \sec \alpha _0)} = \frac{2 \sin ^2 \alpha _0\sqrt{\tan{\alpha _0}}}{\ell (4 + \sec \alpha _0)} \label{eq:kval}
\end{equation}
where the second expression uses $c_0 = 2\ell \sqrt{\cot \alpha_0 }$, by the area constraint. Thus we obtain a dimensionless parameter
\begin{equation}\label{eq:qdef}
    q(\alpha_0) \equiv \ell k = \frac{2 \sin ^2 \alpha _0\sqrt{\tan{\alpha _0}}}{ 4 + \sec \alpha _0}.
\end{equation}
which is shown in Figure \ref{fig:k(alpha0)}. Notice that each value of $q$ is attained at two values of $\alpha_0$, except its maximum value, which we denote by $q(\alpha_0^*)$.

\section{Center of mass description}\label{sec:CM}
To make Newton's law into a dynamical system, we express the net force as a function of $(x,y)$, the Cartesian coordinates of the triangle's center of mass.
Let the vertices have Cartesian coordinates $(x_A,0)$, $(x_B,0)$, and $(x_C,y_C)$ respectively. Then we have
\begin{align}\label{eq:chvar1}
 \cos\alpha=(x_C-x_A)/b,& \quad \cos\beta=(x_B-x_C)/a, \nonumber\\
 \quad \sin\alpha=y_C/b, & \quad \sin\beta=y_C/a.
\end{align}
The origin is put at the center of the base, stationary in the lab frame. Let $d$ be the half-width of the base, so $x_A=-d$ and $x_B=d$. Trigonometry and volume constraint lead to  
\begin{equation}
 (x_A,x_B,x_C,y_C) = \left(-\frac{\ell^2}{3y}, \frac{\ell^2}{3y},3x, 3y\right),
\end{equation}
and the lengths of the three sides in terms of $x$ and $y$,
\begin{align}
 a^2 & =\left(\frac{\ell^2}{3y}-3x\right)^2+(3y)^2,\\ 
 b^2 & =\left(\frac{\ell^2}{3y}+3x\right)^2+(3y)^2,\\
 c^2 & =\left(\frac{2\ell^2}{3y}\right)^2.\label{eq:chvar2}
\end{align}
The components of the net force, $\mathbf{F}=(F_x,F_y)$, can now be written:
\begin{align}
 \frac{F_x}{\sigma \ell} & = \frac{1}{a}\left(\frac{\ell^2}{3y}-3x \right) - \frac{1}{b}\left(\frac{\ell^2}{3y}+3x \right) \label{eq:chvar3}\\ 
\frac{F_y}{\sigma \ell} & = -3y\left(\frac{1}{a}+\frac{1}{b}\right)
 + \frac{k \ell^2}{3y}\left(\frac{2(a+b)}{3y} + \frac{ab}{\ell^2}\right).\label{eq:chvar4}
\end{align}
A dimensionless form of Newton's law \eqref{eq:newton1} follows from \eqref{eq:chvar3} and \eqref{eq:chvar4} by scaling lengths by $\ell$,  and time by the inertial time, $t= \hat{t} \sqrt{\rho V/\sigma}$.

\section{Dynamical System}\label{sec:dynsys}
Rewriting \eqref{eq:newton1} using \eqref{eq:chvar3} and \eqref{eq:chvar4} yields a dimensionless second-order ODE for $\hat{x},\hat{y}$ in $\hat{t}$.  We drop all hats and let dots represent time derivatives to record the resulting dynamical system,
\begin{align}
 [\ddot x,\ddot y] &= \left[\frac{1}{a}\left(\frac{1}{3y}-3x \right) - \frac{1}{b}\left(\frac{1}{3y}+3x \right),\right. \nonumber \\
 &\left.\quad-3y\left(\frac{1}{a}+\frac{1}{b}\right)
 + \frac{q(\alpha_0)}{3y}\left(\frac{2(a+b)}{3y} + ab\right)\right] \label{eq:dynsys}\\
 &\equiv [f(x,y),h(x,y)] 
\end{align}
where the non-dimensionalized side lengths are
\begin{align}
 a &= \sqrt{\left(\frac{1}{3 y}-3 x\right)^2+9 y^2} \\
 b &= \sqrt{\left(3 x+\frac{1}{3 y}\right)^2+9 y^2}
\end{align}
and  $q(\alpha_0)$ is defined in Equation \eqref{eq:qdef}.

It is useful to convert the two second-order equations to four first-order equations:
\begin{align}
 \dot x &= w \label{eq:sys1}\\
 \dot w &= f(x,y) \\
 \dot y &= z \\
 \dot z &= h(x,y). \label{eq:sys4}
\end{align}
The system is not Hamiltonian using the natural generalized coordinates $\mathbf{q}=(x,y)$, $\mathbf{p}=(w,z)$, since 
\begin{equation}
    \frac{\partial f}{\partial y} \neq \frac{\partial h}{\partial x}
\end{equation}

However, notice that $f(x,y)$ is odd and $h(x,y)$ is even with respect to $x$:
\begin{equation}
 f(-x,y)=-f(x,y),\quad h(-x,y)=h(x,y).
\end{equation}

This means that the system is reversible, and in particular it is invariant \cite{ROBERTS1992} under the time-reversing phase space involutions
\begin{equation} \label{eq:invo1}
G_1:\{t\mapsto-t,\,(x,w,y,z)\mapsto(-x,w,y,-z)\}    
\end{equation}
and
\begin{equation}\label{eq:invo2}
    G_2:\{t\mapsto-t,\,(x,w,y,z)\mapsto(x,-w,y,-z)\}
\end{equation}
and hence also under the symmetry \cite{lamb}
\begin{equation} \label{eq:sym}
    S=G_2\circ G_1: (x,w,y,z)\mapsto(-x,-w,y,z)
\end{equation}

Although Hamiltonian systems need not be reversible or equivariant, and vice versa, \cite{lamb} these classes of  systems share many dynamical properties. In particular, there are reversible analogues of the Lyapunov center theorem \cite{devaney} and KAM theory. \cite{sevryuk91,SEVRYUK1998}

 These results allow us to characterize the orbits of our system near the neutrally stable equilibrium, discussed below. Our system is time-reversible and hence Lyapunov stability is the relevant definition; below we shall simply refer to Lyapunov stable states, as `stable'.

\section{Equilibria}\label{sec:eqpts}
In analyzing the equilibria of the drop, it is easiest to write the components of the net force in terms of the contact angles. By the symmetry of the system, any equilibrium configuration must be an isosceles triangle, i.e. $\alpha=\beta$. By trigonometry, the scaled (dimensionless) center of mass of such a triangle is
\begin{equation}\label{eq:equil-cond}
    (x,y)  = \left(0,\frac{\sqrt{\tan\alpha}}{3}\right).
\end{equation}
Substituting these values into Equation \eqref{eq:dynsys} gives
\begin{align}
    \ddot x &= 0 \\
    \ddot y &= -2\sin\alpha + q(\alpha_0) \frac{4+\sec\alpha}{\sin\alpha\sqrt{\tan\alpha}} \nonumber \\
    &= -2\sin\alpha\left(1-\frac{q(\alpha_0)}{q(\alpha)}\right).\label{eq:vert}
\end{align}
Thus any equilibrium configuration is an isosceles triangle with contact angle $\alpha$, where
\begin{equation}
    q(\alpha)=q(\alpha_0). \label{eq:equil-q}
\end{equation} 
Since each value of $q(\alpha_0)$, except the maximum, is attained at two values of $\alpha_0$, this means that in general there are two equilibria.

The equilibrium with contact angle $\alpha_0$ is given by
\begin{equation}
 \mathbf{x}_0=(x_0,y_0)=\left(0,\frac{\sqrt{\tan\alpha_0}}{3}\right).
\end{equation}
We will denote the other equilibrium by $\mathbf{x}_1=(0,y_1)$. Using Equations \eqref{eq:qdef} and \eqref{eq:equil-cond} in Equation \eqref{eq:equil-q},
we see that $y_1$ is the second real, positive solution of
\begin{equation}\label{equil2}
 \frac{486 y^5}{\left(81 y^4+1\right) \left(\sqrt{81 y^4+1}+4\right)} = q( \alpha_0)
\end{equation}
There is no closed-form expression for $y_1$, but we can find it numerically for a given value of $\alpha_0$. 

\begin{figure}
\centering
\includegraphics[width=0.4\textwidth]{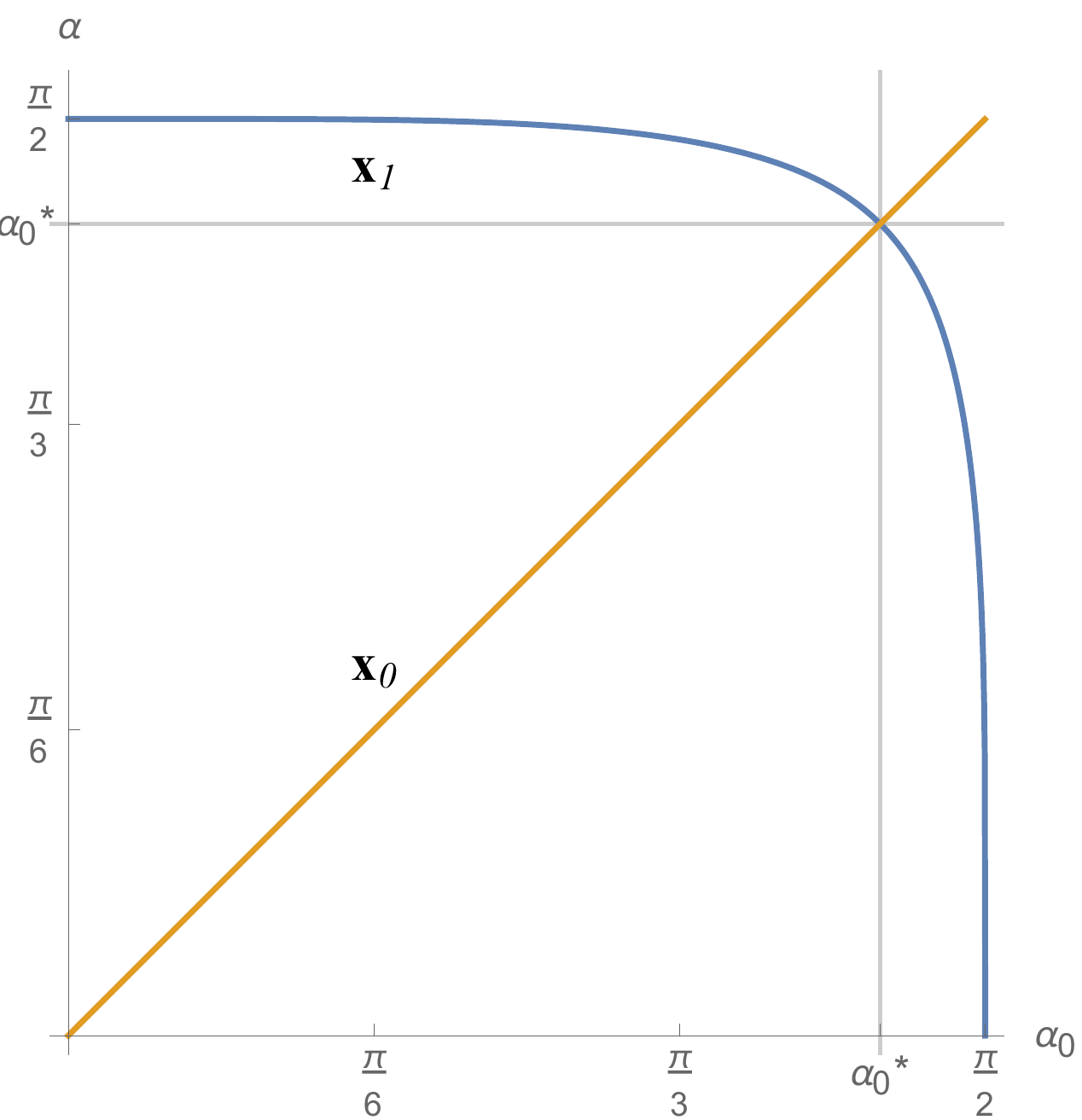}
\caption{\textit{Contact angles $\alpha$ of the two equilibria $\mathbf{x}_0$ (yellow) and $\mathbf{x}_1$ (blue), as functions of $\alpha_0$.}}
\label{fnet(alpha)}
\end{figure}

The parameter $q$ attains its maximum at the critical value
\begin{equation}
    \alpha_0^*\approx 1.391\approx 79.7^\circ. \label{a0val}
\end{equation}
For values of $\alpha_0$ below $\alpha_0^*$, 
$\mathbf{x}_0$ is the lower equilibrium. At $\alpha_0=\alpha_0^*$ the two equilibria coincide, and for $\alpha_0$ above $\alpha_0^*$, $\mathbf{x}_1$ is the lower equilibrium. See Figure \ref{fnet(alpha)}.

\section{Stability of equilibria}\label{sec:stability}

The Jacobian of the system \eqref{eq:sys1}-\eqref{eq:sys4} is given by 
\begin{equation}
 J=\left(
\begin{array}{cccc}
 0 & 1 & 0 & 0 \\
 f_x & 0 & f_y & 0 \\
 0 & 0 & 0 & 1 \\
 h_x & 0 & h_y & 0
\end{array}
\right).
\end{equation}
At $x=0$, we find that $f_y=h_x=0$, so the eigenvalues at either equilibrium are 
\begin{equation}
 \lambda_{1,2}=\pm\sqrt{f_x},\quad
 \lambda_{3,4}=\pm\sqrt{h_y}.
\end{equation}
As predicted by the reversibility of the system, \cite{lamb} for all values of $\alpha_0$ the eigenvalues come in pairs $\{\lambda,-\lambda\}$.

Using Equation \eqref{equil2}, the partial derivatives can be written in terms of the $y$-coordinate of the equilibrium as 
\begin{align}
    f_x &= -\frac{1458 y^5}{\left(81 y^4+1\right)^{3/2}} \\
    h_y &= \frac{18 y }{\left(81 y^4+1\right)^{3/2}}\left(\frac{81 y^4 \left(\sqrt{81 y^4+1}-4\right)}{\sqrt{81 y^4+1}+4}-5\right)
\end{align}
where $y=y_0$ or $y_1$.

At $\mathbf{x}_0$, we have in terms of $\alpha_0$
\begin{align}
 f_x &= -6 \sqrt{\sin ^5 \alpha_0 \cos \alpha_0} \\
 h_y &= -\frac{6 \sqrt{\tan \alpha _0 }}{\sec \alpha _0+4} \big(16 \cos \alpha _0  +3 \cos \left(2 \alpha _0\right) \nonumber \\
 & \qquad +4 \cos \left(3 \alpha _0\right)+2\big)  
\end{align}

For $\alpha_0<\alpha_0^*$, all four eigenvalues of the Jacobian at $\mathbf{x}_0$ are pure imaginary, so the equilibrium of the linearized system is (Lyapunov) stable. By the symmetry of the system, $\mathbf{x}_0$ is a nonlinear center.

\begin{figure}
 \centering
 \includegraphics[width=0.45\textwidth]{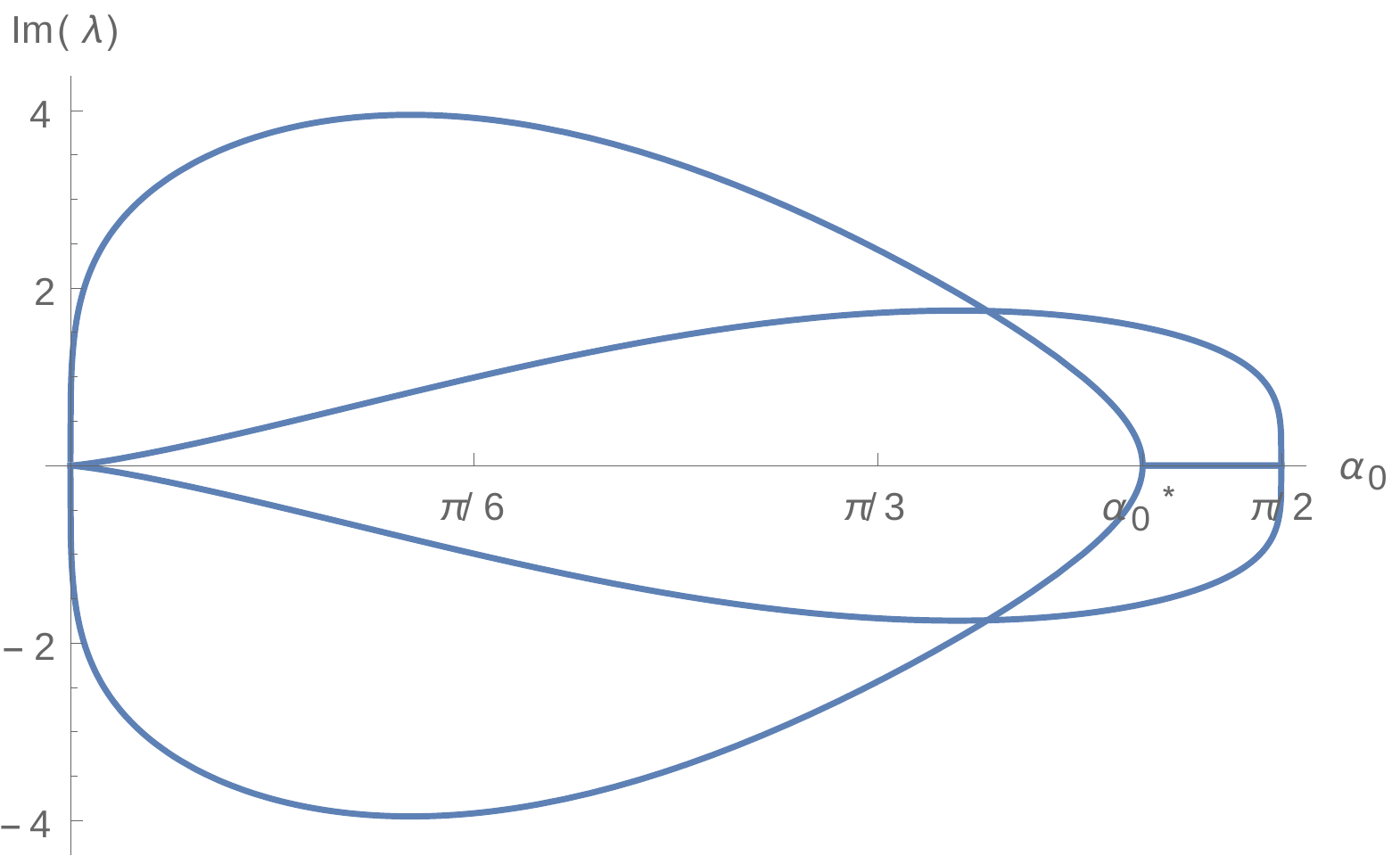}
 \includegraphics[width=0.45\textwidth]{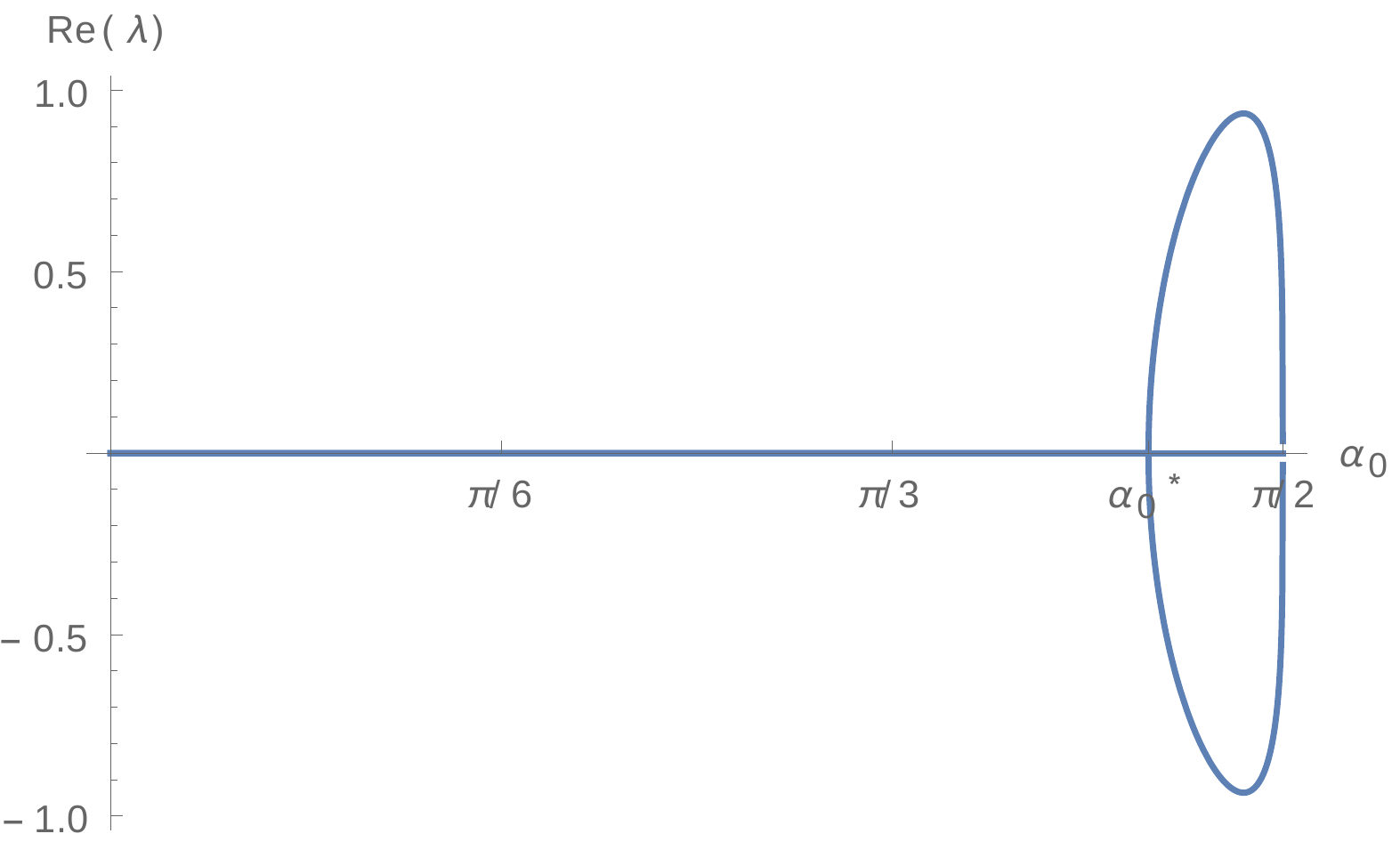}
 \caption{\textit{Stability of equilibrium $\mathbf{x}_0$: (a) Imaginary part of the four eigenvalues, vs $\alpha_0$. (b) Real part of the eigenvalues.}}
 \label{evs-fig}
\end{figure}

At $\alpha_0=\alpha_0^*$, $\lambda_3$ and $\lambda_4$ pass through 0 and become real, so the equilibrium becomes a saddle.
See Figure \ref{evs-fig}.

\begin{figure}
 \centering
 \includegraphics[width=0.45\textwidth]{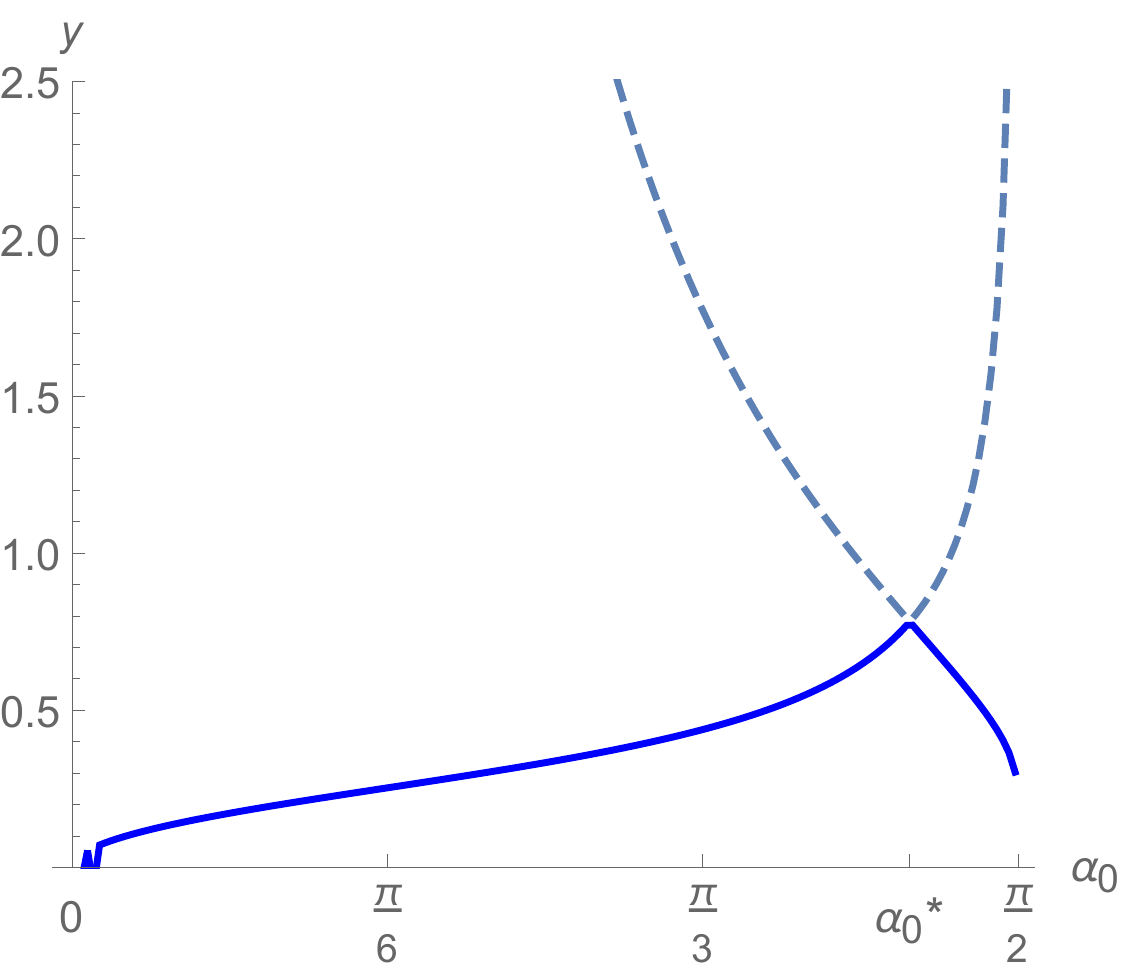}
 \caption{\textit{Bifurcation diagram: branches $\mathbf{x}_0$ and $\mathbf{x}_1$ cross, in this projection of 4D space onto $y$ coordinate of the stable (solid) and unstable (dashed) equilibria as functions of $\alpha_0$. The small wiggles near $\alpha_0=0$ are a numerical artifact.}}
 \label{fig:bif-diag}
\end{figure}

This change of stability occurs in a transcritical (exchange of stability) bifurcation with $\mathbf{x}_1$. See Figure \ref{fig:bif-diag}, in which the $y$ coordinate of the equilibria is given by Equation \eqref{equil2}.

\begin{figure}
 \centering
 \includegraphics[width=0.45\textwidth]{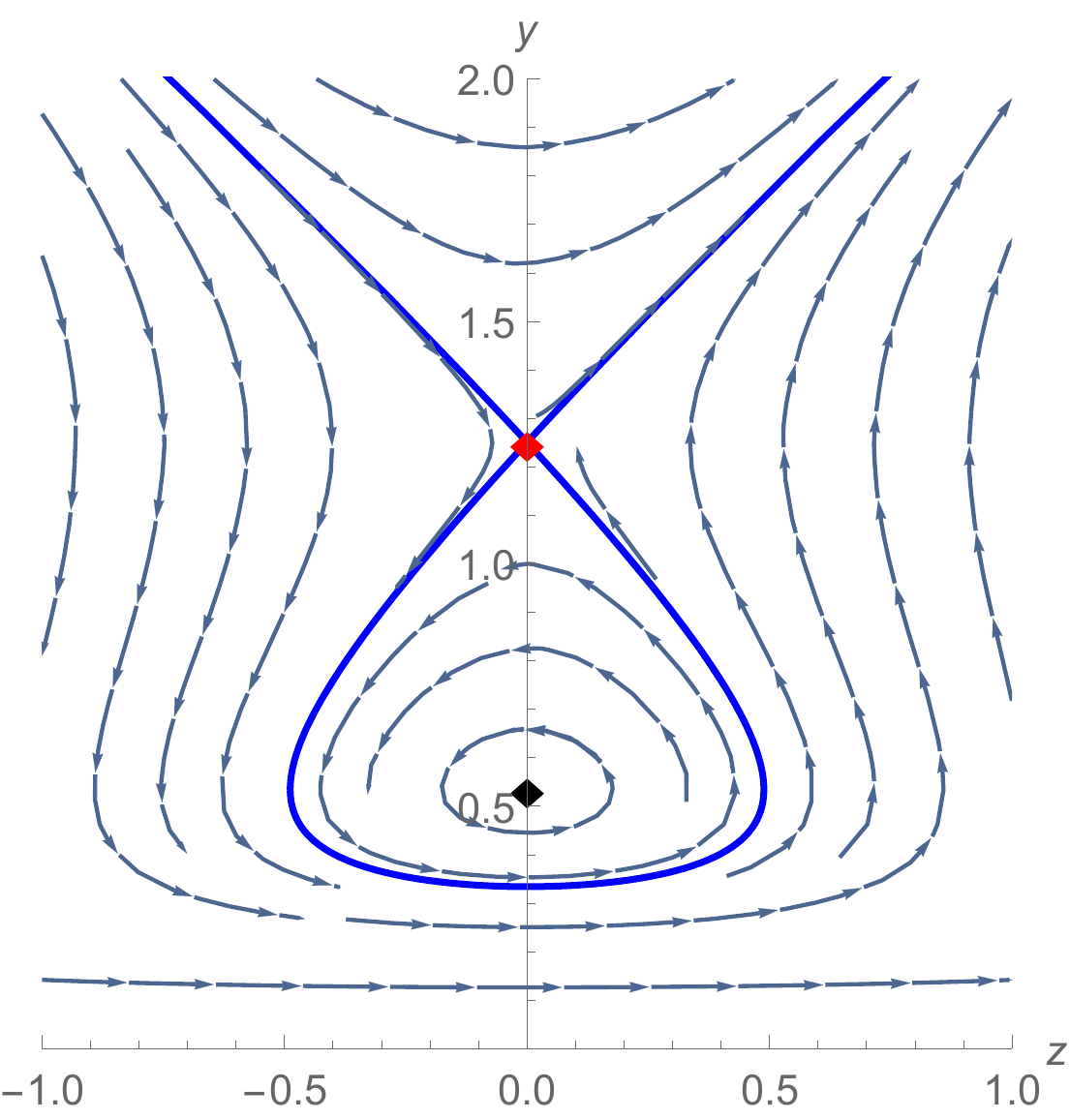}
 \caption{\textit{Phase portrait of the reduced system \eqref{y-red}-\eqref{z-red} in the $(z,y)$ plane for $\alpha_0=1.2$, showing the  stable ($\mathbf{x}_0$, black) and saddle ($\mathbf{x}_1$,red) equilibria, and the stable and unstable manifolds of the saddle (blue). }}
 \label{yz-phase}
\end{figure}

Notice that invariance of the dynamics under $S$ implies that the two-dimensional manifold $x=w=0$ is invariant. On this manifold, we have the reduced system
\begin{align}
    \dot y &= z \label{y-red} \\
    \dot z &= h(0,y). \label{z-red}
\end{align}
Taking $x=0$, we can plot a phase portrait of the reduced system Equation \eqref{y-red}-\eqref{z-red} in the $(z,y)$ plane. Below $\alpha_0^*$, the prescribed equilibrium $\mathbf{x}_0$ is stable and surrounded by an island of periodic orbits, while $\mathbf{x}_1$ is a saddle; at $\alpha_0^*$, there is a single equilibrium, a degenerate saddle; and above $\alpha_0^*$, $\mathbf{x}_0$ is a saddle, while $\mathbf{x}_1$ is stable.  Figure \ref{yz-phase} shows the phase portrait for $\alpha_0 = 1.2$.  Note that, by the phase flow, initial conditions outside the island escape to infinity. Escape may anticipate the droplet jumping that occurs for real droplets \cite{enright,liu}.  For $y=y_0$ there are critical velocities and for $z=0$, critical displacements, beyond which escape occurs.  `Escape displacements' correspond to sufficiently flattened as well as stretched Steiner drops. That these can escape was surprising to us. See Supplementary Material for an animation.

Physically, this bifurcation shows that above a critical contact angle, the stable drop configuration does not adopt
the prescribed contact angle at all, but a lower one.

\section{Symmetry and invariant manifolds}\label{sec:sym}

At the  stable equilibrium, the linearized system has two pairs of pure imaginary eigenvalues. The reversible Lyapunov center theorem \cite{devaney,lamb} asserts that for each of these pairs, there is a one-parameter family of nested symmetric periodic orbits. Therefore, we look for two 2-dimensional invariant manifolds containing these families.

We have already observed the first such manifold: as described in the previous section, there is an invariant manifold $x=w=0$. This is easily seen, since $f(x,y)$
is odd in $x$, so taking $x=0$ reduces the dynamical system to Equations \eqref{y-red}-\eqref{z-red}.
Note that this reduced system is Hamiltonian: if $U(y)$ satisfies $U'(y)=-h(0,y)$ then we have the Hamiltonian
\begin{equation}
    H(y,z) = \frac{1}{2}z^2 + U(y)
\end{equation}
which satisfies
\begin{align}
    \frac{\partial H}{\partial y} &= U'(y) = -h(0,y) = -\dot z \\
    \frac{\partial H}{\partial z} &= z = \dot y
\end{align}
and is constant in time:
\begin{align}
    \frac{dH}{dt} &= \frac{d}{dt}\left(\frac{1}{2}z^2 + U(y)\right) \\
    &= z \dot z - U'(y) \dot y \nonumber \\
    &= z h(0,y) - h(0,y) z = 0
\end{align}
Therefore a region around the center $\mathbf{x}_0$ in this invariant manifold is filled with periodic orbits, corresponding to a vertically oscillating ``bouncing'' mode of the drop. We will refer to the manifold $x=w=0$ as the bouncing manifold for this reason.

The second type of invariant 2-dimensional manifold is less obvious. We find that there are invariant manifolds containing the equilibria, where $y=g(x,w)$. For most values of $\alpha_0$, $g(x,w)$ is relatively flat, so the trajectories in this manifold exhibit a near-horizontally oscillating ``rocking'' motion. The reversible KAM theory\cite{sevryuk91} predicts that a region around the center in this manifold is filled with periodic orbits. We will call these orbits the rocking mode of the drop, and this type of invariant manifold the rocking manifold. See Supplementary Material for animations of bouncing and rocking Steiner drops. 

To find this invariant manifold, we let $y=g(x,w)$, so the dynamical system reduces to
\begin{align}
 \dot x &= w \label{xdot-red}\\
 \dot w &= f(x,g(x,w)). \label{wdot-red}
\end{align}
By the chain rule, we have
\begin{equation}
 z = \dot y = g_x \dot x + g_w \dot w = g_x w + g_w f
\end{equation}
and
\begin{align} \label{zdot}
 \dot z &= f g_x+g_w \left(w f_x+f_y \left(f g_w+w g_x\right)\right) \nonumber \\
 & \quad +f \left(f g_{ww}+w g_{xw}\right)+w \left(f g_{xw}+w g_{xx}\right) 
\end{align}
where $f$ and its derivatives are evaluated at $(x,g(x,w))$.

By hypothesis, 
\begin{equation} \label{zdot-hyp}
 \dot z = h(x,g(x,w)).
 \end{equation}
We express $g$ as a power series,
\begin{equation}
 g(x,w)=y_0 + a_1 x + a_2 w + a_3 x^2 + a_4 x w + a_5 w^2 +\dots
\end{equation}
and expand $f(x,g(x,w))$ and $h(x,g(x,w))$ as series in $x$ and $w$. Then, we set Equations  \eqref{zdot} and \eqref{zdot-hyp} equal and collect like terms, 
and thus solve for the coefficients $a_i$.

The linear and quadratic coefficients are 

 \begin{align}
 &a_1 = a_2 = a_4 = 0 \\
 &a_3 = \nonumber\\
 &-\frac{3 \left(14 \cos \alpha _0 +4 \cos \left(2 \alpha _0\right)+6 \cos \left(3 \alpha _0\right)+1\right) \left(8 \cos \alpha _0 +7 \cos \left(2 \alpha _0\right)+8 \cos \left(3 \alpha _0\right)+1\right) \sqrt{\cot \alpha _0} \csc ^2 \alpha _0}{4 \left(\left(14 \cos \alpha _0 +4 \cos \left(2 \alpha _0\right)+6 \cos \left(3 \alpha _0\right)+1\right){}^2 \csc ^4 \alpha _0 -4 \left(4 \cos \alpha _0+1\right)^2\right)} \\
 &a_5 = \frac{\sin \alpha _0 \left(4 \cos \alpha _0 +1\right) \left(8 \cos \alpha _0 +7 \cos \left(2 \alpha _0\right)+8 \cos \left(3 \alpha _0\right)+1\right)}{2  \left(208 \cos \alpha _0 +388 \cos \left(2 \alpha _0\right)+148 \cos \left(3 \alpha _0\right)+191 \cos \left(4 \alpha _0\right)+44 \cos \left(5 \alpha _0\right)+32 \cos \left(6 \alpha _0\right)+239\right)} 
\end{align}

\begin{figure}
 \centering
 \includegraphics[width=0.45\textwidth]{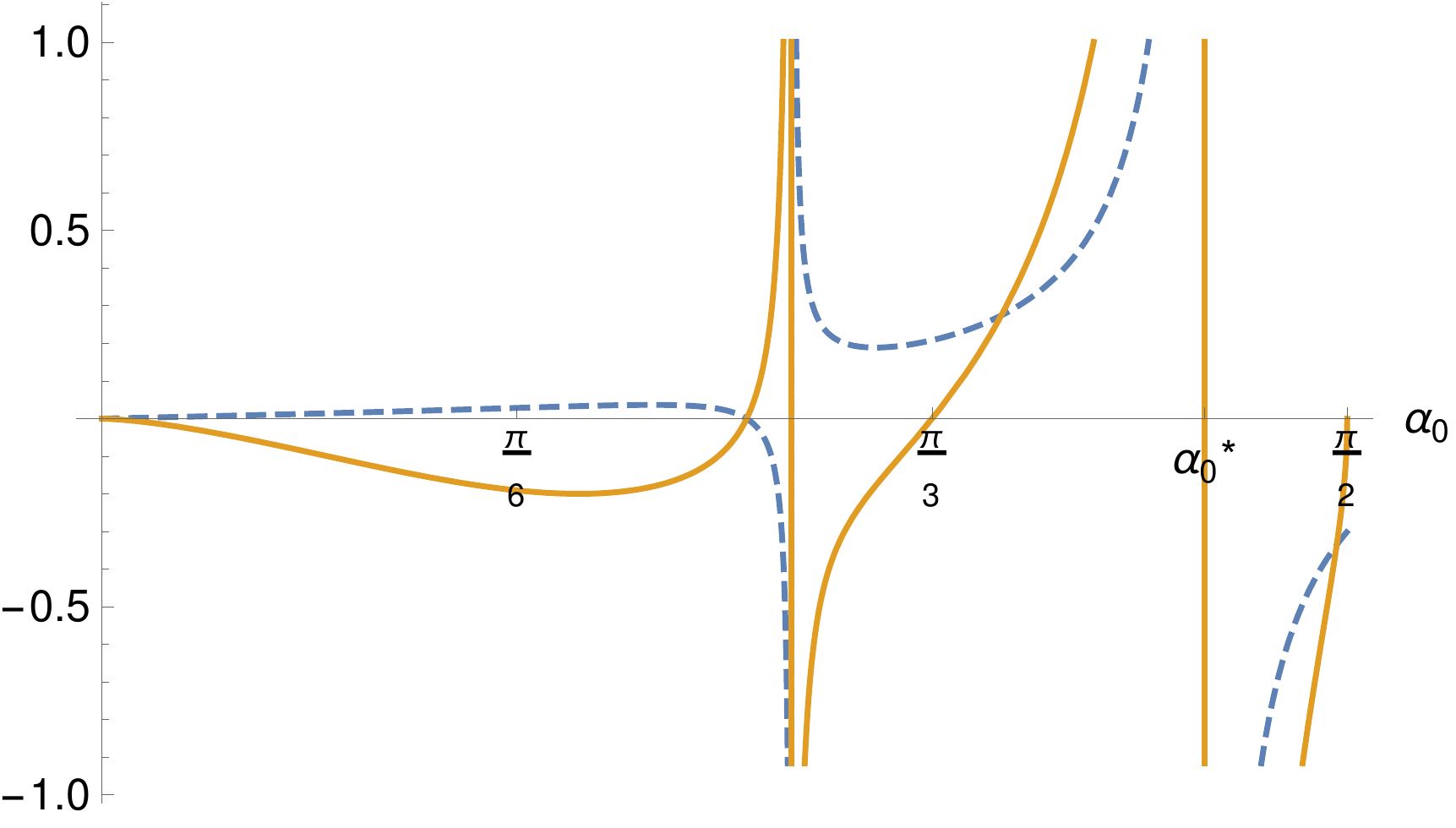}
 \caption{\textit{Coefficients of the $x^2$ (solid curve) and $w^2$ (dashed) terms of the rocking invariant manifold $y=g(x,w)$. Singularities arise at $\alpha_0 = \alpha_0^*$ and $\alpha_0^\dag$.}}
 \label{a35-fig}
\end{figure}

See Figure \ref{a35-fig}. Notice that these have singularities at $\alpha_0=\alpha_0^*$, which is to be expected, and at $\alpha_0^\dag\approx 0.870\approx 49.8^\circ$, which is surprising.

\begin{figure}
 \centering
 \includegraphics[width=0.45\textwidth]{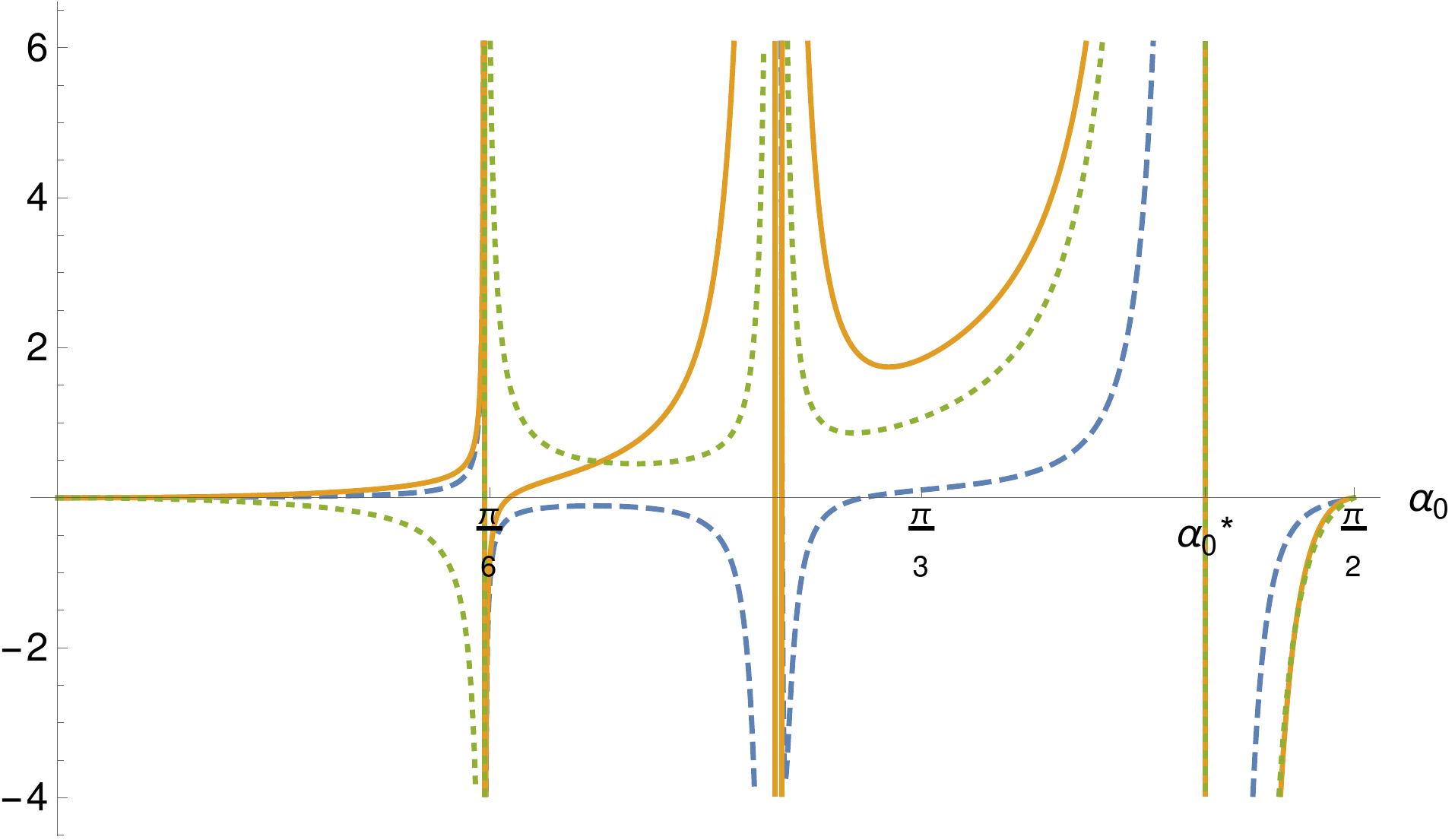}
 \caption{\textit{Coefficients of the $x^4$ (solid curve), $x^2 w^2$ (dotted), and $w^4$ (dashed) terms of the rocking invariant manifold $y=g(x,w)$. Singularities arise at $\alpha_0 = \alpha_0^*$ and $\alpha_0^\dag$, and at the new singular value $\alpha_0\approx 0.517\approx 29.6^\circ$.  }}
 \label{fig:more-sing}
\end{figure}

Computing the coefficients of higher-degree terms, we find that the cubic terms vanish, as do the $xw^3$ and $x^3 w$ terms. The coefficients of the $x^4$, $x^2 w^2$, and $w^4$ terms are cumbersome but can be computed numerically. Like the quadratic coefficients, the quartic coefficients have singularities at $\alpha_0^*$ and  $\alpha_0^\dag$; there is also a singular value $\alpha_0\approx 0.517\approx 29.6^\circ$. See Figure \ref{fig:more-sing}.

\begin{figure}
 \centering
 \includegraphics[width=0.45\textwidth]{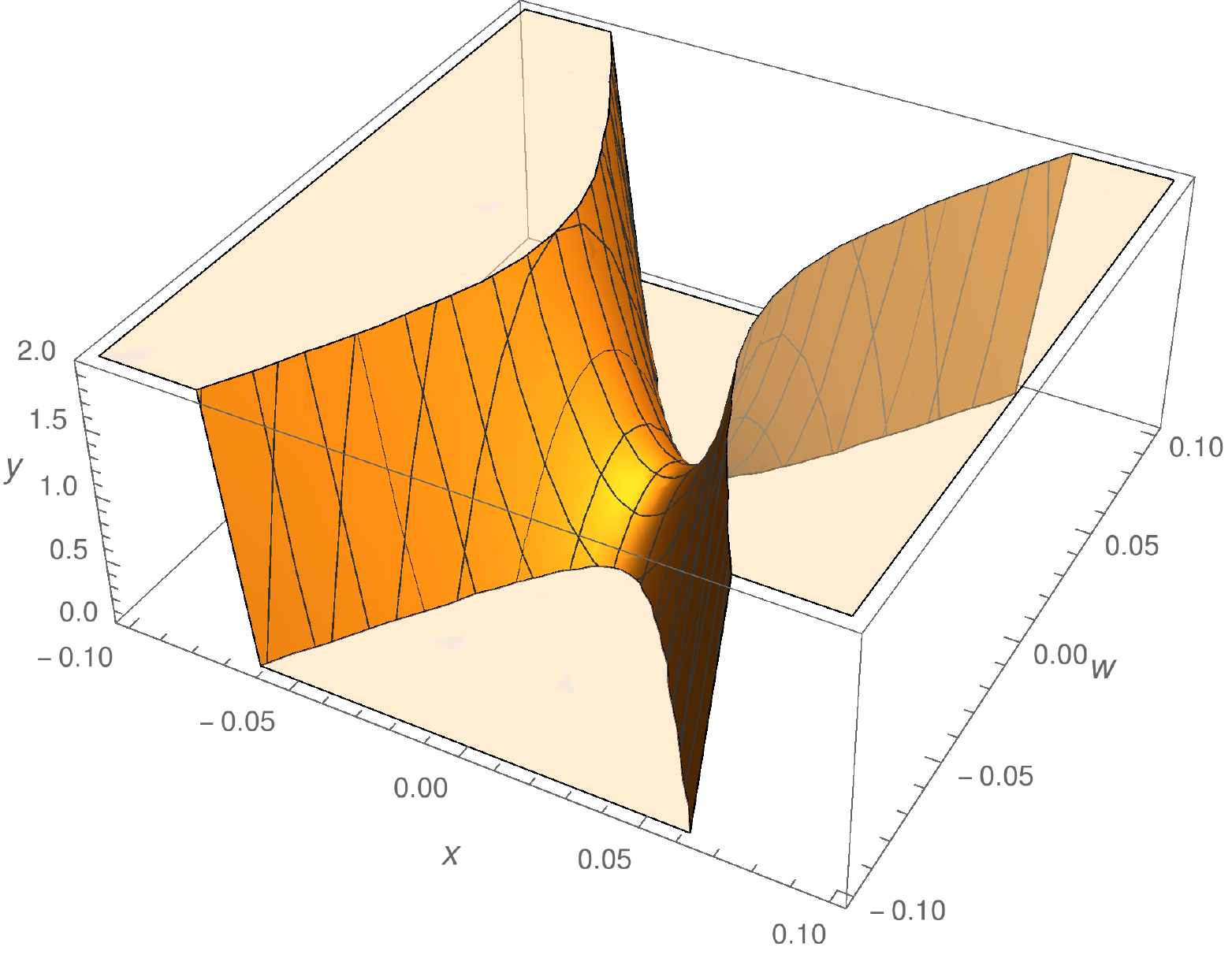}
 \includegraphics[width=0.45\textwidth]{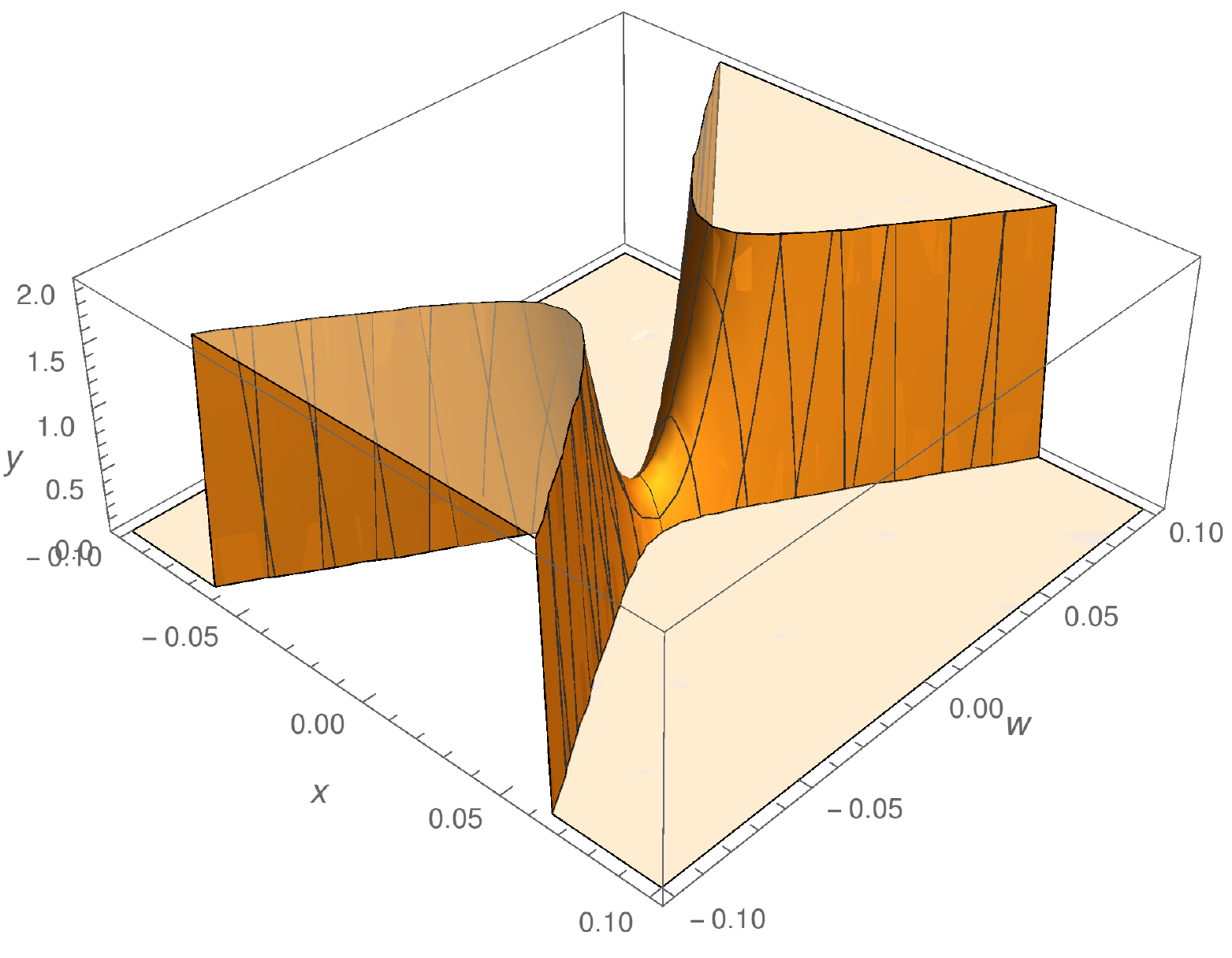}
 \caption{\textit{Plots of the rocking invariant manifold $y=g(x,w)$ for values of $\alpha_0$ near the singular value $\alpha_0^\dag\approx0.870\approx 49.8^\circ$. 
 As $\alpha_0$ approaches this value, the manifold is a saddle (hyperbolic paraboloid) whose curvature at the origin goes to $-\infty$. 
 (a) $\alpha_0=\alpha_0^\dag -0.000004$. (b) $\alpha_0=\alpha_0^\dag + 0.000002$.}}
\end{figure}

\begin{figure}
 \centering
 \includegraphics[width=0.45\textwidth]{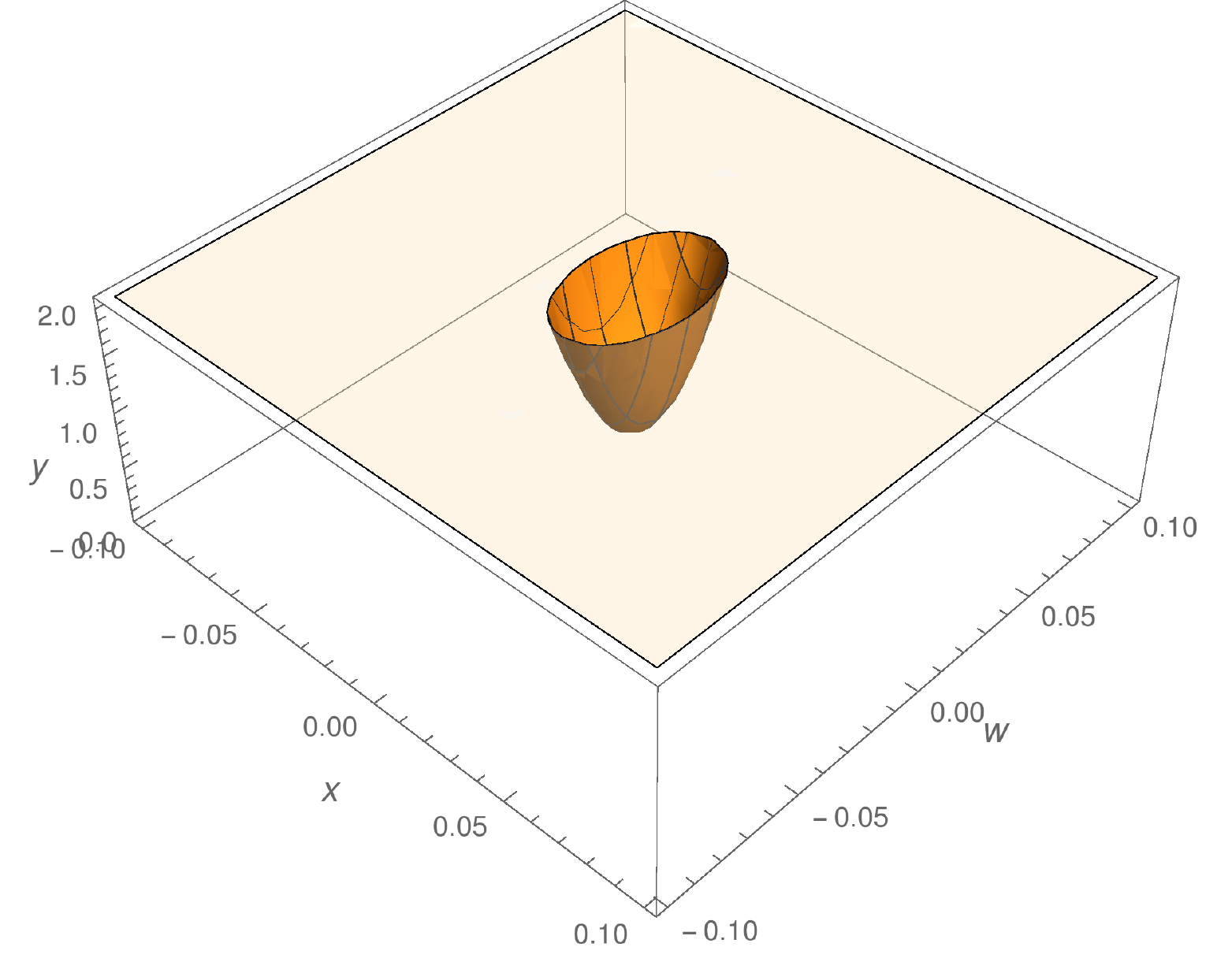}
 \includegraphics[width=0.45\textwidth]{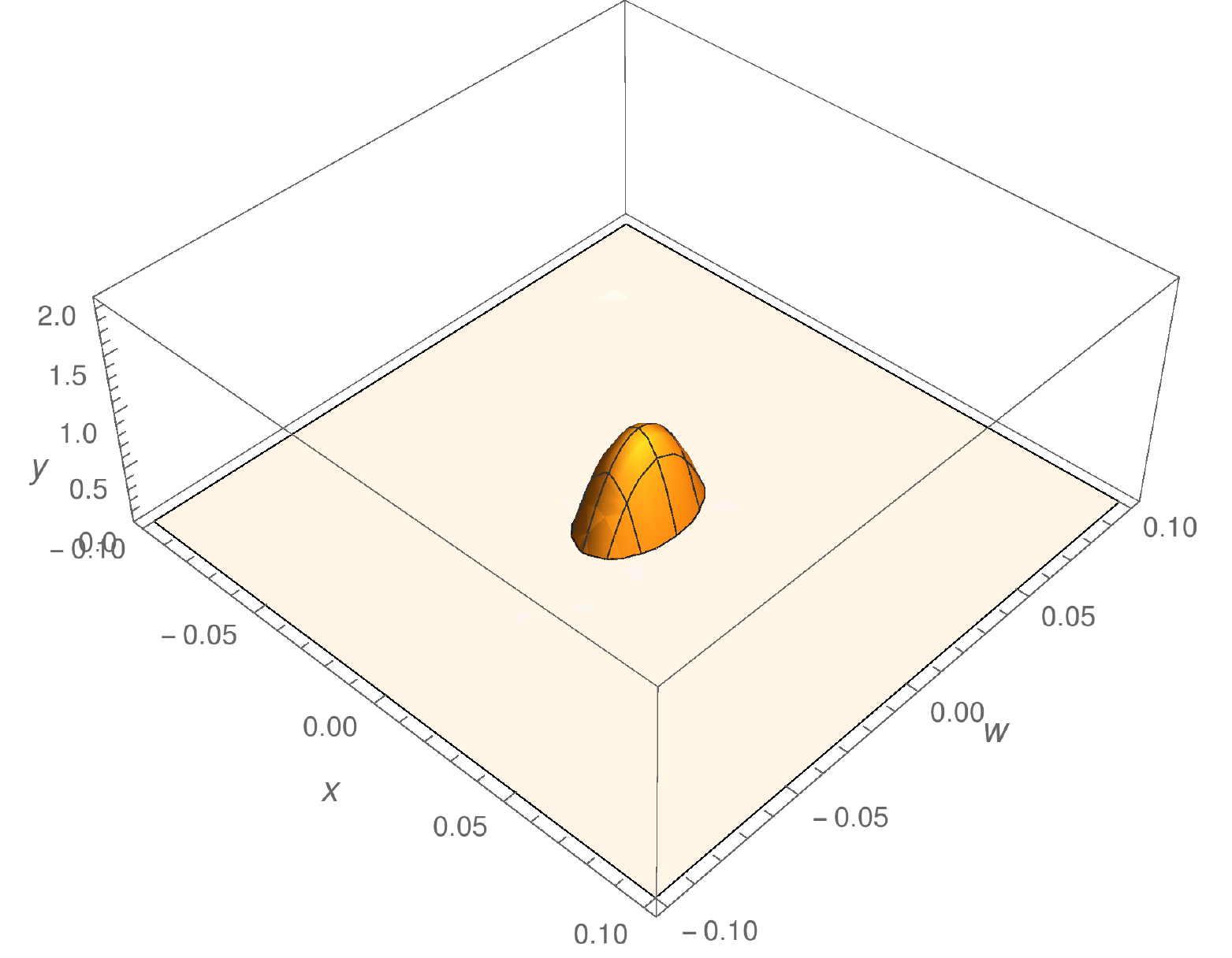}
 \caption{\textit{Plots of the rocking invariant manifold $y=g(x,w)$ for values of $\alpha_0$ near the singular value $\alpha_0^*\approx1.391\approx 79.7^\circ$. 
 As $\alpha_0$ approaches this value, the manifold is a paraboloid whose curvature at the origin goes to $\infty$.
 (a) $\alpha_0=\alpha_0^* -0.00005$. (b) $\alpha_0=\alpha_0^* + 0.00005$.}}
\end{figure}

\begin{figure}
 \centering
 \includegraphics[width=0.38\textwidth]{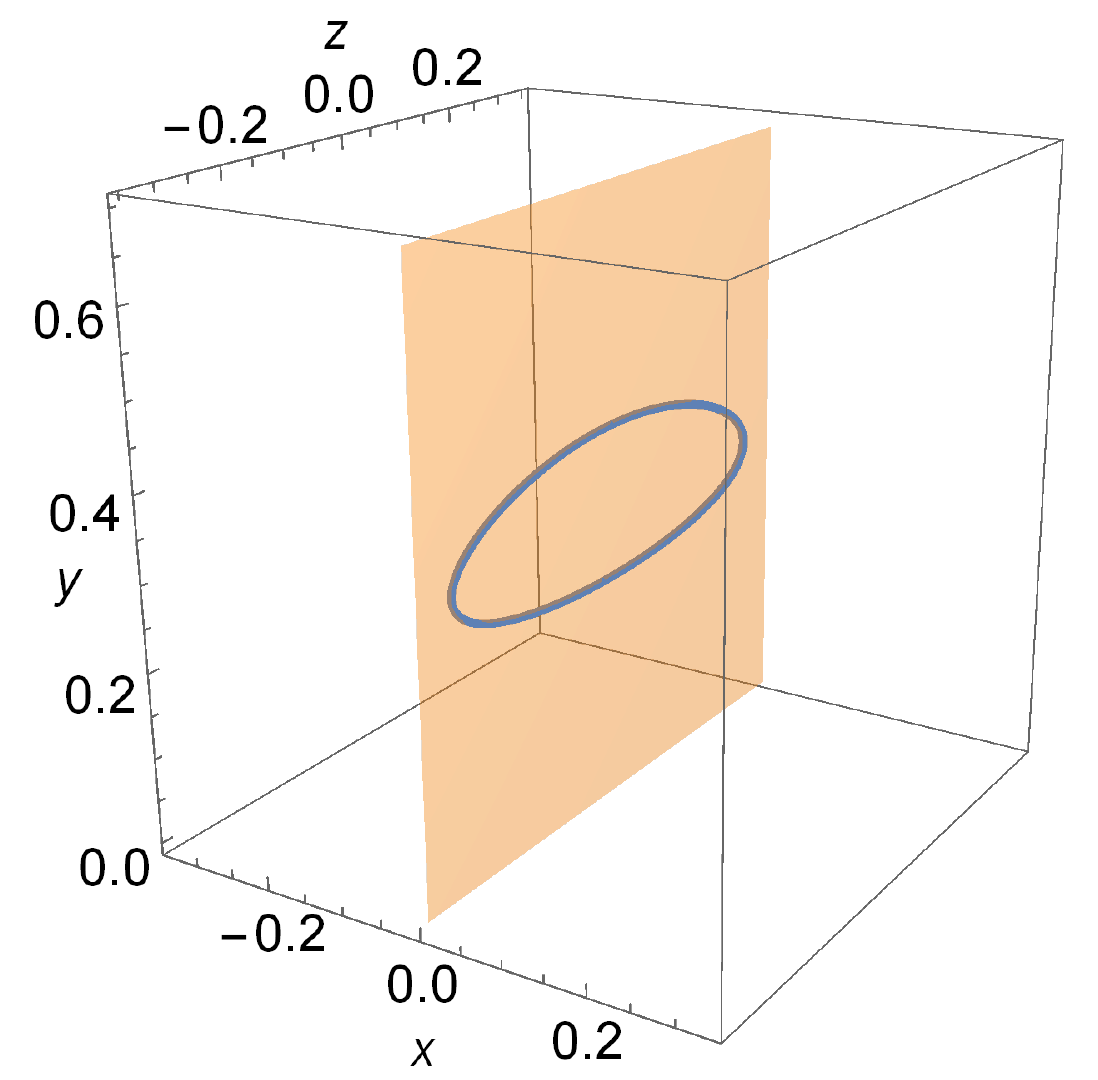}
  \includegraphics[width=0.45\textwidth]{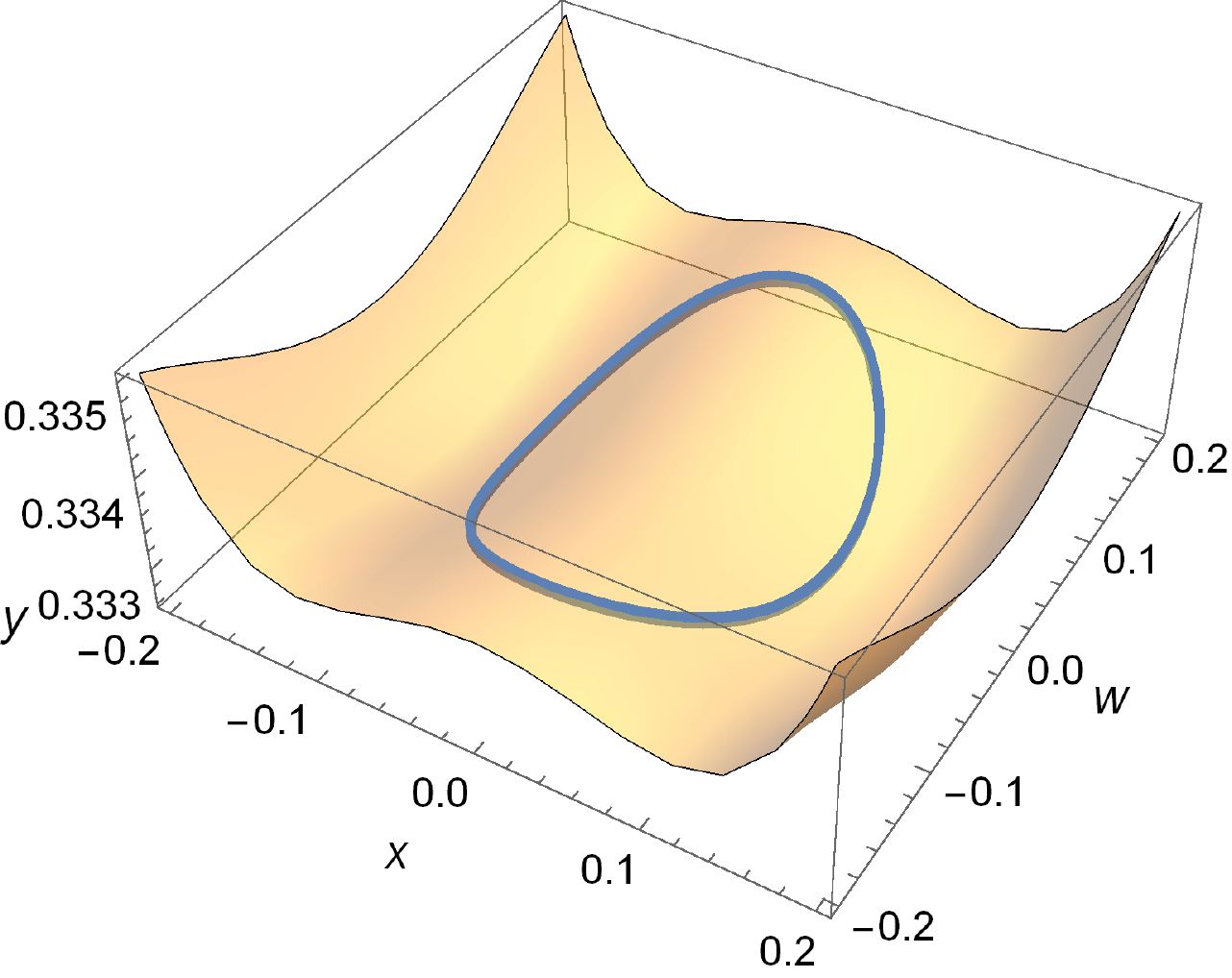}
 \caption{\textit{Numerically integrated trajectories for $\alpha_0=\pi/4$. (a) A trajectory lying in the bouncing
 invariant manifold $x=0$. (b) A trajectory lying in the rocking invariant manifold $y=g(x,w)$}}
 \label{fig:bouncing}
\end{figure}

It makes sense that $g$ should be even with respect to both variables. The time-reversing involutions $G_1$ and $G_2$ applied to the reduced system give
\begin{align}
    \dot x &= w \\
    -\dot w &= f(-x,g(-x,w)) \nonumber \\
    &= -f(x,g(-x,w))
\end{align}
and
\begin{align}
    -\dot x &= -w  \\
    \dot w &= f(x,g(x,-w))
\end{align}
respectively. If $g$ is even with respect to both variables, the system constrained to the invariant manifold inherits the full system's invariance under $G_1$ and $G_2$.

Away from the singular values of $\alpha_0$, the symmetry of the reduced system Equations \eqref{xdot-red}-\eqref{wdot-red} lying in the manifold $y=g(x,w)$ shows that a region around the origin in this manifold is filled with periodic orbits.

We use this same method to find the invariant manifold containing the second equilibrium $\mathbf{x}_1$, for values of $\alpha_0$ above the transcritical bifurcation at $\alpha_0=\alpha_0^*$. Note that since there is no closed-form expression for $\mathbf{x}_1$, we must find the invariant manifold numerically for a specified value of $\alpha_0$.
For example, the result when $\alpha_0=1.45$ is
\begin{equation}
 y=g(x,w) = 0.6504+2.1064 x^2+0.91936 w^2+\dots.
\end{equation}

\section{Numerical solutions: Torus trajectories}\label{sec:torus}

The reversible KAM theory \cite{sevryuk91} says that in the $2n$-dimensional phase space of a smooth dynamical system that is reversible with respect to an involution which fixes an $n$-dimensional submanifold, for each $0\leq m\leq n$ there is an $m$-parameter family of invariant $m$-tori. The union of these tori has positive $2m$-dimensional measure. 

In particular, in our 4-dimensional phase space with two involutions which each fix a 2-dimensional subspace, there is a 0-dimensional manifold of 0-tori (the stable equilibrium); 2-dimensional manifolds of 1-tori (the periodic orbits discussed in the previous section) and a 4-dimensional region of 2-tori. The quasiperiodic frequencies of trajectories on these 2-tori are strongly incommensurable, \cite{sevryuk91} meaning that each such trajectory is dense in its invariant torus.

\begin{figure}
 \centering
 \includegraphics[width=.45\textwidth]{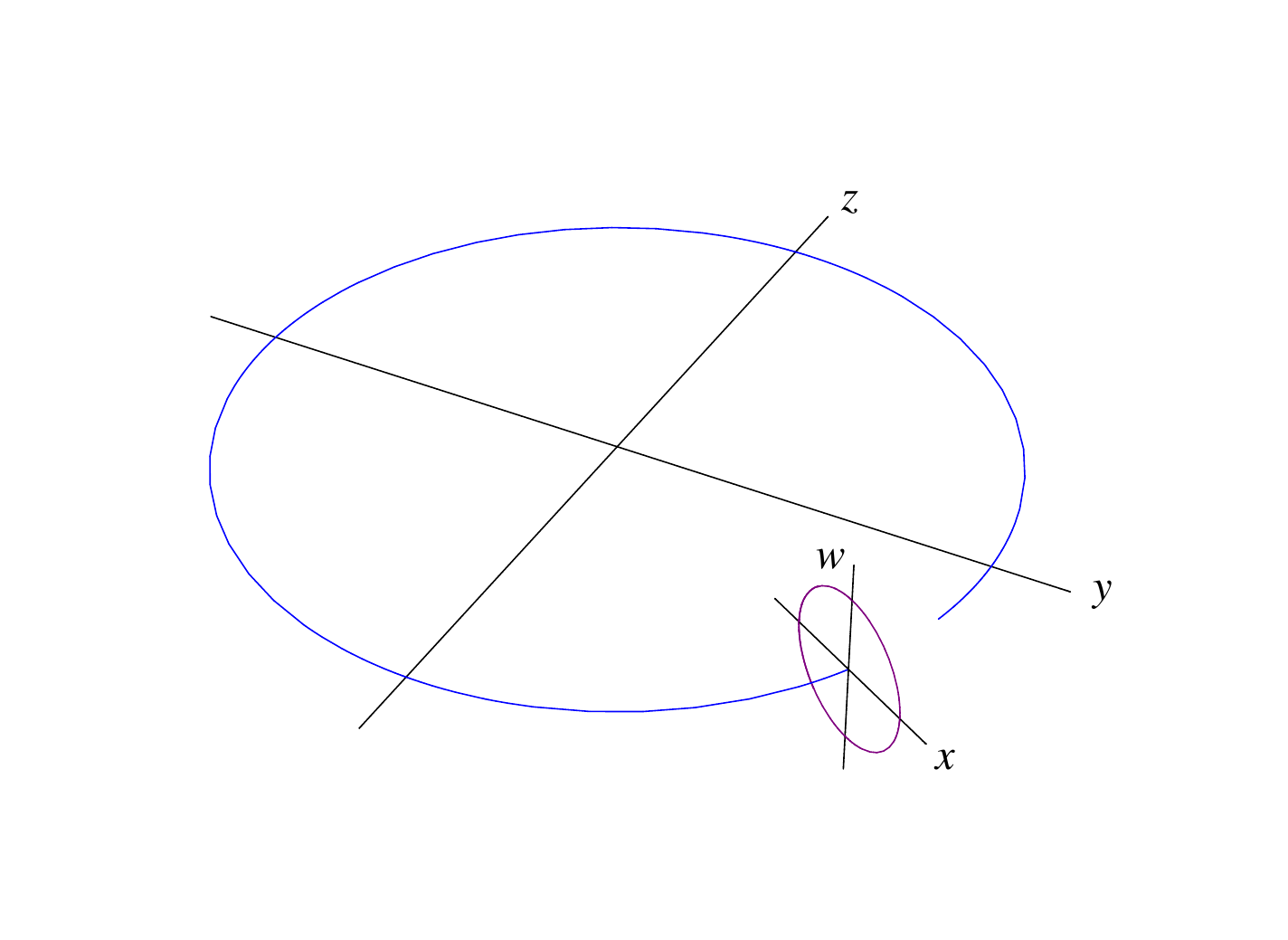}
 \caption{\textit{Coordinate sketch for the torus visualization of the trajectories in 4D.}}
 \label{torus-coords}
\end{figure}

Physically, these torus trajectories combine both vertical and horizontal motion. They can be visualized in three dimensions as lying on tori where one direction of revolution is 
in the $y$-$z$ plane, and the other is in a traveling $x$-$w$ plane. See Figure \ref{torus-coords}.

To show the trajectories in this way, we parameterize them via
\begin{equation}
 \left[y+x\frac{y}{\sqrt{y^2+z^2}},z+x\frac{y}{\sqrt{y^2+z^2}},w\right].
\end{equation}
Several torus trajectories are shown in the Appendix. See Supplementary Material for an animation of this.

\begin{figure}
 \centering
 \includegraphics[width=0.31\textwidth]{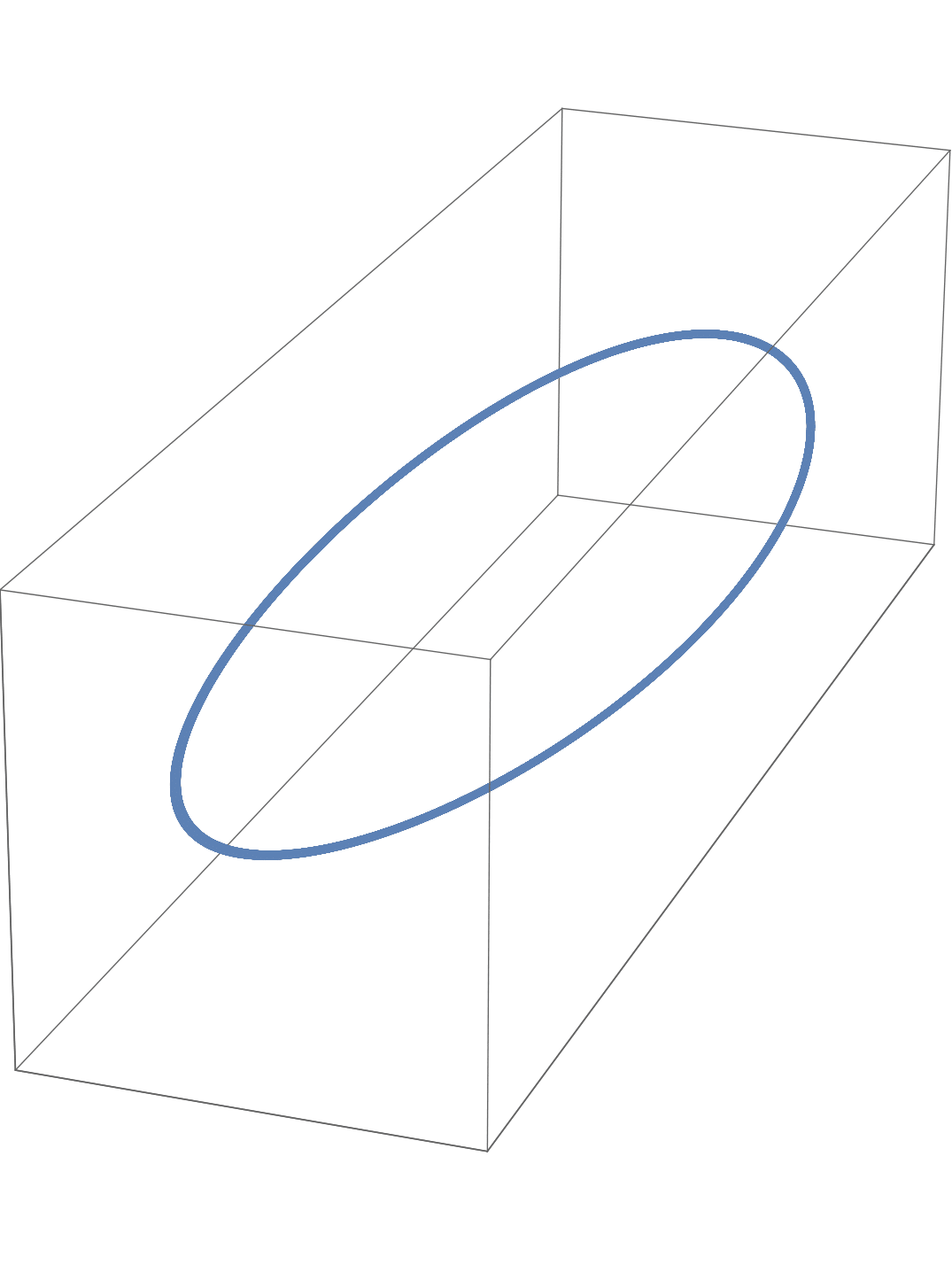}
 \includegraphics[width=0.31\textwidth]{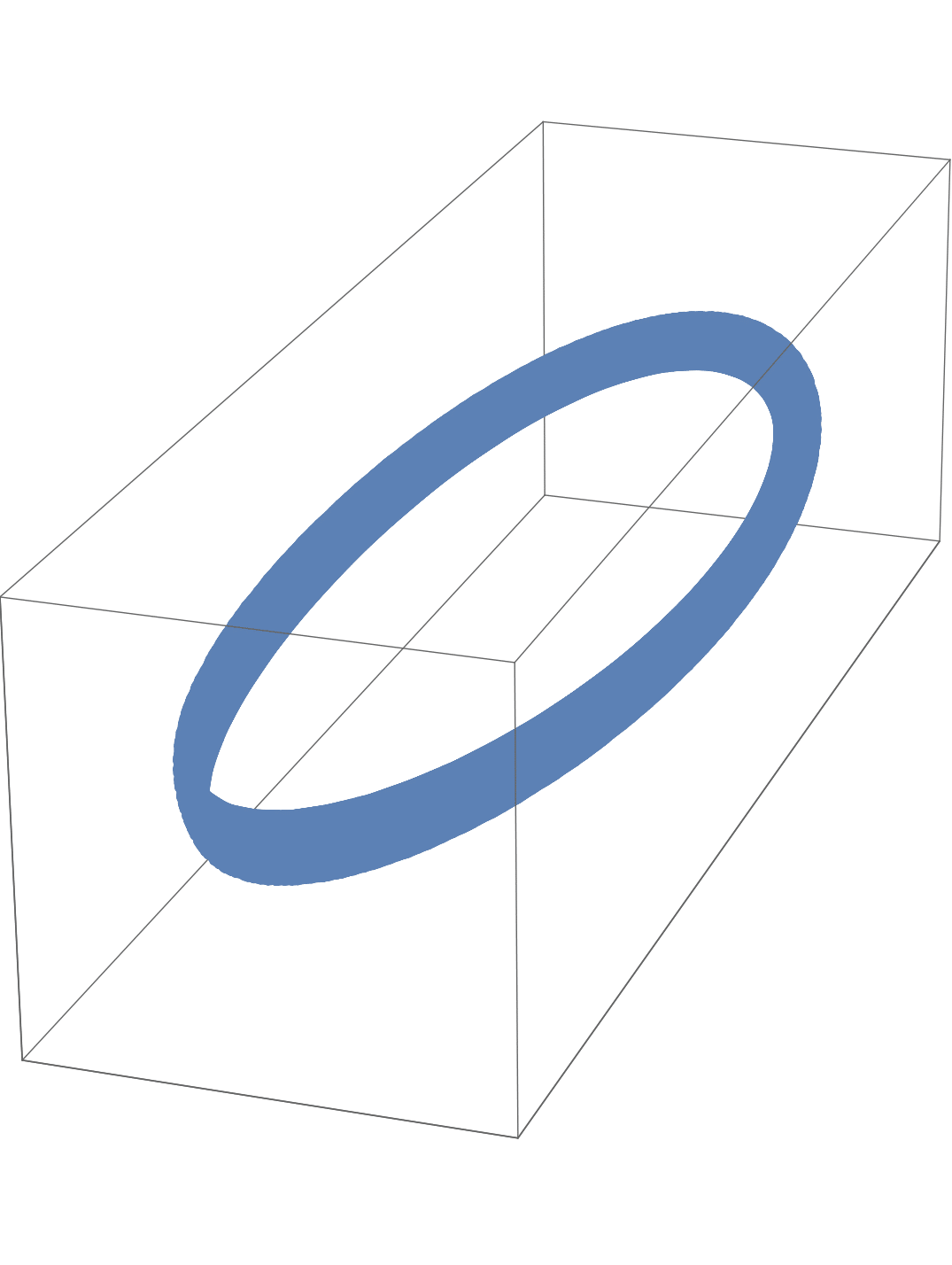}
 \includegraphics[width=0.31\textwidth]{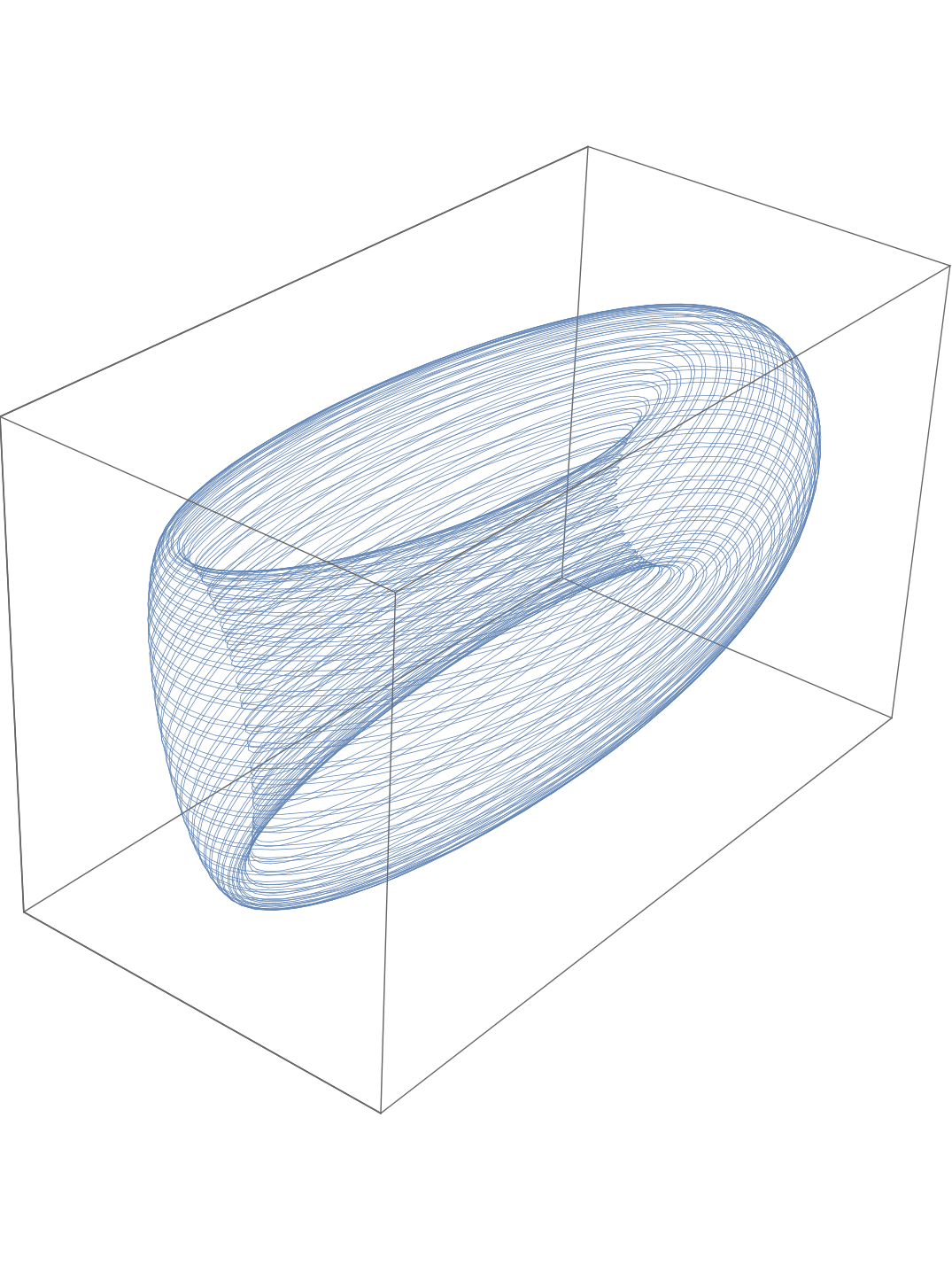} \\
 \includegraphics[width=0.31\textwidth]{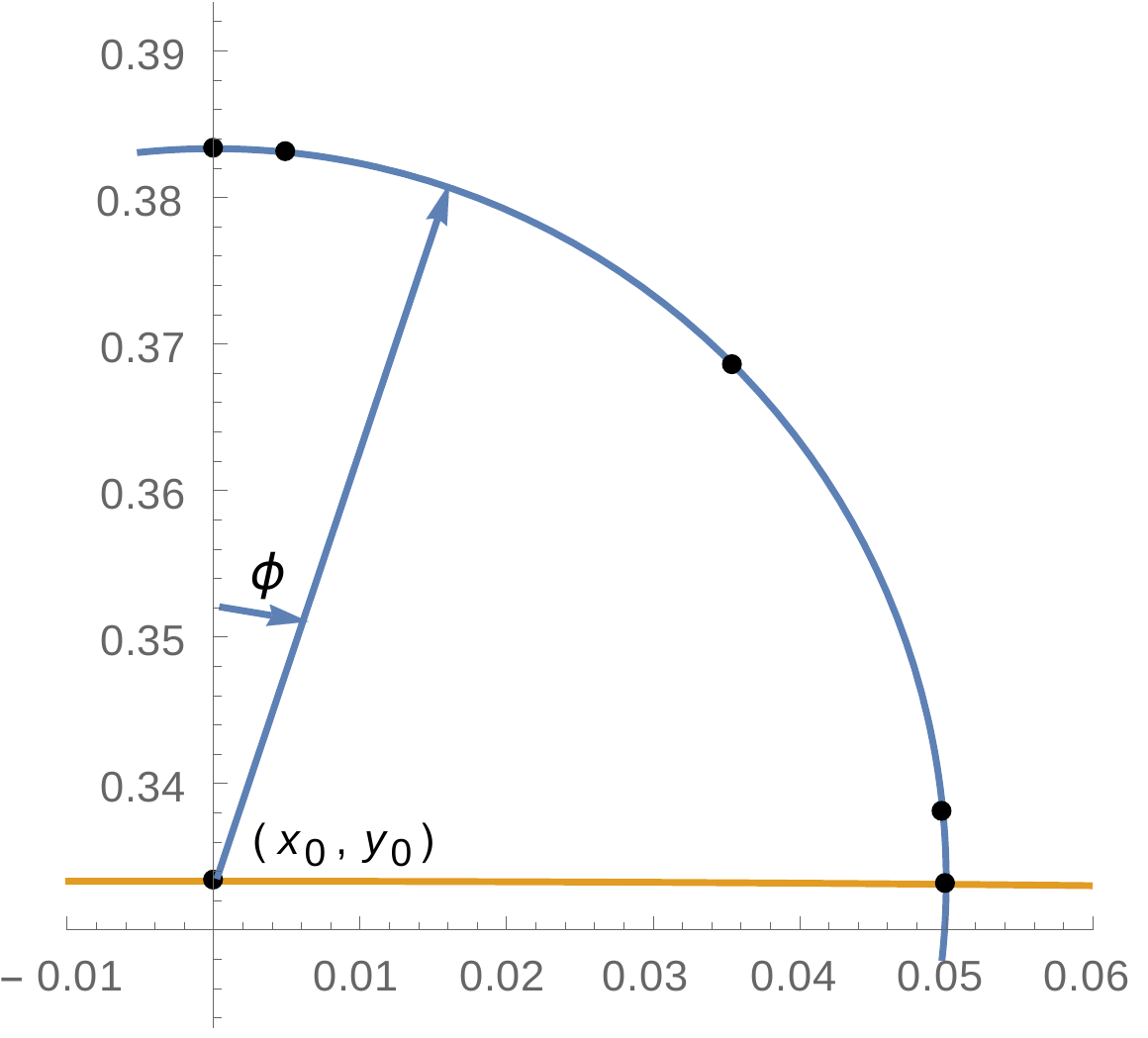}
 \includegraphics[width=0.31\textwidth]{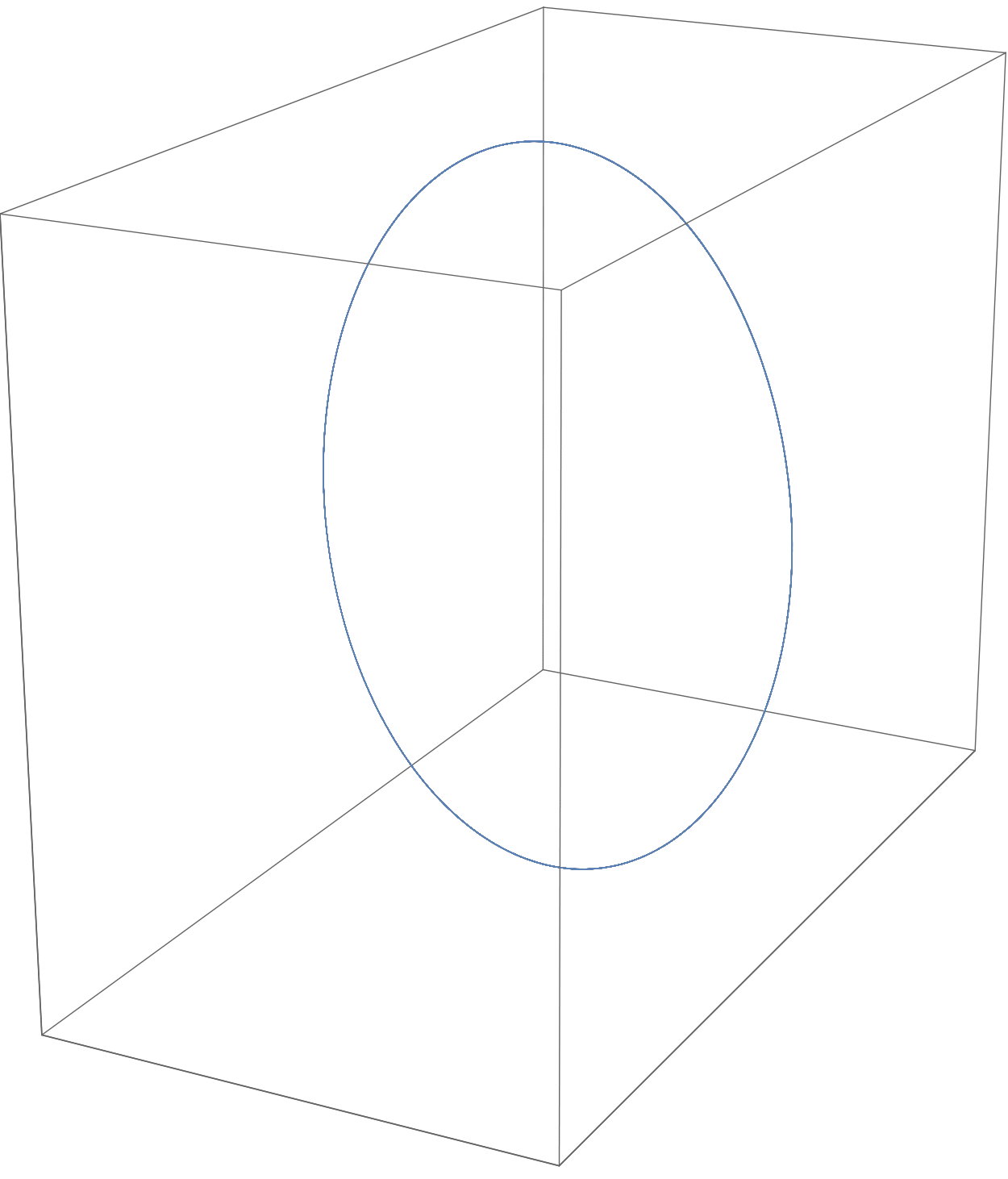}
 \includegraphics[width=0.31\textwidth]{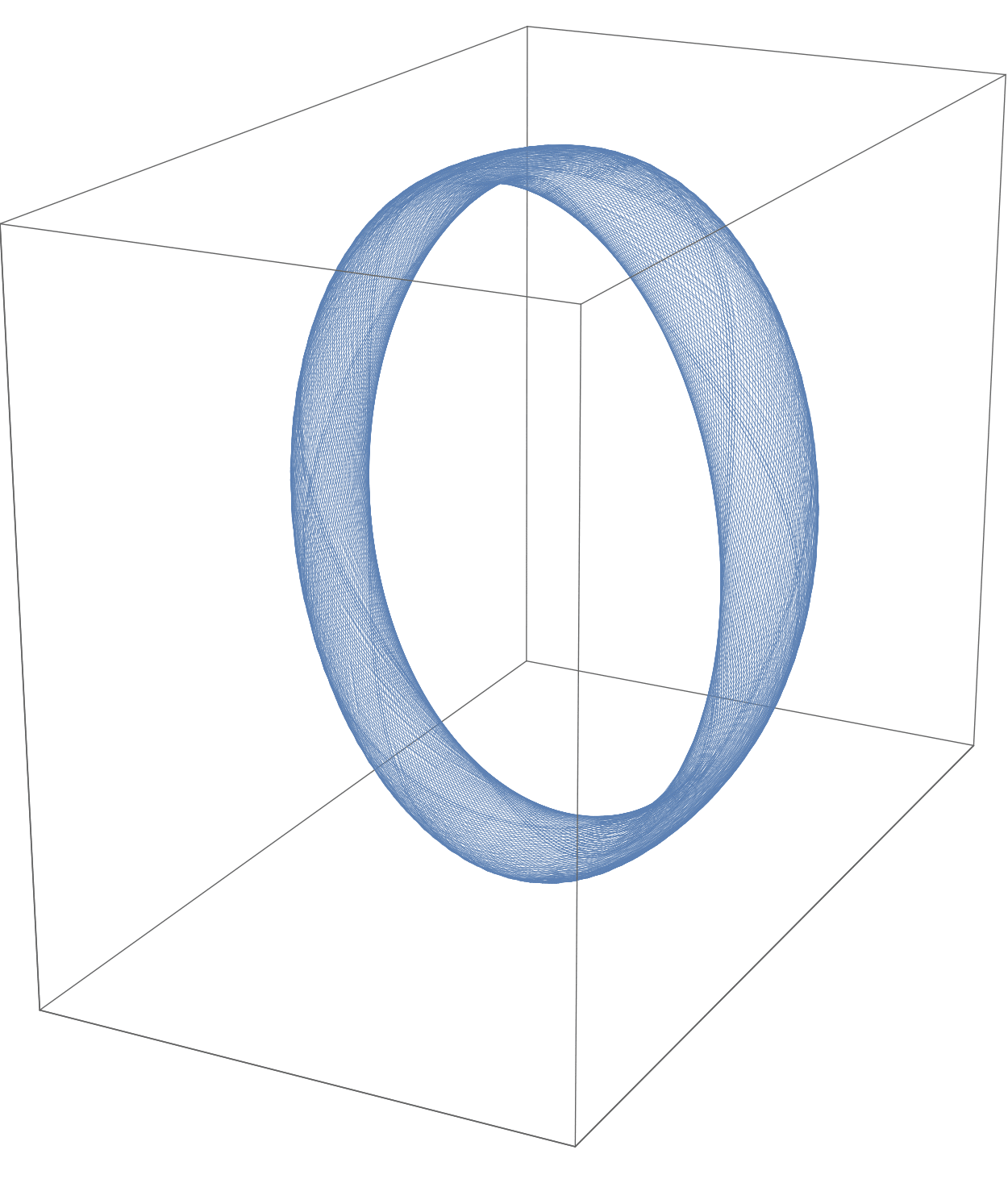}
 \caption{\textit{Lower left: Map of initial conditions shown for torus trajectories for $\alpha_0=\pi/4$. The yellow curve is the graph  of $y=g(x,0)$ in the invariant manifold. 
 Clockwise from upper left: Bouncing mode $\phi=0$; $\phi=0.1$; $\phi=\pi/4$; $\phi=1.475$; rocking mode $\phi=1.575$}}
 \label{torusfig-1}
\end{figure}

\begin{figure}
 \centering
  \includegraphics[width=0.31\textwidth]{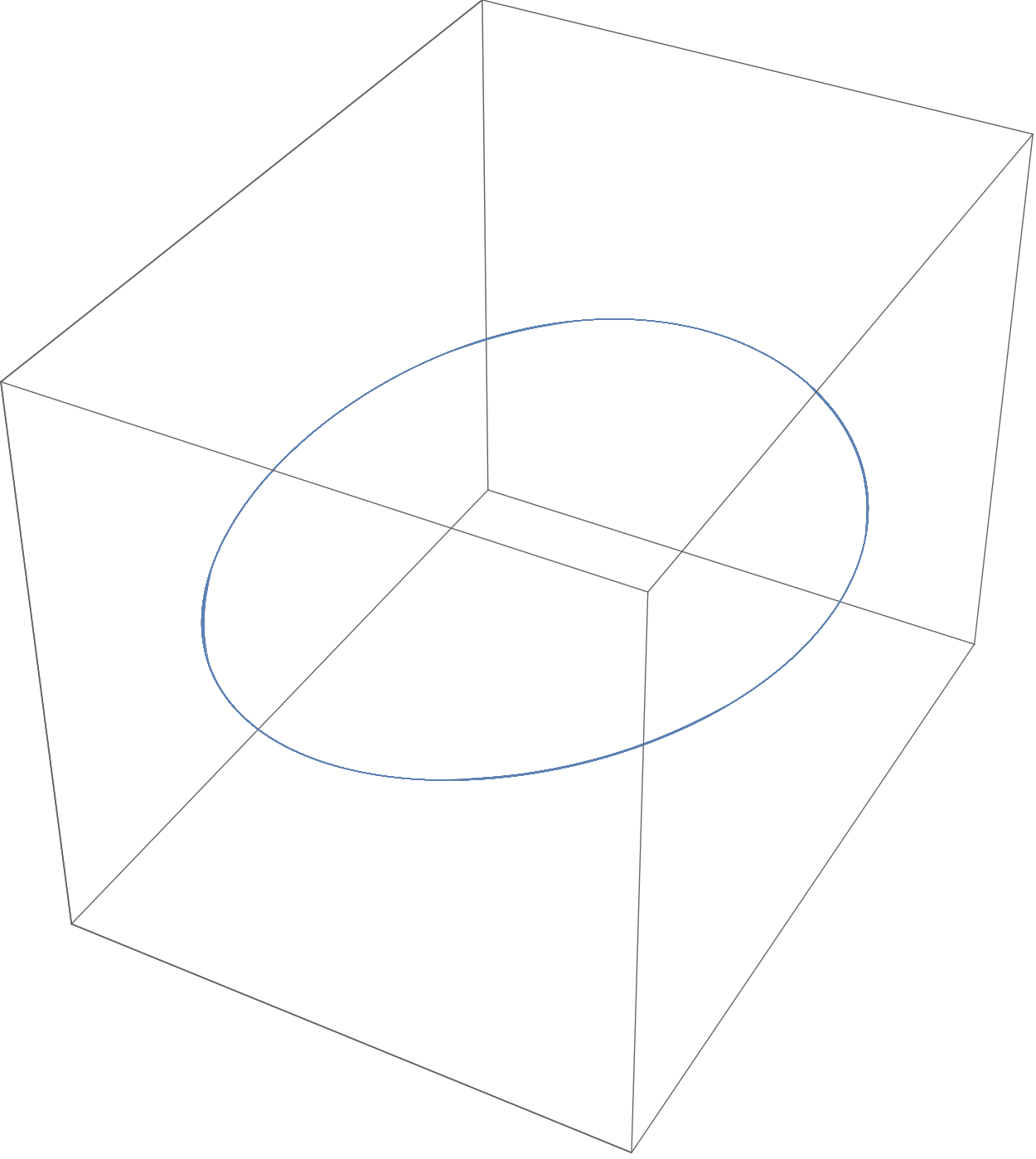}
 \includegraphics[width=0.31\textwidth]{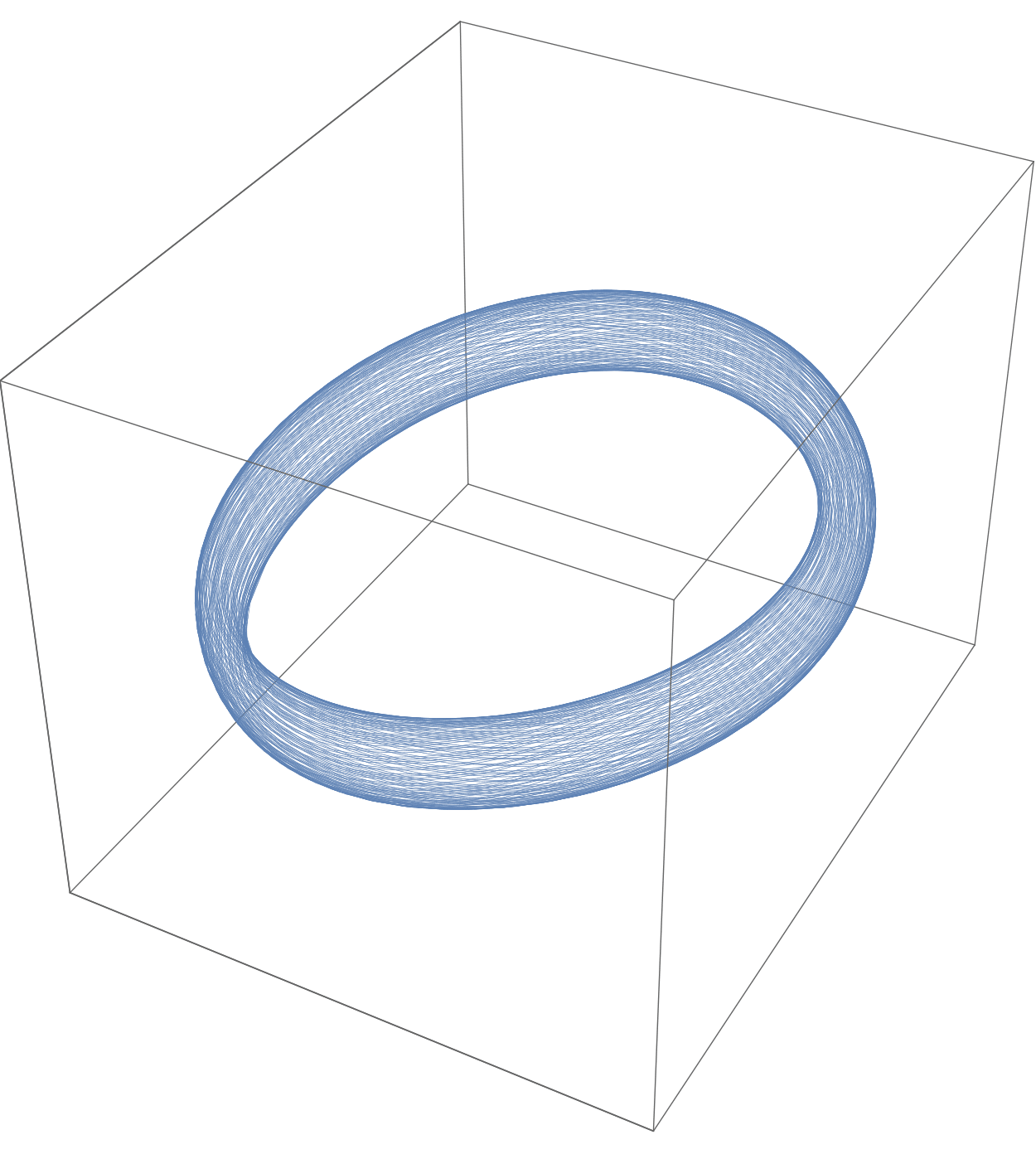}
 \includegraphics[width=0.31\textwidth]{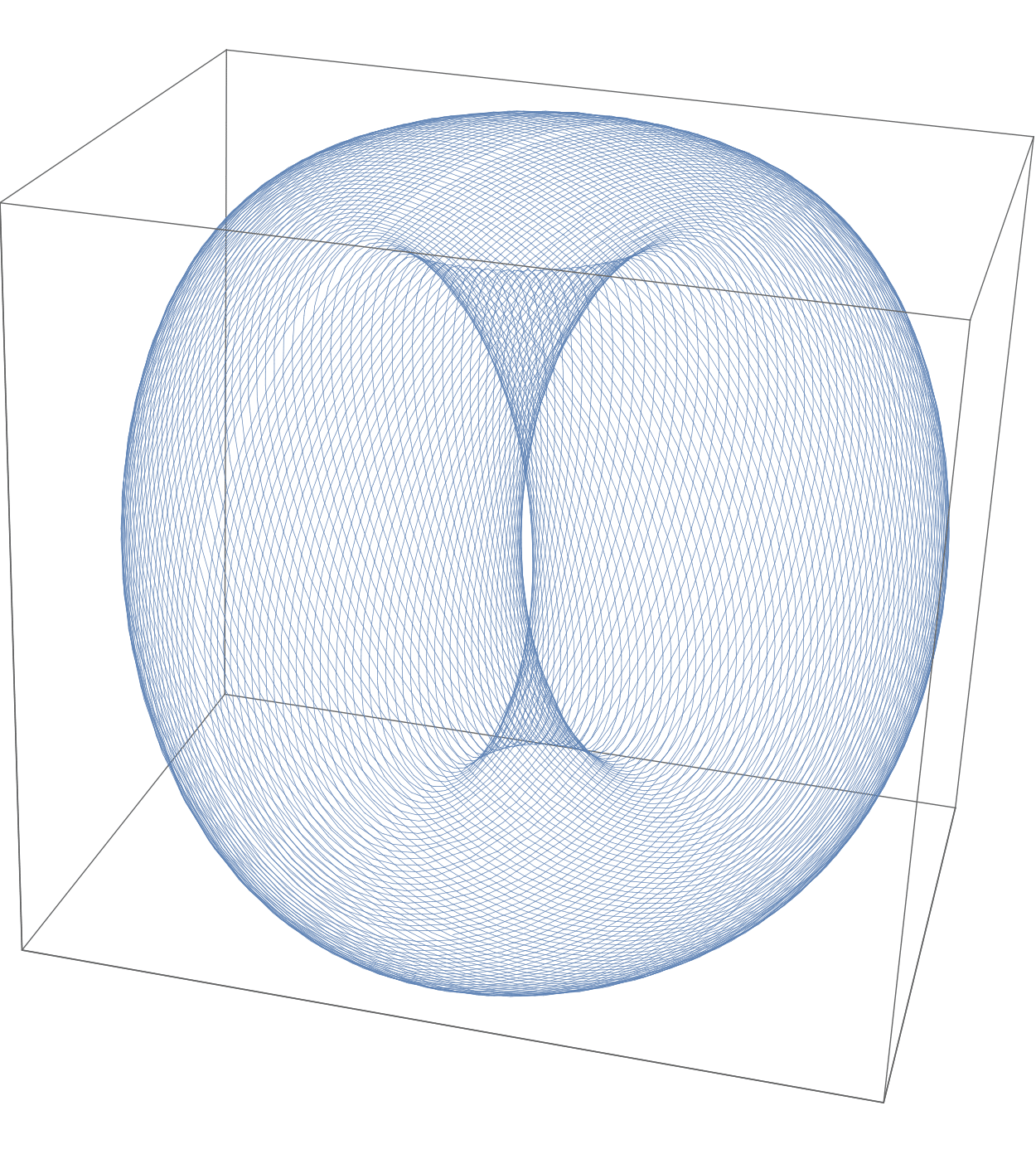} \\
 \includegraphics[width=0.31\textwidth]{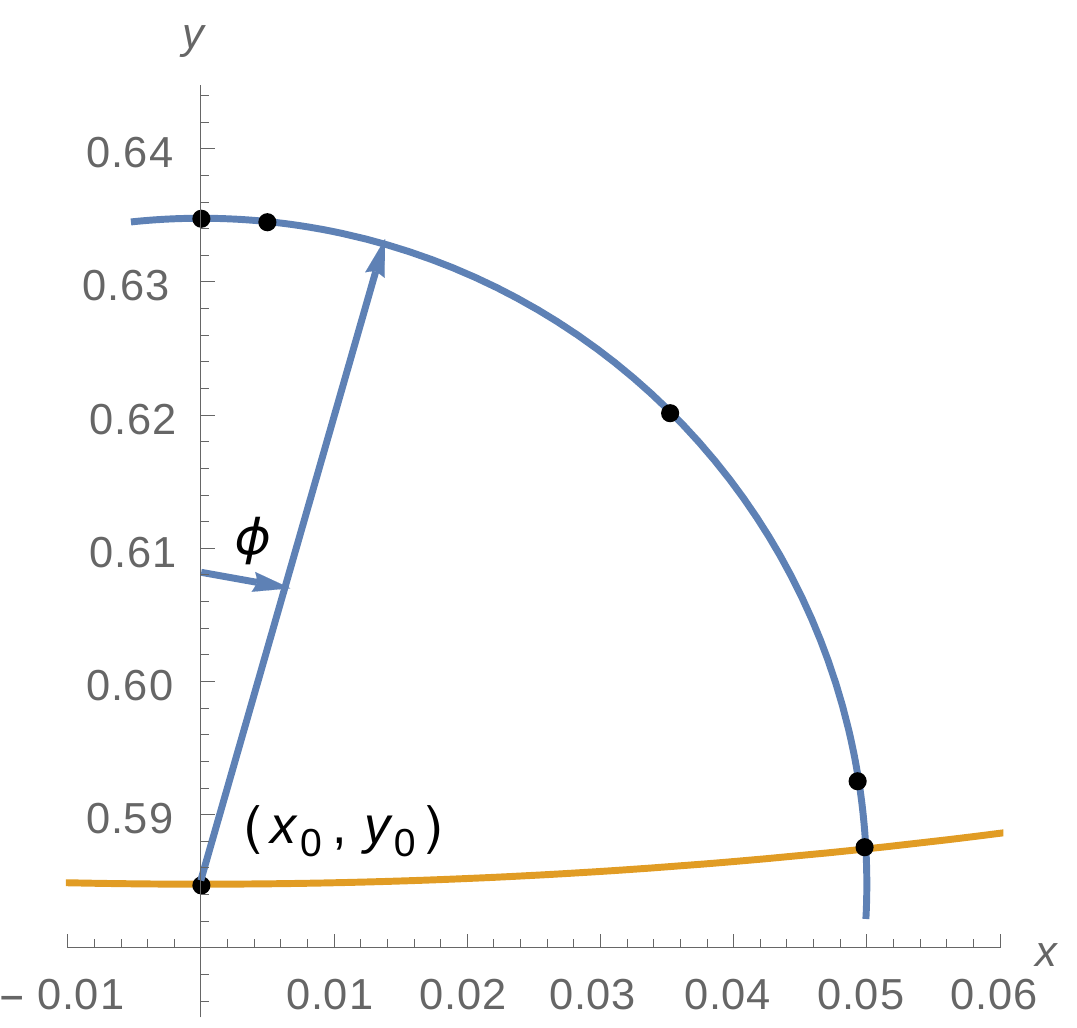}
 \includegraphics[width=0.31\textwidth]{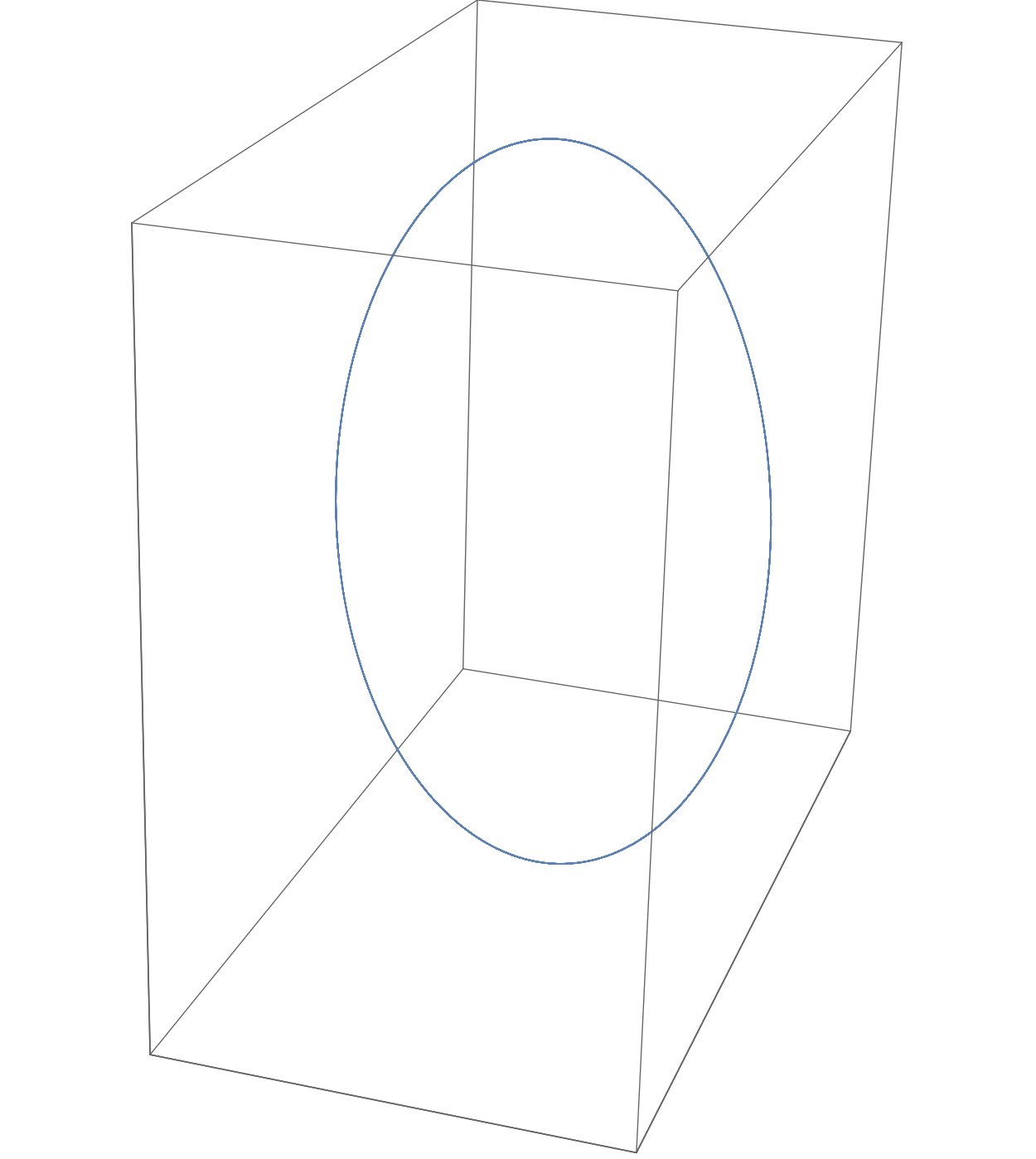}
 \includegraphics[width=0.31\textwidth]{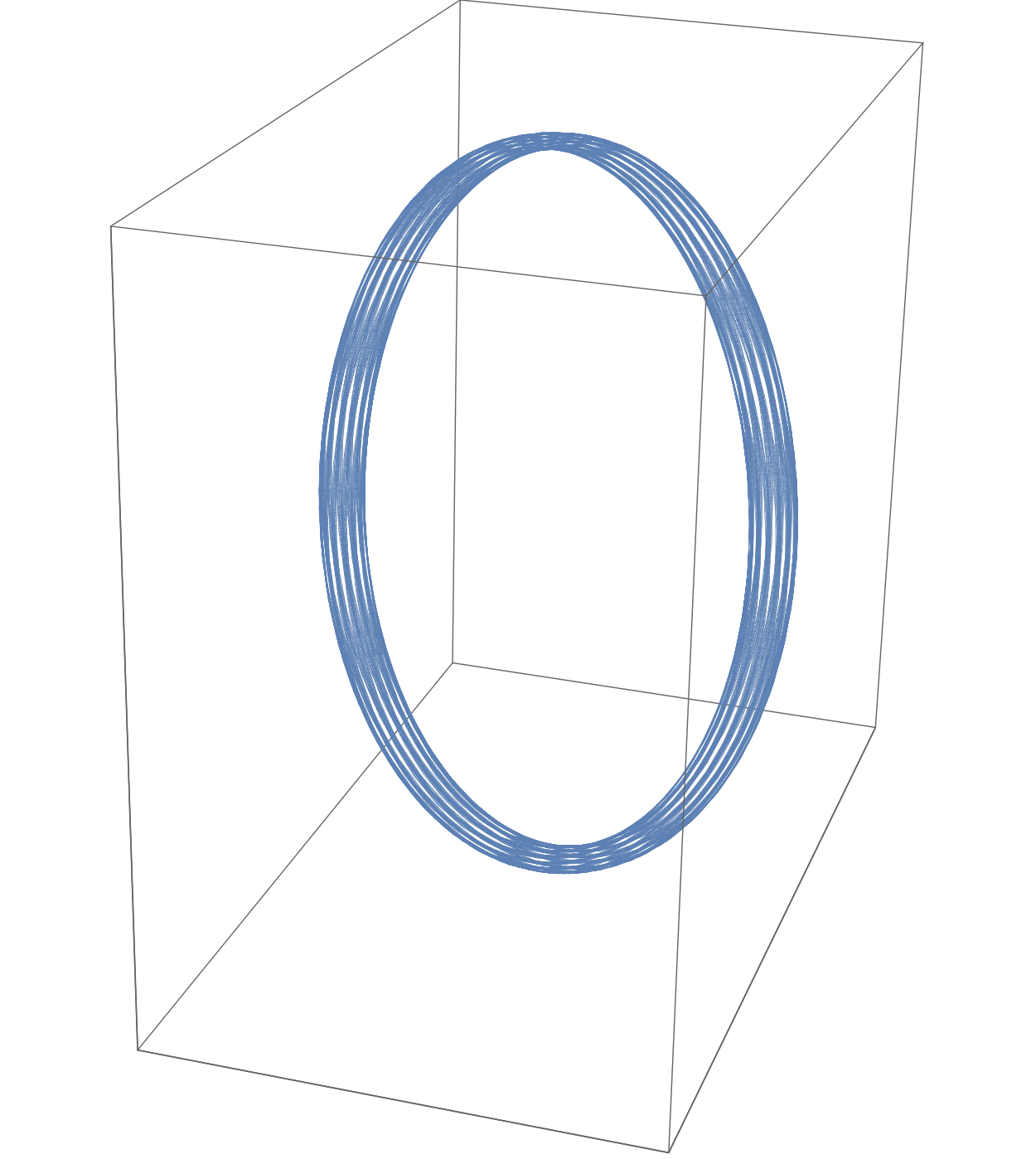}
 \caption{\textit{Lower left: Map of initial conditions shown for torus trajectories for $\alpha_0=2\pi/5$. The yellow curve is the graph of $y=g(x,0)$ in the invariant manifold. 
 Clockwise from upper left: Bouncing mode $\phi=0$; $\phi=0.1$;  $\phi=\pi/4$; $\phi=1.418$; rocking mode $\phi=1.518$}}
 \label{torusfig-3}
\end{figure}

Figures \ref{torusfig-1}, where $\alpha_0 = \pi/4$, and \ref{torusfig-3}, where $\alpha_0 = 2\pi/5$, show typical sets of trajectories. In all of these, the initial condition is perturbed from $(x_0,y_0)$ by a distance of 0.05. The figures show a range of angles of the perturbation, from vertical (angle $\phi=0$), lying in the bouncing manifold, to lying in the $y=g(x,w)$
rocking manifold).

\begin{figure}
 \centering
 \includegraphics[width=0.31\textwidth]{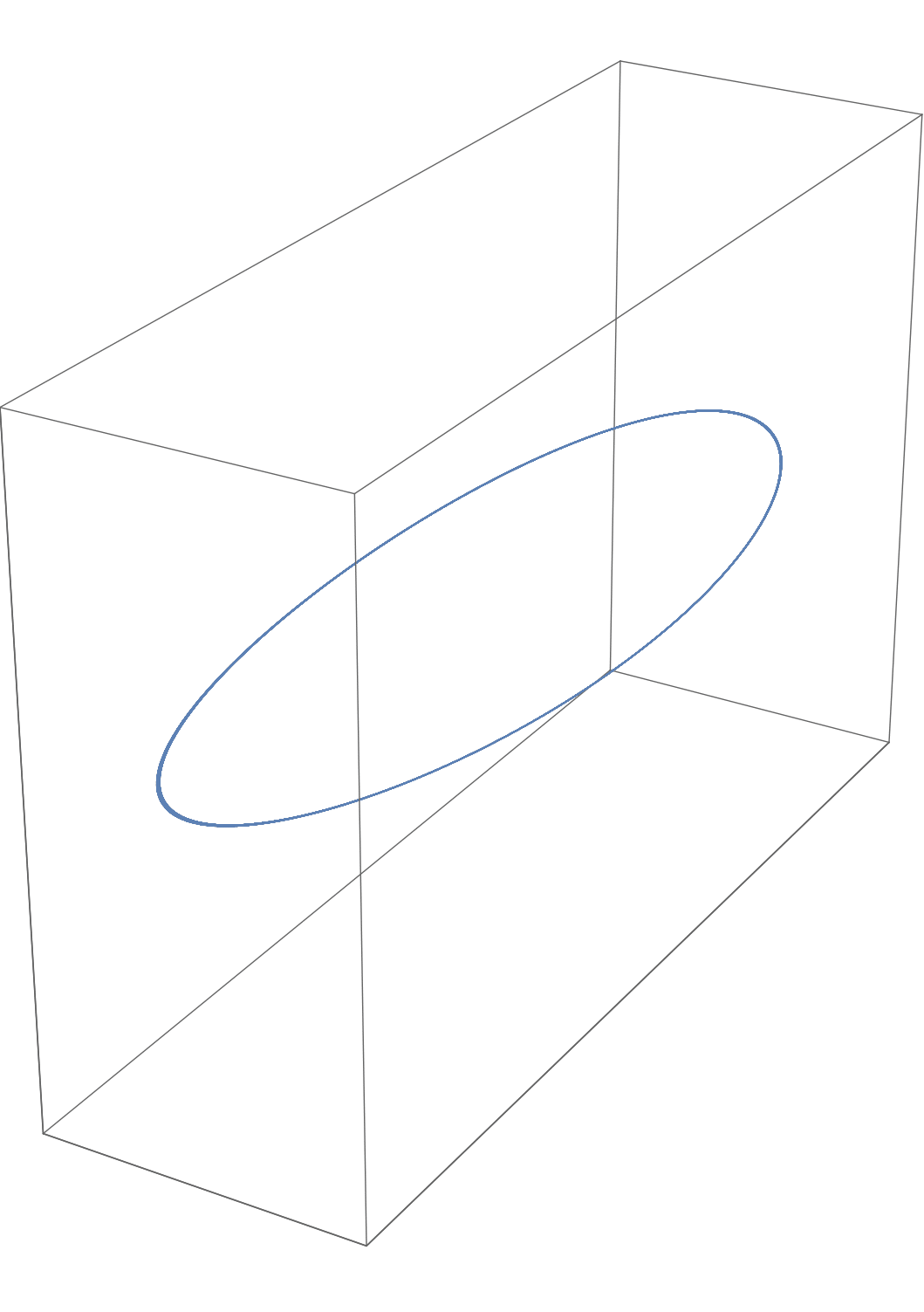}
 \includegraphics[width=0.31\textwidth]{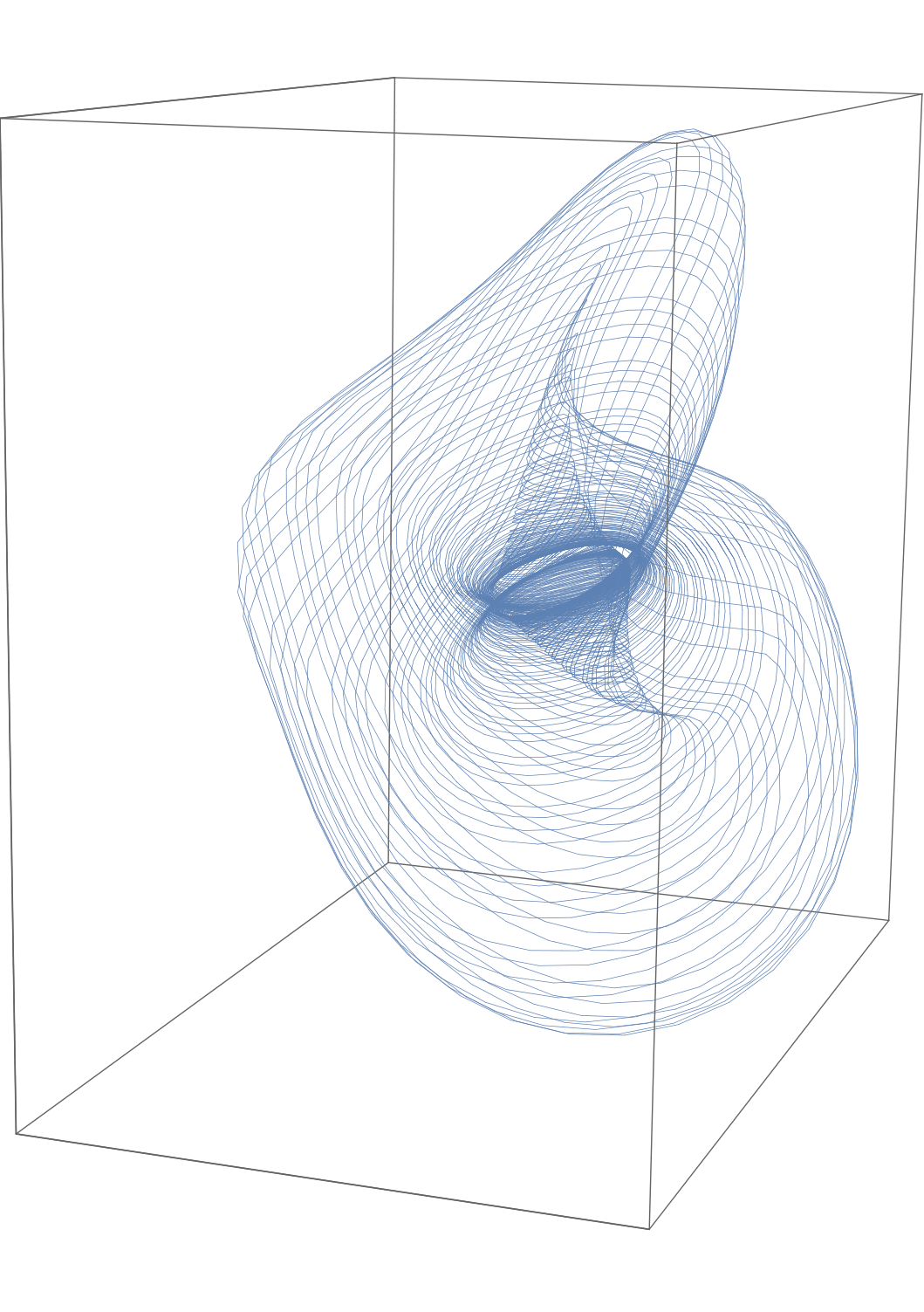}
 \includegraphics[width=0.31\textwidth]{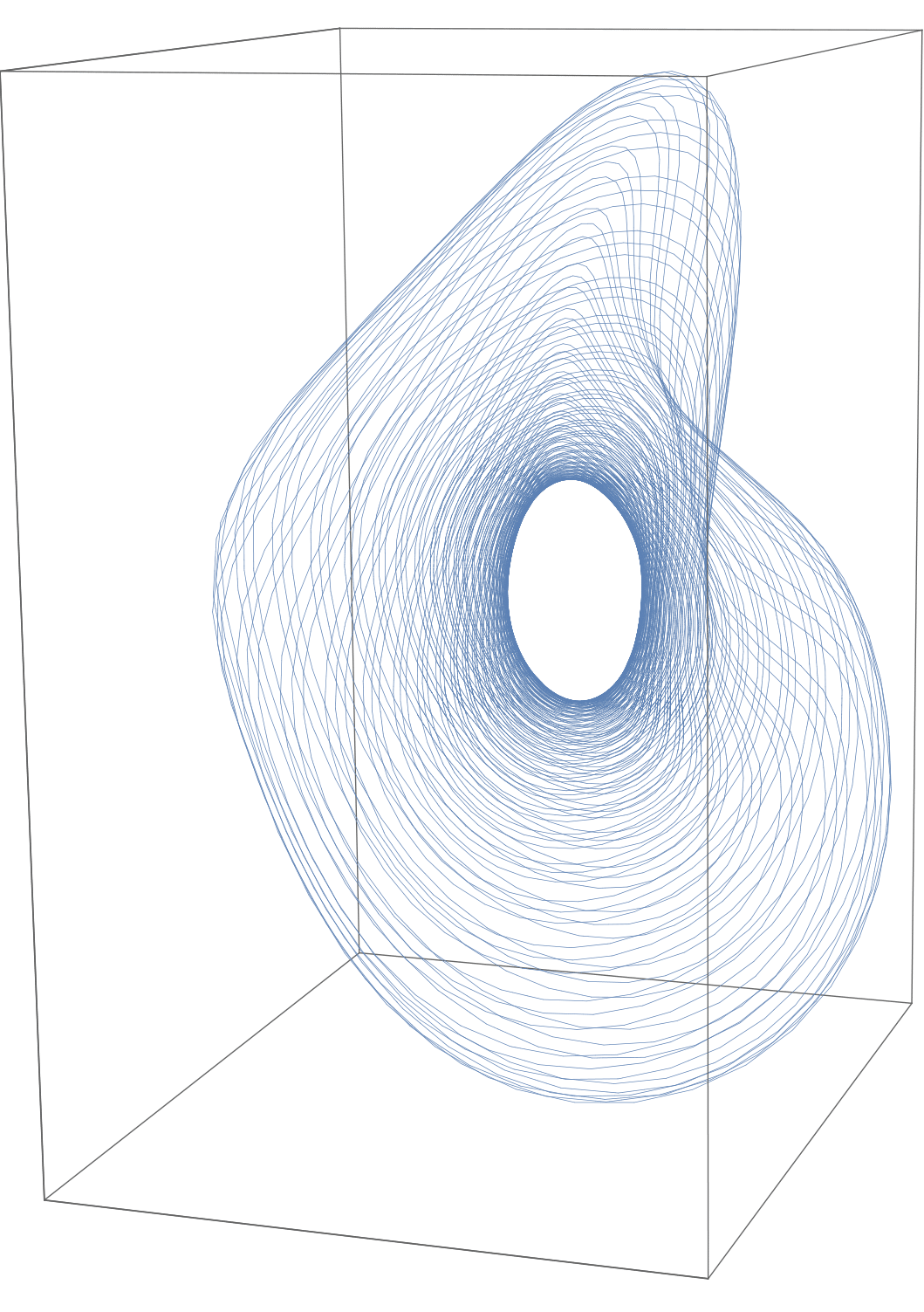} \\
 \includegraphics[width=0.38\textwidth]{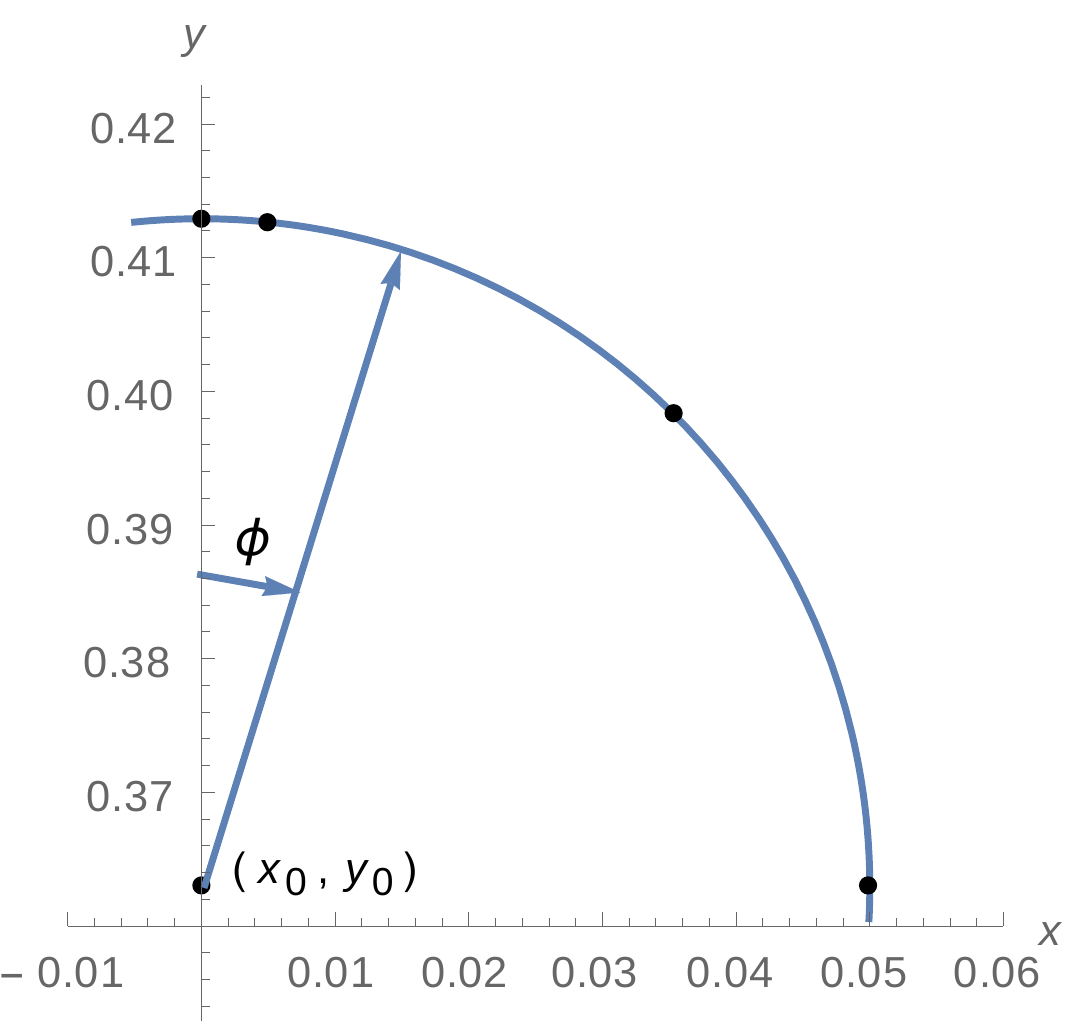}
 \includegraphics[width=0.31\textwidth]{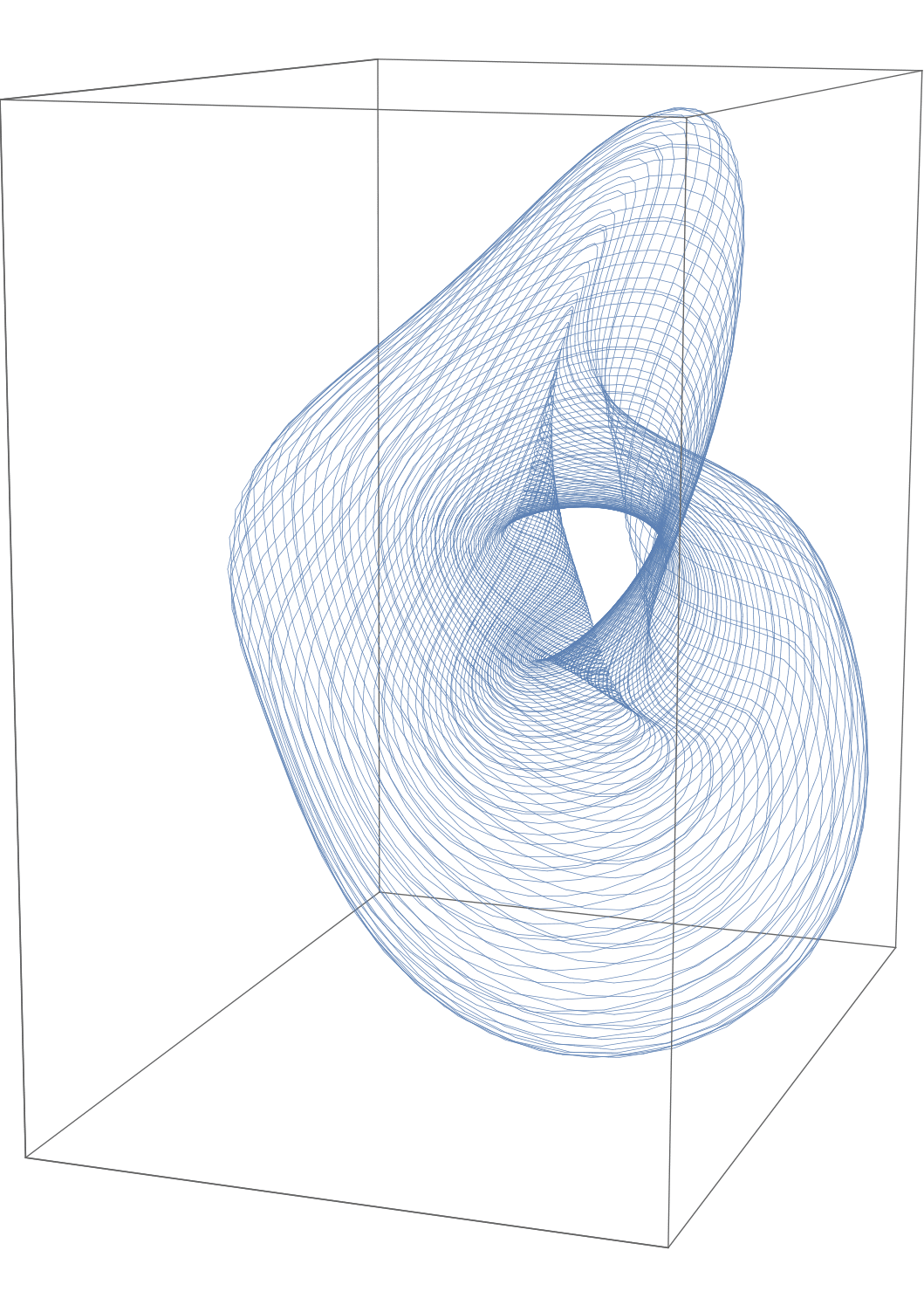}
 \caption{\textit{Lower left: Map of initial conditions shown for torus trajectories for $\alpha_0=\alpha_0^\dag$. For this value of $\alpha_0$ there is no differentiable rocking mode manifold. 
 Clockwise from upper left: Bouncing mode $\phi=0$; $\phi=0.1$; $\phi=\pi/4$; $\phi=\pi/2$}}
 \label{torusfig-2}
\end{figure}

Figure \ref{torusfig-2} shows a set of trajectories where $\alpha_0=\alpha_0^\dag\approx0.870\approx 49.8^\circ$ is the critical value at which the invariant manifold 
has a singularity. Again, the initial condition is perturbed from $(x_0,y_0)$ by a distance of 0.05. We see that no periodic orbit corresponding to $y=g(x,w)$ 
can be found.

\begin{figure}
 \centering
 \includegraphics[width=0.31\textwidth]{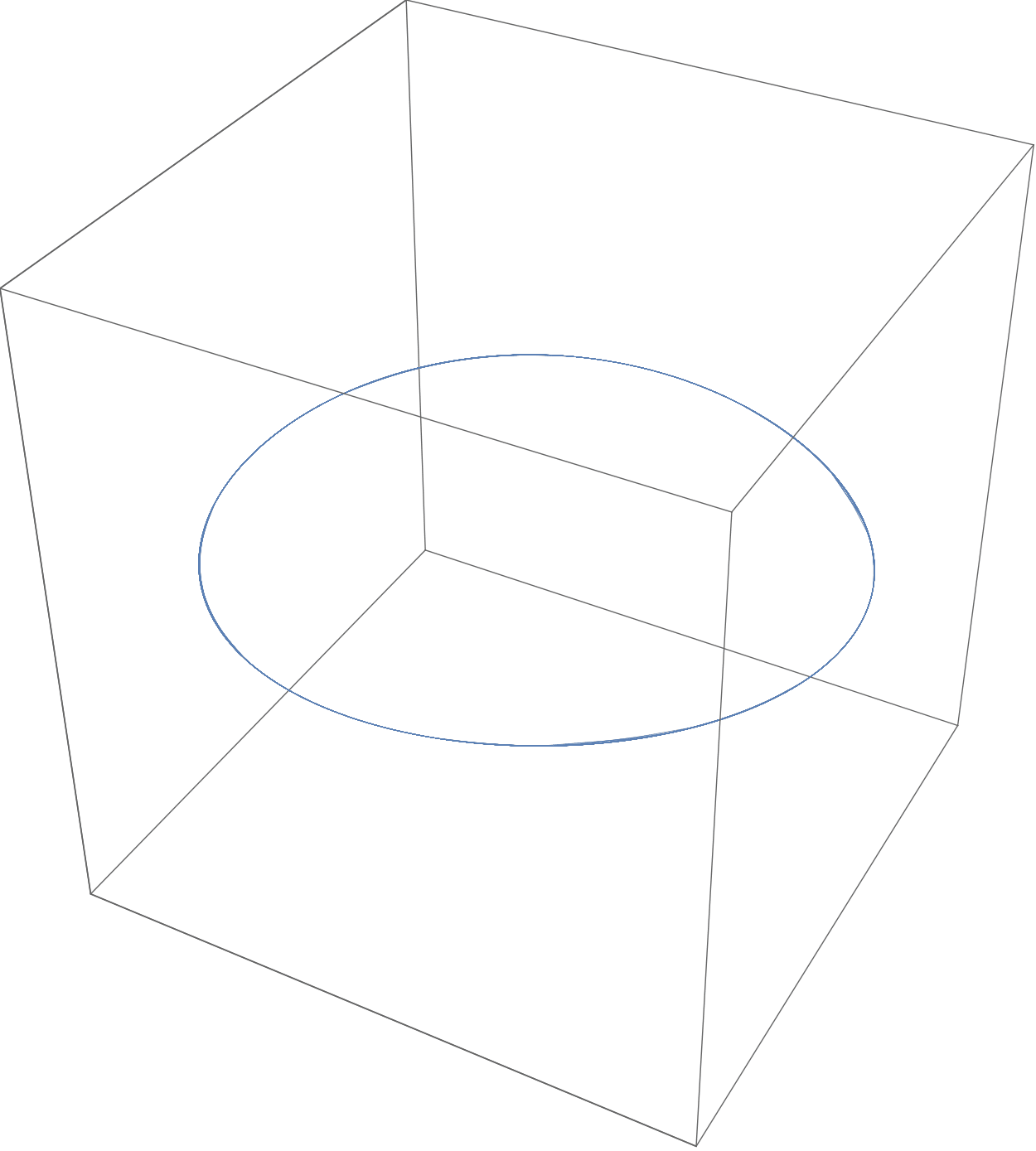}
 \includegraphics[width=0.31\textwidth]{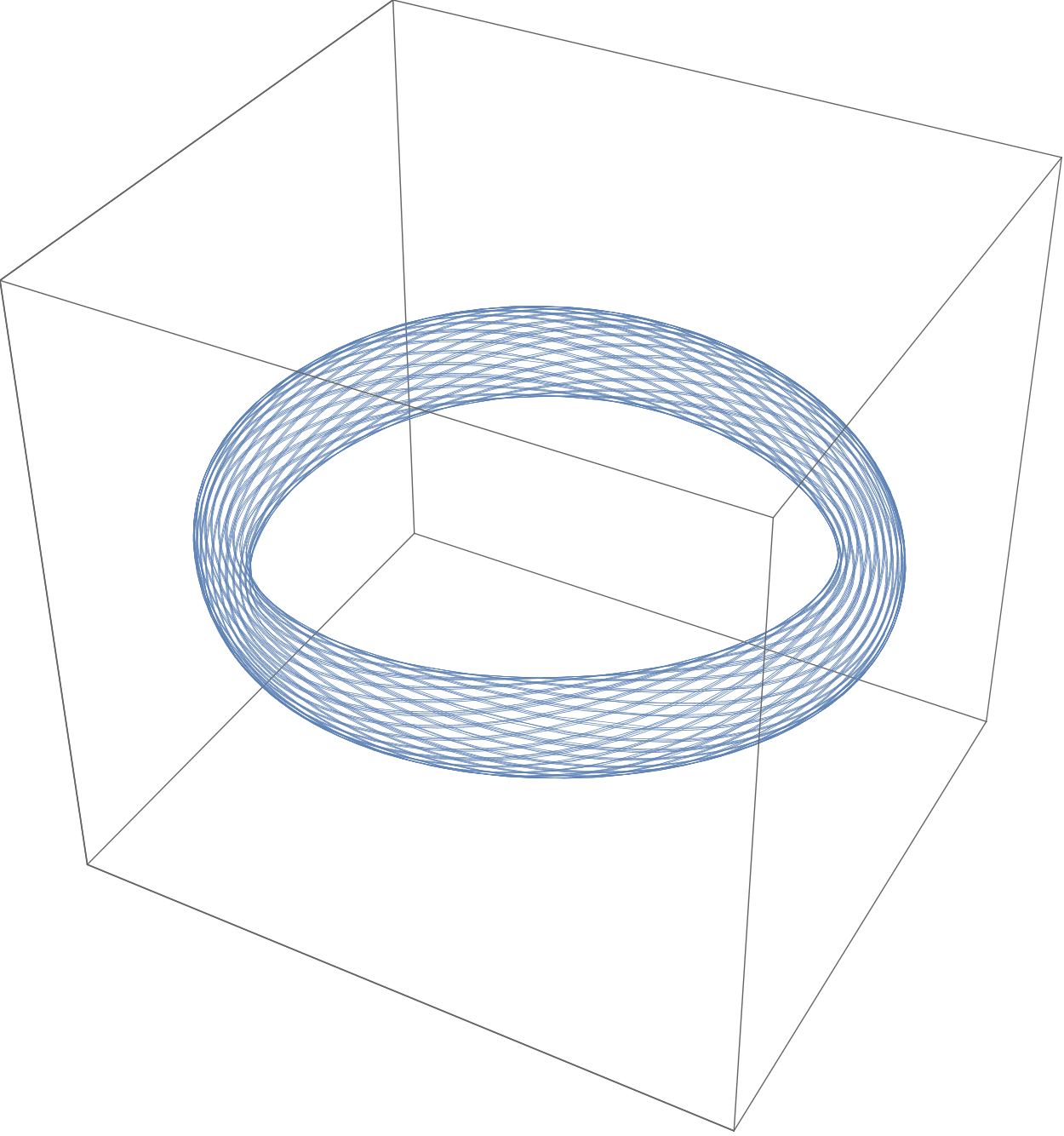}
 \includegraphics[width=0.31\textwidth]{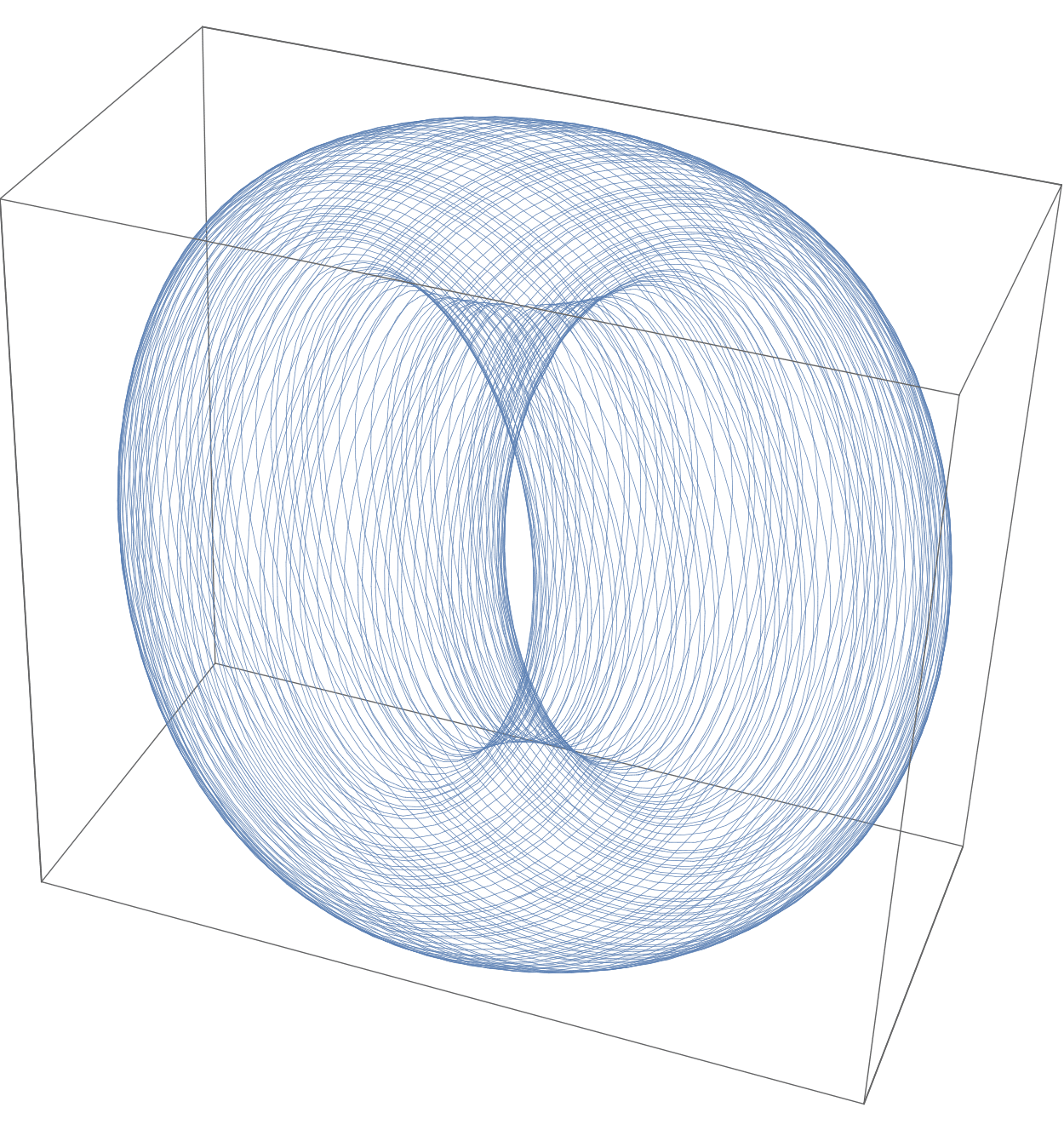} \\
  \includegraphics[width=0.31\textwidth]{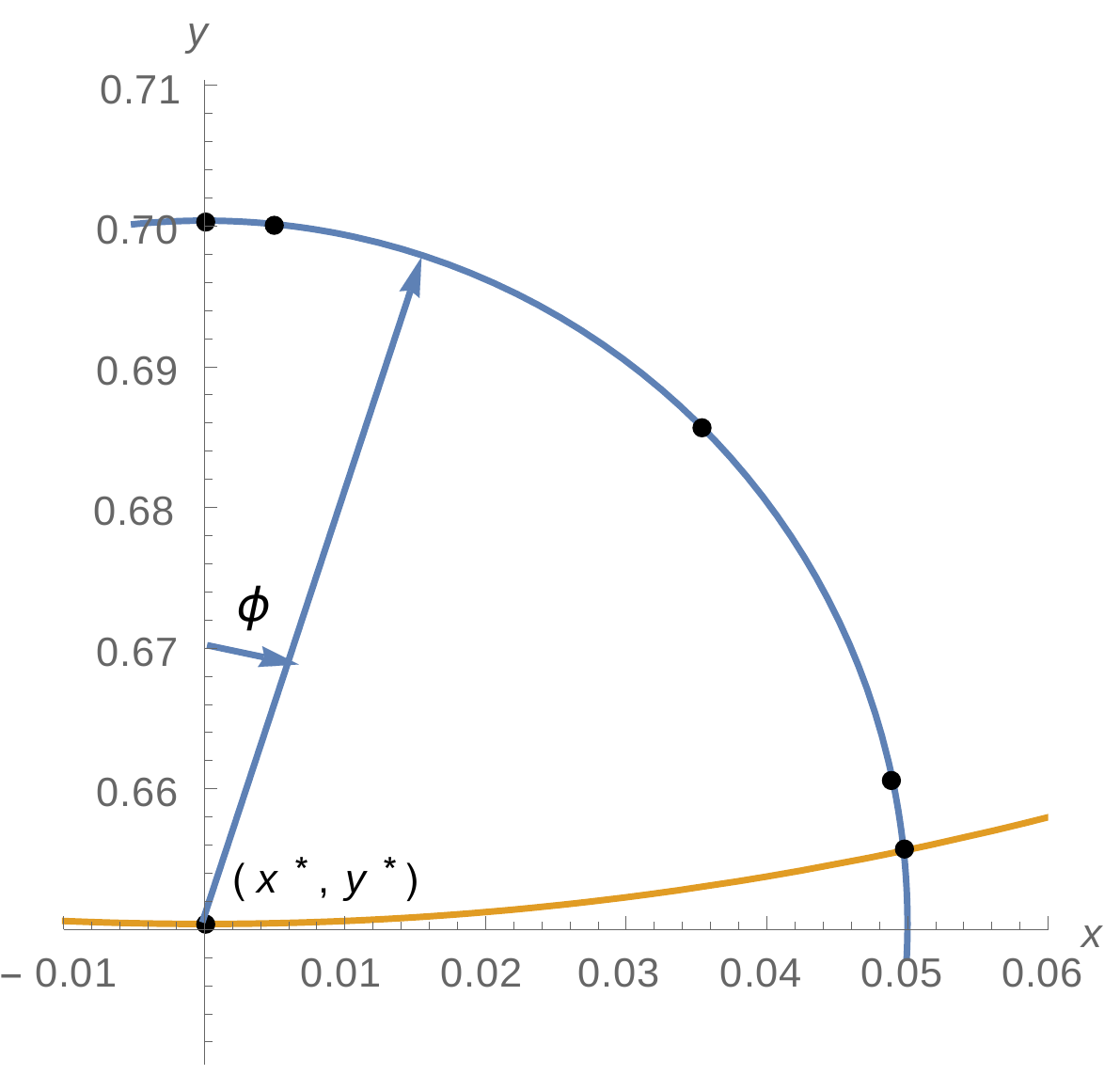}
  \includegraphics[width=0.31\textwidth]{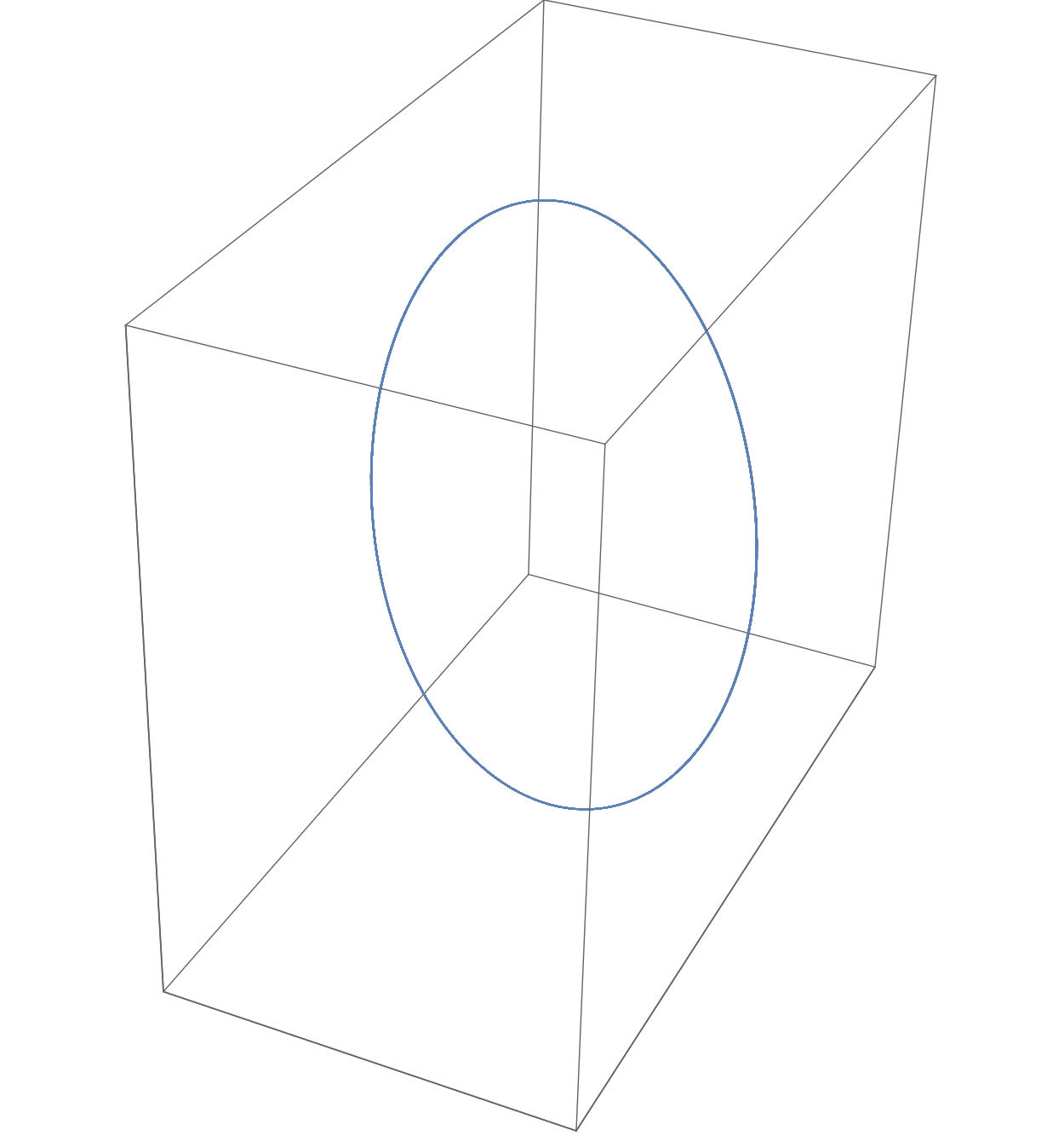}
  \includegraphics[width=0.31\textwidth]{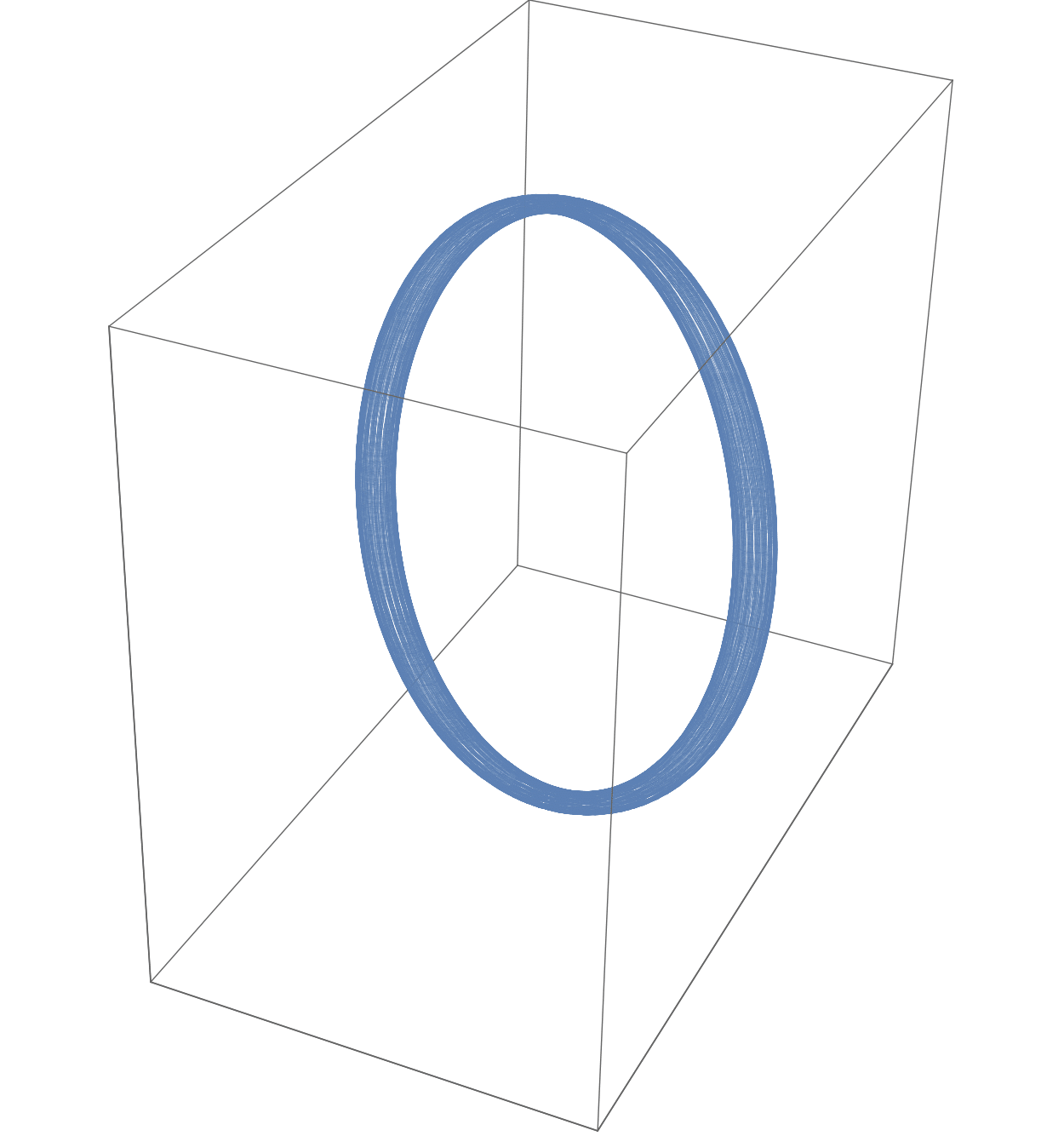}
 \caption{\textit{Lower left: Map of initial conditions shown for torus trajectories for $\alpha_0=1.45>\alpha_0^*$. The yellow curve is the graph of $y=g(x,0)$ in the invariant manifold. 
 Clockwise from upper left: Bouncing mode $\phi=0$; $\phi=0.1$; $\phi=\pi/4$; $\phi=1.366$; rocking mode $\phi=1.466$}}
\label{torusfig-4}
\end{figure}

Figure \ref{torusfig-4} shows a set of trajectories where $\alpha_0 = 1.45 > \alpha_0^*$. In this case the stable equilibrium
is $\mathbf{x}_1$, for which we have no closed-form expression. Thus the invariant manifold must be computed numerically.
Again, the initial condition is perturbed from $(x^*,y^*)$ by a distance of 0.05.

\section{Comparison to continuum spherical-cap drop model predictions}\label{sec:fluid}

Continuum models consider deformations with infinite degrees of freedom. Studies of these models typically solve, in some limit, the Navier-Stokes equations with  boundary conditions that accommodate moving contact lines \cite{sui2014numerical}. Bostwick and Steen \cite{bostwick} perform a linear stability analysis of a 2-parameter family of spherical-cap base states to find the resonant mode shapes and frequencies in the Euler equation limit.   These continuum predictions have guided  Chang et al.\cite{chang2013substrate} to the laboratory discovery of the first 35 modes. The space-time symmetries represented by these various solutions are non-trivial \cite{chang2015dynamics, steen2019droplet}. In this section we post process the mode shapes presented by Bostwick and Steen \cite{bostwick} to find their center-of-mass dynamics in order to compare with the bouncing and rocking modes reported above.

\begin{figure}
    \centering
    \includegraphics[width=0.45\textwidth]{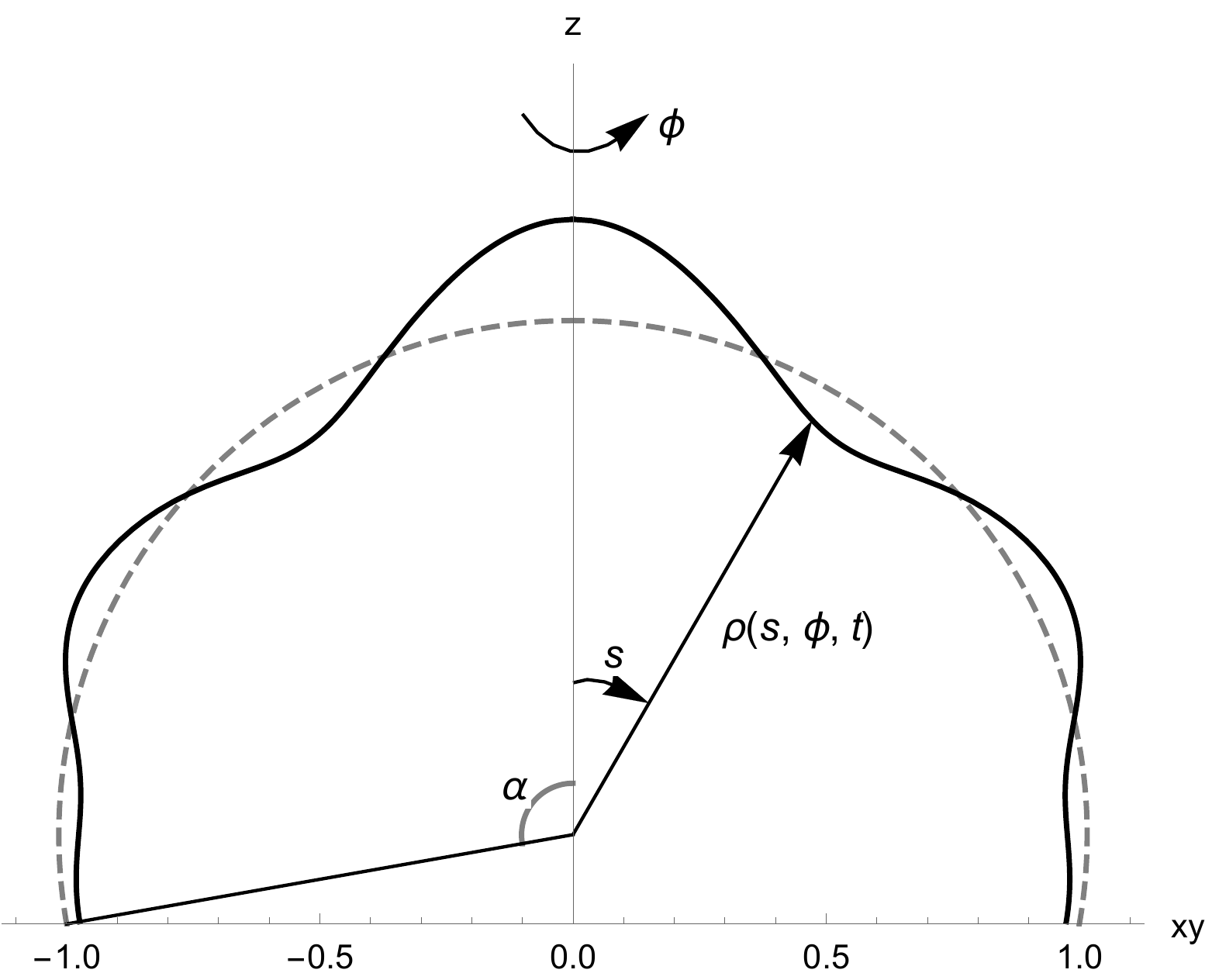}
    \caption{\textit{Definition sketch of liquid sessile drop in cross-section.}}
    \label{fig:jbb-drop}
\end{figure}

A static spherical-cap drop with equilbrium contact angle $\alpha$ is scaled so that the contact line has radius 1; thus the radius of the unperturbed drop is $\csc(\alpha)$. Normal modes for perturbations of the surface are parameterized by azimuthal wavenumber $l$ and frequency $\Omega$, which depends on the polar wavenumber $k$:
\begin{equation}
    \rho(s,\phi,t)=\csc(\alpha)+\epsilon \xi(s)\cos(l\phi)\cos(\Omega t)
\end{equation}
where $s$ is the polar angle. See Figure \ref{fig:jbb-drop}.
This ansatz holds for all $(k,l)$ except $(1,1)$, in which case it is inappropriate since there is instability \cite{bostwick}. To compute the normal modes, Bostwick and Steen reduce the dynamic pressure balance to an operator equation on a function space satisfying Laplace's equation and a no-penetration boundary condition. The normal modes are the eigenfunctions of this operator equation, given as sums of the harmonic basis functions\cite{bostwick}.

To calculate the center of mass of the perturbed drop, we integrate over two regions, taking the origin to be the center of the unperturbed drop: the cone whose base is the wetted disk and apex is the origin, and the sector with $0\leq s\leq \alpha$.

Note that if $\alpha<\pi/2$ the cone volume will be negative. Also, no domain perturbation is necessary for $\epsilon\ll 1$, the errors introduced will be higher order, that is, $O(\epsilon^2)$.

We compute
\begin{equation}
    (\bar x,\bar y, \bar z) = \frac{1}{M}\iiint_{drop} (x,y,z)dV
\end{equation}
where $M$ is the volume of the drop
\begin{align}
    M &= \iiint_{cone} dV + \iiint_{sector} dV \nonumber \\
      &= -\frac{\pi}{3}\cot(\alpha) + \int_0^{2\pi}\int_0^\alpha \int_0^{\rho(s,\phi,t)} r^2 \sin(s)\, dr\, ds\, d\phi 
\end{align}
and the integrals of $x$, $y$ and $z$ over the drop are similarly taken in spherical coordinates.

Evaluating the inner integral, we find that there are terms whose denominators vanish at $l=0$ and $l=1$. We treat these cases separately. As the integrand is a bounded function of $l$, we use the dominated convergence theorem to justify taking limits inside the integral using L'H\^{o}pital's rule.

\subsection{L=0 modes}
\begin{figure}
    \centering
    \includegraphics[width=0.35\textwidth]{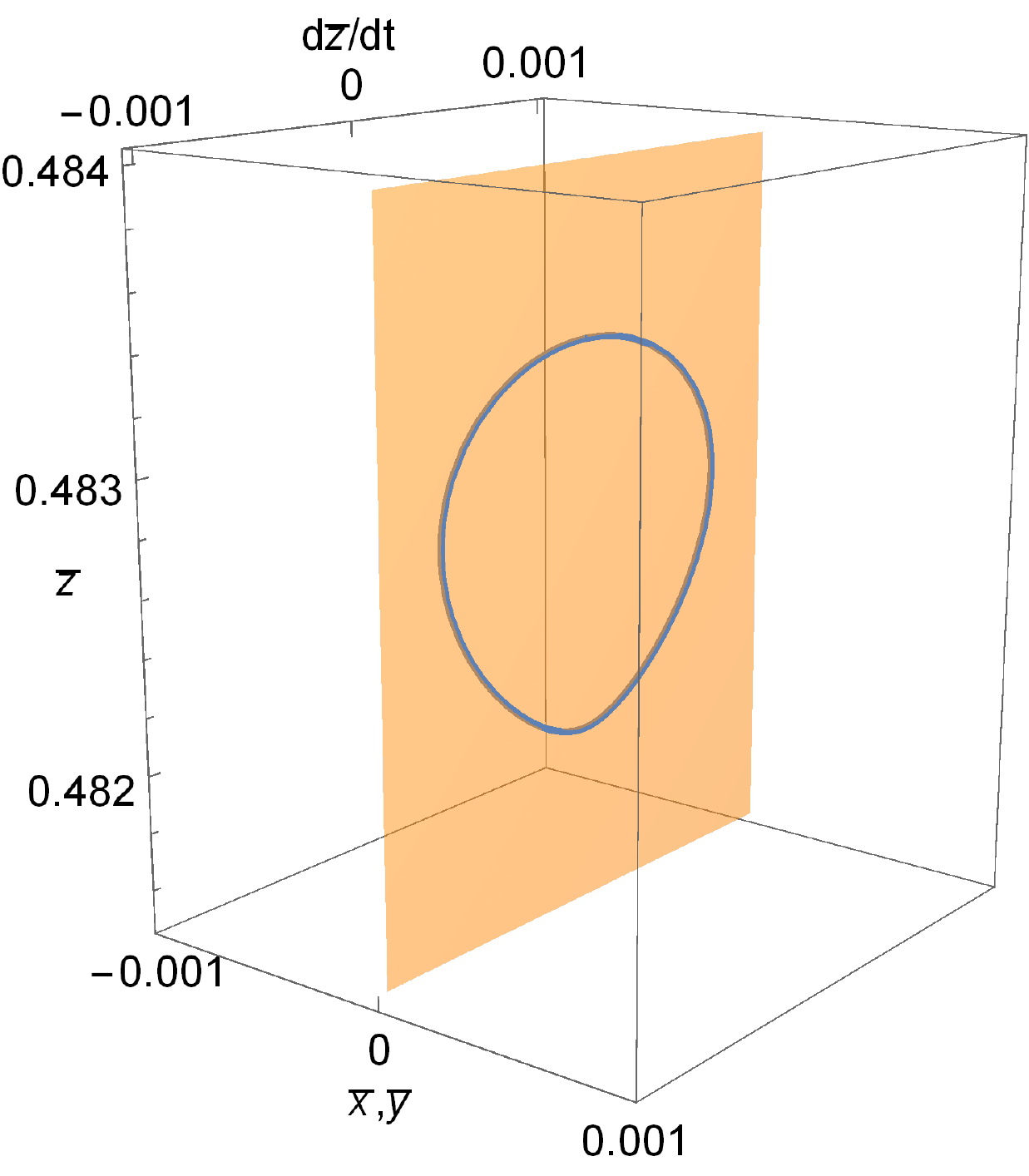}
    \caption{\textit{Parametric plot of center-of-mass trajectory of a sessile fluid drop with $\alpha=4\pi/9$, $l=0$, polar wavenumber $k=6$, $\epsilon=0.01$, natural contact line boundary conditions. It is pictured in the plane of all such trajectories for varying $\epsilon$.}}
    \label{fig:jbb-bounce}
\end{figure}

We find that for $l=0$
\begin{align}
    M &= \frac{\pi}{6}(\cos(\alpha )+2) \tan \left(\frac{\alpha }{2}\right) \sec ^2\left(\frac{\alpha }{2}\right) \nonumber \\
    & + \epsilon  \cos (t \Omega ) \int_0^{\alpha } 2 \pi  \csc ^2(\alpha ) \xi(s) \sin (s) \, ds + O(\epsilon^2) 
\end{align}
but the $O(\epsilon)$ term integrates to 0 when $\xi(s)$ is one of the eigenfunctions computed by Bostwick and Steen\cite{bostwick}. The unscaled coordinates of the center of mass are 
\begin{align}
    M\bar x &= 0 \\
    M\bar y &= 0 \\
    M\bar z &= \frac{\pi }{4}+\epsilon  \cos (t \Omega ) \int_0^{\alpha } \pi  \csc ^3(\alpha ) \xi(s) \sin (2 s) \, ds+O(\epsilon^2)
\end{align}
which tells us that the center-of-mass motion of the $l=0$ modes is vertical, corresponding to the bouncing mode of the triangular drop. These integrals can be computed numerically, taking $\xi(s)$ to be an eigenfunction. For example, see Figure \ref{fig:jbb-bounce}, and compare to Figure \ref{fig:bouncing} (a).

\subsection{L=1 modes}
\begin{figure}
    \centering
    \includegraphics[width=0.4\textwidth]{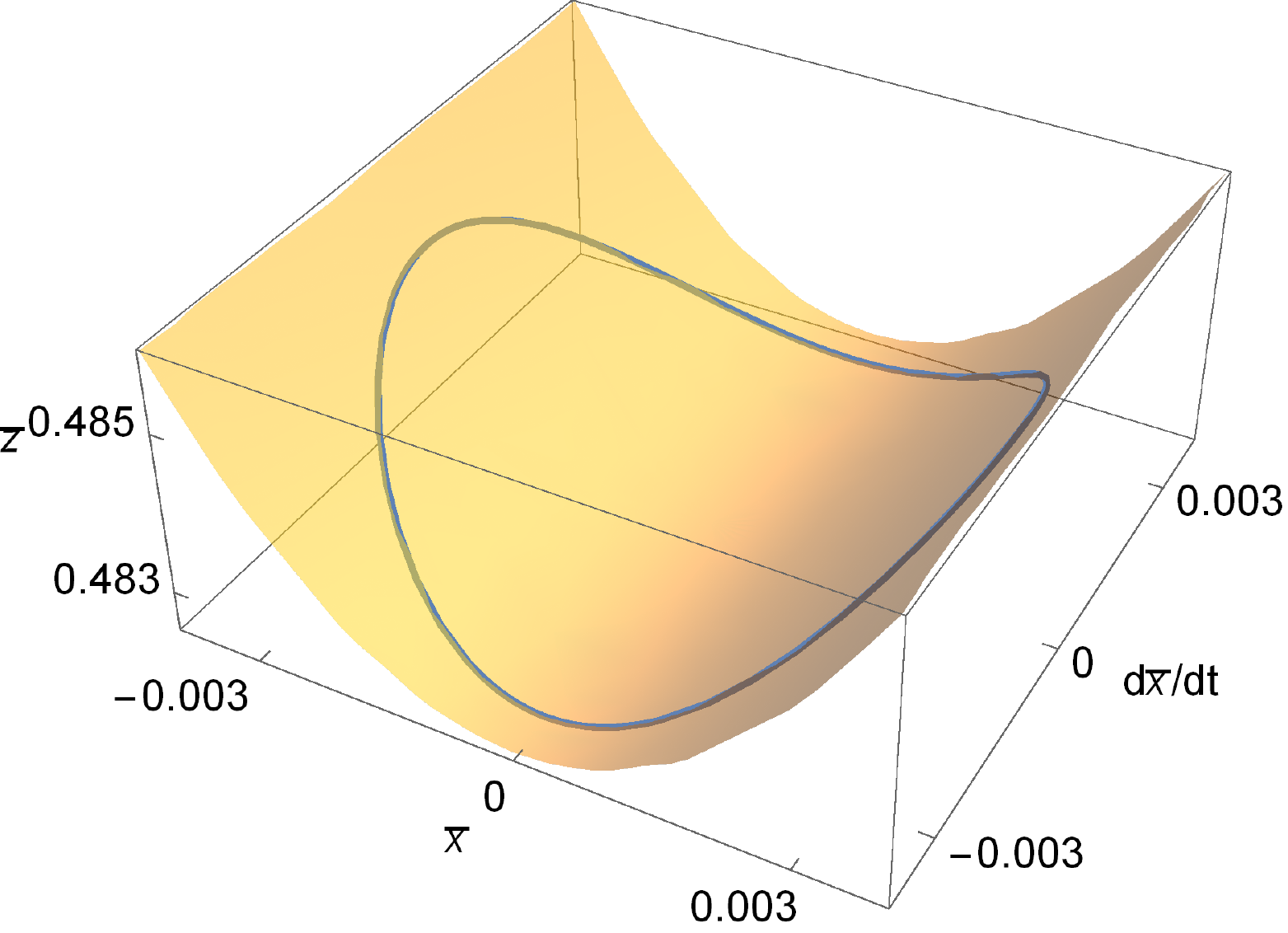}
    \caption{\textit{Parametric plot of center-of-mass trajectory of a sessile fluid drop with $\alpha=4\pi/9$, $l=1$, polar wavenumber $k=5$, $\epsilon=0.01$, pinned contact line boundary conditions. It is pictured in the manifold of all such trajectories for varying $\epsilon$.}}
    \label{fig:sessile-rocking}
\end{figure}

For $l=1$ (and $k\neq 1$) the volume and unscaled center-of-mass coordinates are
\begin{align}
    M &= \frac{\pi}{6}(\cos(\alpha )+2) \tan \left(\frac{\alpha }{2}\right) \sec ^2\left(\frac{\alpha }{2}\right) + O(\epsilon^2) \\
    M\bar x &= \epsilon  \cos (t \Omega ) \int_0^{\alpha } \pi  \csc ^3(\alpha ) \xi(s) \sin ^2(s) \, ds + O(\epsilon^2) \\
    M\bar y &= 0 \\
    M\bar z &= \frac{\pi }{4} + O(\epsilon^2)
\end{align}
which shows that the center-of-mass motion of the $l=1$ modes is mostly horizontal, and confined to a vertical plane determined by the phase of the $\phi$ term. We can draw a qualitative correspondence to the rocking mode of the triangular drop. A numerically integrated trajectory, taking $\xi(s)$ to be an eigenfunction, is shown in Figure \ref{fig:sessile-rocking}. Compare to Figure \ref{fig:bouncing} (b).

\subsection{L$>$1 modes}
For all other whole-number values of $l$ there are no $O(\epsilon)$ terms in the center-of-mass motion:
\begin{align}
    M &= \frac{\pi}{6}(\cos(\alpha )+2) \tan \left(\frac{\alpha }{2}\right) \sec ^2\left(\frac{\alpha }{2}\right) + O(\epsilon^2) \\
    M\bar x &= 0 \\
    M\bar y &= 0 \\
    M\bar z &= \frac{\pi }{4} + O(\epsilon^2)
\end{align}
Thus for $l>1$ the center of mass of the sessile drop is stationary to $O(\epsilon)$.

We conclude that the triangular drop model in fact qualitatively predicts all of the $O(\epsilon)$ pure-mode motions of the fluid drop's center of mass. In general, the motion of a fluid drop is a linear combination of modes, each with its own frequency. A combination of one $l=0$ and one $l=1$ mode may be identified with the torus trajectories of the triangular drop. 

\section{Conclusion}
We have introduced the sliding Steiner triangular drop in the spirit of a minimal model.  We are unaware of any prior studies of triangular drop models.  The dynamical system is surprisingly rich, organized by two 1-parameter families of fixed points, with the stable equilibrium at the intersection of two invariant manifolds.  As the parameter varies, the fixed points exchange stabilities at a critical value where the ``rocking" manifold folds up into the ``bouncing" manifold. Evidence presented pre- and post-collision suggests a complicated transition, whose exact nature remains an open question. Off the manifolds, but in a nested structure, exist 2-tori where quasi-periodic dynamics occur.  

To our surprise, the model captures space-time symmetry-breaking motions observed in real drops and may anticipate nonlinear behavior like torus dynamics and droplet jumping, currently under active study for real drops \cite{vahabi2018coalescence}.

The dynamical system  resides in a four dimensional phase space.  Although not Hamiltonian, the system is invariant under two time-reversing phase-space involutions and by these spatial symmetries the system falls under a theorem of Sevryuk that predicts 2 families of nested $m$-tori up to $2m=4$, where $2m$ is the dimension of the full phase space and $m$ is the dimension of the subspace fixed by each involution. That is, for our $4D$ phase space, nested 0-tori (fixed points), 1-tori (periodic orbits), and 2-tori are predicted. These are identified. An open question is whether chaotic orbits also exist: to our knowledge, the reversible KAM theory does not rule them out\cite{SEVRYUK1998} but we have not observed them.

The dynamical system depends on a single parameter, the rest state contact angle $\alpha_0$.  For every stable fixed point in this  1-parameter family of isosceles triangles, there is an unstable (saddle) fixed point. In a bifurcation diagram, these two families of fixed points cross transcritically at $\alpha_0 = \alpha_0^*$. We identify two 2-dimensional invariant manifolds $\{x=w=0\}$ and $\{y=g(x,w),z=d g(x,w)/dt\}$ where the stable fixed points are surrounded by periodic orbits. These manifolds are flat and nearly flat and constitute the bouncing and rocking modes of the drop, respectively.  A region around these periodic orbits in phase space is filled with nested invariant tori, which correspond to quasiperiodic orbits.  We observe several of these numerically.

The rocking manifold is apparently not smooth for all values of $\alpha_0$ as higher-order approximations reveal singular points where the power-series approximation fails. The nature of this manifold near these points remains an open question.


Here we remind the reader of the spirit of our minimal model.  Most engineering studies seek models that can make detailed predictions of observables.  We are not seeking such a model, rather a model with the fewest degrees of freedom that can predict qualitative features of the motions of droplets with moving contact lines.  Our Steiner triangular drop predicts bouncing and rocking motions.  Our model is minimal in the sense that, if any one degree of freedom is removed, our Steiner drop no longer functions as a drop.  For example, pinning either basal vertex confines the apex to a hyperbola while pinning the apex fixes the base length owing to the area constraint. In both cases, the dynamics become trivial. Like removing one leg of a three-legged stool, removing one degree of freedom from our model leaves it non-functional.

As a minimal model, adding features goes against its minimalist nature.  If one chooses nonetheless to do so, extensions in any number of directions are possible.   Contact line damping or bulk viscous dissipation could be added.  Extensions to three dimensions, say, to polyhedral `drops', would add deformational degrees of freedom. Tetrahedral drops might be related to the Steiner ellipsoid, if that aspect of the model is of interest.  The pressure constitutive relationship we introduce is purposely {\it ad hoc} and one could postulate other such relationships.  To study homogeneous (same $\alpha_0$) populations of interacting Steiner drops, one might add coalescence, jumping or ejection behavior depending on the phenomenon of interest.  Alternatively, for heterogeneous (different $\alpha_0$) populations, disparities in size and/or chemical make-up might be included.

\section*{Supplementary material}

See the supplementary material for videos of Steiner triangular drops: bouncing, rocking, torus trajectory, and unbounded trajectory.

\section*{Acknowledgments}
ENW acknowledges support from Cornell University as the inaugural CAM multidisciplinary postdoctoral fellow.  PHS acknowledges support by NSF Grant CBET-1637960.


\section*{Appendix: Torus trajectory figures}

\begin{thebibliography}{99}

\bibitem{basaran}
Basaran, Osman A. and David W. DePaoli. ``Nonlinear oscillations of pendant drops'', \textit{Physics of Fluids}, 6.9 (1994), pp. 2923--2943.

\bibitem{bostwick}
Bostwick, Joshua and Steen, Paul. ``Dynamics of sessile drops. Part 1. Inviscid theory'',
\textit{Journal of Fluid Mechanics} 760 (2014), pp. 5--38.

\bibitem{bostwick2015stability}
Bostwick, JB and Steen, PH. ``Stability of constrained capillary surfaces'',
\textit{Annual Review of Fluid Mechanics} 47 (2015), pp. 539--568.

\bibitem{chang2013substrate}
Chang, Chun-Ti and Bostwick, Joshua B and Steen, Paul H and Daniel, Susan. 
``Substrate constraint modifies the Rayleigh spectrum of vibrating sessile drops'',
\textit{Physical Review E} 88.2, 023015 (2013).

\bibitem{chang2015dynamics}
Chang, Chun-Ti and Bostwick, JB and Daniel, Susan and Steen, PH.
``Dynamics of sessile drops. Part 2. Experiment'',
\textit{Journal of Fluid Mechanics} 768 (2015), pp. 442--467.

\bibitem{davis}
Davis, Stephen H. ``Moving contact lines and rivulet instabilities. Part 1. The static rivulet'', \textit{Journal of Fluid Mechanics}, 98.2 (1980), pp. 225--242.  
  
\bibitem{deGennes}
De Gennes, Pierre-Gilles. ``Wetting: statics and dynamics'', \textit{Reviews of modern physics}, 57.3 (1985), p. 827. 
  
\bibitem{devaney},
Devaney, Robert L. ``Reversible Diffeomorphisms and Flows'',
\textit{Transactions of the American Mathematical Society} 218 (1976), pp. 89--113.

\bibitem{enright}
Enright, Ryan et al. ``How coalescing droplets jump'', \textit{ACS nano} 8.10 (2014), pp. 10352--10362.

\bibitem{fayzra}
Fayzrakhmanova, Irina S. and Arthur V. Straube. ``Stick-slip dynamics of an oscillated sessile drop'', \textit{Physics of fluids}, 21.7 (2009), pp. 072104. 
 
\bibitem{fowlkes2011self}
Fowlkes, Jason D and Kondic, Lou and Diez, Javier and Wu, Yueying and Rack, Philip D.
``Self-assembly versus directed assembly of nanoparticles via pulsed
laser induced dewetting of patterned metal films'',
\textit{Nano letters} 11.6 (2011), pp. 2478--2485.

\bibitem{josserand2016drop}
Josserand, Christophe and Thoroddsen, Sigurdur T. ``Drop impact on a solid surface'',
\textit{Annual review of fluid mechanics} 48 (2016), pp. 365-391.

\bibitem{kumar2015liquid}
Kumar, Satish. ``Liquid transfer in printing processes: liquid bridges with moving contact lines'',
\textit{Annual Review of Fluid Mechanics} 47 (2015), pp. 67--94.
  
\bibitem{lamb}
Jeroen S.W. Lamb and John A.G. Roberts. "Time-reversal symmetry in dynamical systems: A survey",
\textit{Physica D: Nonlinear Phenomena} 112.1 (1998), pp. 1--39.
  
\bibitem{liu}
Liu, Fangjie et al. ``Numerical simulations of self-propelled jumping upon drop coalescence on non-wetting surfaces'', \textit{Journal of Fluid Mechanics} 752 (2014), pp. 39--65.  
  
\bibitem{lyubimovNonAxi}
Lyubimov, D. V. and T. P. Lyubimova and S. V. Shklyaev. ``Non-axisymmetric oscillations of a hemispherical drop'', \textit{Fluid Dynamics}, 39.6 (2004), pp. 851--862.  

\bibitem{lyubimovBehavior}
Lyubimov, Dmitry V. and Tatyana P. Lyubimova and Sergey V. Shklyaev. ``Behavior of a drop on an oscillating solid plate'', \textit{Physics of fluids}, 18.1 (2006), pp. 012101.
  
\bibitem{macner2014condensation}
Macner, Ashley M and Daniel, Susan and Steen, Paul H.
``Condensation on surface energy gradient shifts drop size distribution toward small drops'',
\textit{Langmuir} 30.7 (2014), pp. 1788--1798.  

\bibitem{noblin}
Noblin, X., A. Buguin, and F. Brochard-Wyart. ``Vibrated sessile drops: Transition between pinned and mobile contact line oscillations'', \textit{The European Physical Journal E}, 14.4 (2004), pp. 395--404.
  
\bibitem{rayleigh18}
Rayleigh, Lord. ``On the capillary phenomena of jets'', \textit{Proc. R. Soc. London}, 29.196-199 (1879), pp. 71-97.  
  
\bibitem{ROBERTS1992}
J.A.G. Roberts and G.R.W. Quispel. 
"Chaos and time-reversal symmetry. Order and chaos in reversible dynamical systems",
\textit{Physics Reports} 216.2 (1992), pp. 63--177.

\bibitem{sevryuk91}
M.B. Sevryuk. ``Lower‐dimensional tori in reversible systems'',
\textit{Chaos: An Interdisciplinary Journal of Nonlinear Science} 1.2 (1991), pp. 160--167.

\bibitem{SEVRYUK1998}
M.B. Sevryuk. "The finite-dimensional reversible KAM theory",
\textit{Physica D: Nonlinear Phenomena} 112.1 (1998), pp. 132--147.

\bibitem{sharp}
Sharp, James S. ``Resonant properties of sessile droplets; contact angle dependence of the resonant frequency and width in glycerol/water mixtures'', \textit{Soft Matter}, 8.2 (2012), pp. 399--407.

\bibitem{sharpEtAl}
Sharp, James S. and David J. Farmer and James Kelly, ``Contact angle dependence of the resonant frequency of sessile water droplets'', \textit{Langmuir}, 27.15 (2011), pp. 9367--9371.

\bibitem{snoeijer2013moving}
Snoeijer, Jacco H and Andreotti, Bruno. ``Moving contact lines: scales, regimes, and dynamical transitions'',
\textit{Annual review of fluid mechanics} 45 (2013), pp. 269--292.

\bibitem{steen2019droplet}
Steen, Paul H and Chang, Chun-Ti and Bostwick, Joshua B.
``Droplet motions fill a periodic table'',
\textit{Proceedings of the National Academy of Sciences} 116.11 (2019), pp. 4849--4854.

\bibitem{steiner1829}
Steiner, J. ``Géométrie pure. Développement d'une série de théorèmes relatifs aux sections coniques'',
\textit{Annales de Mathématiques pures et appliquées} 19 (1828-1829), pp. 37--64.

\bibitem{sui2014numerical}
Sui, Yi and Ding, Hang and Spelt, Peter DM. ``Numerical simulations of flows with moving contact lines'',
\textit{Annual Review of Fluid Mechanics} 46 (2014).

\bibitem{tanaka}
Tanaka, H. ``Measurements of Drop-Size Distributions during Transient Dropwise Condensation'',
\textit{J. Heat Transfer} 97.3 (1975), pp. 341--346.

\bibitem{vahabi2018coalescence}
Vahabi, Hamed and Wang, Wei and Mabry, Joseph M and Kota, Arun K.
``Coalescence-induced jumping of droplets on superomniphobic surfaces with macrotexture'',
\textit{Science advances} 4.11, eaau3488 (2018).

\bibitem{weisstein}
Weisstein, Eric W. ``Steiner Circumellipse'', from MathWorld -- A Wolfram Web Resource.
http://mathworld.wolfram.com/SteinerCircumellipse.html

\bibitem{xia}
Xia, Yi, and Paul H. Steen. ``Moving contact-line mobility measured'', \textit{Journal of Fluid Mechanics}, 841 (2018), pp. 767--783.

\bibitem{xu2017collective}
Xu, Chenglong and Yu, Haitao and Peng, Shuhua and Lu, Ziyang and
Lei, Lei and Lohse, Detlef and Zhang, Xuehua. 
``Collective interactions in the nucleation and growth of surface droplets'',
\textit{Soft matter} 13.5 (2017), pp. 937--944.

\end{thebibliography}
\end{document}